\documentclass[a4paper,12pt,fleqn]{article}
\usepackage{libertine}
\usepackage{amsmath}
\usepackage{amsthm}
\usepackage{mathrsfs}
\usepackage{amssymb}
\usepackage{framed}
\usepackage[inline]{enumitem}
\usepackage{mathtools}
\usepackage{mdframed}
\usepackage[all]{xy}
\usepackage[dvipsnames]{xcolor}
\usepackage{color}
\usepackage[utf8]{inputenc}
\usepackage{setspace}
\usepackage{currfile}
\usepackage{comment}
\usepackage[unicode,colorlinks,linkcolor=blue]{hyperref}
\usepackage{bm}
\usepackage{tensor}
\usepackage{thmtools}
\usepackage{xfrac}

\setlist[enumerate]{topsep=0pt,itemsep=-1ex,partopsep=1ex,parsep=1ex}

\usepackage[backend=biber,style=alphabetic,sorting=nyt,giveninits=true,language=english,maxbibnames=10]{biblatex}

\renewbibmacro{in:}{%
  \ifentrytype{article}{}{\printtext{\bibstring{in}\intitlepunct}}}

\DeclareSourcemap{
  \maps[datatype=bibtex, overwrite]{
    \map{
      \step[fieldset=language, null]
       \step[fieldset=shorthand, null]
    }
  }
}

\bibliography{arxiv_v2.bib}

\newcommand{\fell}{\mathfrak{l}}

\DeclareMathOperator{\Fl}{Fl}
\DeclareMathOperator{\ext}{ext}
\DeclareMathOperator{\pr}{pr}
\DeclareMathOperator{\Sym}{Sym}

\DeclareMathOperator{\Riem}{Riem}
\DeclareMathOperator{\curv}{{\boldsymbol{\mathrm{R}}}}
\DeclareMathOperator{\supp}{supp}
\DeclareMathOperator{\GL}{GL}
\DeclareMathOperator{\ev}{ev}
\DeclareMathOperator{\carr}{carr}
\DeclareMathOperator{\id}{id}

\newcommand{\sym}{\mathrm{sym}}
\newcommand{\cP}{\mathcal{P}}
\newcommand{\bN}{\mathbb{N}}
\newcommand{\csub}{\subset\subset}
\newcommand{\ud}{\mathrm{d}}
\newcommand{\loc}{\mathrm{loc}}
\newcommand{\eqcol}{=\mathrel{\mathop:}}
\newcommand{\coleq}{\mathrel{\mathop:}=}
\newcommand{\catVec}{\bm{\mathrm{Vec}}}
\DeclareMathOperator{\cacx}{\overline{acx}\,}
\DeclareMathOperator{\Hom}{\mathrm{Hom}}
\DeclareMathOperator{\Lin}{\mathrm{L}}
\newcommand{\Lint}{\mathrm{L}}
\DeclareMathOperator{\Vol}{Vol}
\newcommand{\D}{\mathrm{D}}
\newcommand{\cA}{\mathcal{A}}

\newcommand{\fX}{\mathfrak{X}}

\newcommand{\bR}{\mathbb{R}}

\newcommand{\bE}{\mathbb{E}}
\newcommand{\bF}{\mathbb{F}}
\newcommand{\sF}{\mathscr{F}}

\newcommand{\cD}{\mathcal{D}}
\newcommand{\cE}{\mathcal{E}}

\newcommand{\Lie}{\mathrm{L}}

\newcommand{\cN}{\mathcal{N}}

\newcommand{\cG}{\mathcal{G}}
\newcommand{\cU}{\mathcal{U}}
\newcommand{\cT}{\mathcal{T}}

\newcommand{\VSO}{\mathrm{VSO}}
\newcommand{\SK}{\mathrm{SK}}
\newcommand{\SO}{\mathrm{SO}}
\newcommand{\TO}{\mathrm{TO}}
\newcommand{\tang}{\mathrm{T}}

\newcommand{\csn}{\mathrm{csn}}
\newcommand{\bvarphi}{{\boldsymbol\varphi}}
\newcommand{\bpsi}{{\boldsymbol\psi}}
\newcommand{\bPhi}{{\boldsymbol\Phi}}
\newcommand{\bPsi}{{\boldsymbol\Psi}}

\newcommand{\pd}{\partial}
\newcommand{\bTheta}{{\boldsymbol\Theta}}
\newcommand{\ba}{{\boldsymbol a}}

\newcommand{\bb}{{\boldsymbol b}}
\newcommand{\bA}{{\boldsymbol A}}

\newcommand{\bB}{{\boldsymbol B}}

\newcommand{\unterstrich}{{-}}

\newcommand{\abso}[1]{\left|#1\right|}
\newcommand{\norm}[1]{\lVert#1\rVert}

\declaretheoremstyle[headfont=\scshape]{mystyle}
\declaretheorem[numberwithin=section,style=mystyle,name=Theorem,refname={theorem,theorems},Refname={Theorem,Theorems}]{theorem}
\declaretheorem[name=Definition,style=mystyle,refname={definition,definitions},Refname={Definition,Definitions},sibling=theorem]{definition}
\declaretheorem[name=Proposition,style=mystyle,refname={proposition,propositions},Refname={Proposition,Propositions},sibling=theorem]{proposition}
\declaretheorem[name=Corollary,style=mystyle,refname={corollary,corollaries},Refname={Corollary,Corollaries},sibling=theorem]{corollary}
\declaretheorem[name=Lemma,style=mystyle,refname={lemma,lemmas},Refname={Lemma,Lemmas},sibling=theorem]{lemma}
\declaretheorem[name=Remark,refname={remark,remarks},Refname={Remark,Remarks},sibling=theorem]{remark}

\newcommand{\e}{\varepsilon}

\title{Spacetimes with distributional semi-Riemannian metrics and their curvature}

\author{E.~A.~Nigsch\footnote{Faculty of Mathematics, University of Vienna, Oskar-Morgenstern-Platz 1, 1090 Vienna, Austria. email: \href{mailto:eduard.nigsch@univie.ac.at}{eduard.nigsch@univie.ac.at}}}
\date{\today}

\begin{document}

\maketitle

\begin{abstract}
We develop a comprehensive geometric framework for defining spaces $\cG(M,E)$ of nonlinear generalized sections of vector bundles $E \to M$ containing spaces of distributional sections $\cD'(M, E)$. Our theory incorporates classical differential geometric operations (like tensor products, covariant derivatives and Lie derivatives), is localizable and fully compatible with smooth and distributional tensor calculus. As an application to the treatment of singular metrics, we calculate the curvature of the conical metric used to describe cosmic strings.
\end{abstract}

\section{Introduction}\label{section.1}
While the theory of distributions of L.~Schwartz has proved to be extremely powerful for a large variety of linear problems, its applicability in nonlinear situations is more restricted. A common situation where difficulties appear is given by singular Riemannian metrics occurring in the study of space-times describing concentrated sources like point particles or cosmic strings \cite{gerochtraschen}. The underlying reason is that the nonlinear operations involved in the calculation of the curvature tensor are not compatible with differentiation and the linear structure of distribution theory.

In this article we will build a theory that is able to consistently handle such problems on the base of regularization methods. As an application, we will focus on the conical metric
\begin{equation}\label{equation.1.1}
\ud s^2 = \ud r^2 + \alpha^2 r^2 \ud \phi^2 \qquad (\abso{\alpha} < 1)
\end{equation}
and calculate its curvature, which will be proportional to the Dirac delta distribution. In fact, it will turn  out that for any ``reasonable'' regularization of the conical metric, the scalar curvature of the regularization converges to $4 \pi(1-\alpha)\delta$. The main virtue of our framework lies in the fact that the result is valid independently of the coordinate system used and that it holds for a large class of possible regularizations; moreover, we have full compatibility with situations where the curvature makes sense classically.

There is considerable interest in spacetimes with metrics of of low regularity (see, e.g., \cite{DahiaRomero,JVdistgen,genrel,LeeLeFloch,Lott2016,2012CQGra,zbMATH06592071}). Approaching this problem by working with regularized version of these metrics poses two fundamental questions: first, do the obtained results depend on the particular choice of regularization? And second, to what extent is this approach compatible with metrics were the classical theory still applies? Taking into account Schwartz' result about the impossibility of multiplication of distributions \cite{Schwartz}, an optimal solution to both questions is given by Colombeau algebras of generalized functions \cite{ColNew, ColElem, MOBook, GKOS}. While on open subsets of $\bR^n$ the development of Colombeau algebras proceeded rapidly and brought numerous applications mainly in the field of non-linear partial differential equations, the development of a geometric theory on manifolds faced many technical and conceptual difficulties. It took a number of attempts to establish a diffeomorphism invariant scalar theory in the local setting: in \cite{AB}, a variant of the so-called special algebra on manifolds was given where, however, the embedding of the space of distributions is not canonical and where the Lie derivative does not commute with this embedding; in \cite{colmani} an approach to give a diffeomorphism invariant theory was given, which was shown to be incomplete in \cite{Jelinek}, where an improved version of the theory was suggested. Only in \cite{found1}, the question of existence of a diffeomorphism invariant Colombeau algebra on open subsets of $\bR^n$ was settled.

While it was possible to extend the diffeomorphism invariant local theory of \cite{found1} to manifolds in \cite{global}, it became apparent rather soon that this setting was too inflexible for an extension to a theory of vector-valued nonlinear generalized functions; the best one was able to obtain from it is the theory of generalized tensor fields in \cite{global2} where not only the sheaf property failed to hold, but where it also was impossible to define a reasonable notion of covariant derivative. This incited a number of fundamental investigations and changes of perspective, resulting in a completely new foundation of the theory \cite{papernew,specfull,vecreg,bigone}. The culmination of this development in the form of the construction of a global theory of generalized tensor fields which on the one hand has all desirable abstract properties (diffeomorphism invariance, Lie and covariant derivatives, tensor calculus, canonical embeddings and the sheaf property) and on the other hand is usable in concrete calculations is presented in this article.

The structure of this article is as follows: in Section \ref{section.2} we will introduce the notation and terminology used for the rest of the exposition and the necessary background on smooth and distributional differential geometry; in Sections \ref{section.3}--\ref{section.5} we detail the theoretical framework which will enable us to calculate the curvature of metrics such as \eqref{equation.1.1}; this calculation will finally be performed in Section \ref{section.6}.

\section{Preliminaries}\label{section.2}

As index set used throughout we fix $I \coleq (0,1]$; note that we could use $\bN$ instead, but in practice this makes no difference because $I$ has $\{1/n: n \in \bN\}$ as countable cofinal subset. We frequently use the Landau notation
\[ f(\e) = O(g(\e)) \qquad (\e \to 0) \]
signifying that
\[ \exists C,\e_0>0\ \forall \e<\e_0: \abso{f(\e)} \le C g(\e). \]

By $\cU_x(X)$ we denote the filter base of open neighborhoods of an element $x$ of a topological space $X$. For two subsets $A,B$ of a topological space we write $A \csub B$ if $A$ is compact and contained in the interior $B^\circ$ of $B$. Given a function $f \colon X \to \bR$, $\carr f \coleq \{ x \in X : f(x) \ne 0 \}$ is called the \emph{carrier} of $f$, and the closure $\supp f = \overline{\carr f}$ its \emph{support}.

The projection on the $i$th factor of a product $U_1 \times \dotsm \times U_n$ is denoted by $\pr_i$ or $\pr_{U_i}$. By $\id_X$ or simply $\id$ we denote the identity on a set $X$. The equivalence class of an object $x$ with respect to some equivalence relation or a quotient is denoted by $[x]$.

$B_r(x)$ denotes the open Euclidean ball of radius $r>0$ with center $x \in \bR^n$. The total differential of a function $f \in C^\infty(\bE, \bF)$ (where $\bE, \bF$ are finite-dimensional vector spaces) is denoted by $\D f \in C^\infty(\bE, \Lin(\bE, \bF))$.

For two $R$-modules $M_1$ and $M_2$, $\Hom_R(M_1, M_2)$ denotes the set of $R$-module homomorphisms from $M_1$ to $M_2$.

Locally convex spaces will always be assumed to be Hausdorff; we use \cite{zbMATH03230708,Schaefer} as standard references. For any locally convex space $\bE$ we denote by $\csn(\bE)$ the set of continuous seminorms on $\bE$. Given two locally convex spaces $\bE$, $\bF$, the space of continuous linear mappings from $\bE$ to $\bF$ is denoted by $\Lin(\bE, \bF)$. By default, this space is endowed with the topology of uniform convergence on the family of all bounded subsets of $\bE$. For calculus in infinite-dimensional locally convex spaces we employ \emph{convenient calculus} of \cite{KM}. In this context, the differential of a smooth mapping $f \in C^\infty(\bE, \bF)$ is denoted by $\ud f \in C^\infty(\bE, \Lin(\bE, \bF))$; the differential with respect to the $i$th variable of $f \in C^\infty(\bE_1 \times \dotsm \times \bE_n, \bF)$ is denoted by $d_i f \in C^\infty(\bE_1 \times \dotsm \times \bE_n, \Lin(\bE_i, \bF))$.

\subsection{Differential geometry}

In general, our terminology for differential geometric notions mostly follows that of \cite{Lee}, except for the fact that (as in \cite{marsden}) for vector bundles the typical fiber may be any (finite dimensional) vector space $\bE$ instead of only $\bR^n$ -- this makes working with tensor products more comfortable, as this way we do not have to choose a basis of $\bE$. A manifold will always mean a Hausdorff paracompact smooth manifold. We will only employ real, smooth vector bundles and smooth trivializations thereof. The projection of a given vector bundle $E$ will usually be denoted by $\pi$ or $\pi_E$. If $E \to M$ is a vector bundle with trivialization $(U,\Phi)$ we set $\Phi_x \coleq \Phi|_{E_x} \colon E_x \to \{x\} \times \bE \cong \bE$, where $E_x$ is the fiber of $E$ over $x$. $\Riem(M)$ is the set of Riemannian metrics on a manifold $M$. For $g \in \Riem(M)$, $B^g_r(x)$ denotes the open metric ball of radius $r>0$ centered at $x\in M$ with respect to $g$. By $\fX(M)$ we denote the space of smooth vector fields on $M$. $\tang M$ is the tangent bundle of $M$, and $\tang_x M$ the tangent space at $x \in M$.

In the following, suppose we are given manifolds $M,N$ as well as vector bundles $E \to M$ and $F \to N$ with typical fibers $\bE$ and $\bF$, respectively.

Let $G \colon E \to F$ be a (smooth) bundle morphism covering $g \colon M \to N$. Given charts $(U,\varphi)$ of $M$ and $(V,\psi)$ of $N$ such that $g(U) \subseteq V$ as well as local trivializations $\Phi \colon \pi^{-1}(U) \to U \times \bE$ and $\Psi \colon \psi^{-1}(V) \to V \times \bF$, the \emph{local expression} of $G$ is given by the mapping
\begin{equation}\label{equation.2.2}
\begin{aligned}
\varphi(U) \times \bE & \to \psi(V) \times \bF \\
(x,v) & \mapsto \left((\psi \times \id) \circ \Psi \circ G \circ \Phi^{-1} \circ (\varphi^{-1} \times \id)\right) (x, v) \\
& \qquad = (g_0(x), G_0(x)\cdot v)
\end{aligned}
\end{equation}
with $g_0 = \psi \circ g \circ \varphi^{-1} \in C^\infty(\varphi(U), \psi(V))$ and $G_0 \in C^\infty(\varphi(U), \Lin(\bE, \bF))$.

If $G$ is a vector bundle isomorphism with inverse $H$ we accordingly have $h_0(y) = g_0^{-1}(y)$ and $H_0(y) = (G_0(g^{-1}(y)))^{-1}$ for $y \in \psi(g(U))$.

The \emph{external tensor product} of $E$ and $F$ is given by
\[ E \boxtimes F \coleq \bigsqcup_{(x,y) \in M \times N} E_x \otimes F_y \]
with projection $\pi(v \otimes w) = (\pi_E(v), \pi_F(w))$ for $v \in E$ and $w \in F$. For each $(x,y) \in M \times N$ choose neighborhoods $V$ of $x$ in $M$ and $W$ of $y$ in $N$ such that there exist local trivializations $(V, \Phi)$ of $E$ and $(W, \Psi)$ of $F$. Define the trivialization $\Theta \colon \pi^{-1}(V \times W) \to V \times W \times (\bE \otimes \bF)$ by
\begin{equation}\label{equation.2.3}
\Theta(v \otimes w) = \bigl(\pi_E(v), \pi_F(w), \pr_{\bE} (\Phi(v) ) \otimes \pr_{\bF} ( \Psi(w))\bigr).
\end{equation}

Given another pair of local trivializations $(\widetilde V, \widetilde \Phi)$ and $(\widetilde W, \widetilde \Psi)$, let $\tau \colon V \cap \widetilde V \to \GL(\bE)$ and $\sigma \colon W \cap \widetilde W \to \GL(\bF)$ be the corresponding transition functions satisfying
\[
 (\widetilde \Phi \circ \Phi^{-1})(x,v) = (x, \tau(x) v), \qquad (\widetilde \Psi \circ \Psi^{-1})(y,w) = (y, \sigma(y) w).
\]
Then the corresponding transition function for $E \boxtimes F$ is given by
\[ (\widetilde \Theta \circ \Theta^{-1}) (x,y,v \otimes w) = (x, y, \theta(x,y) \cdot (v \otimes \omega) ) \]
where $(x,y) \mapsto \theta(x,y) \coleq \tau(x) \otimes \sigma(y)$ defines a smooth map $(U \cap \widetilde U) \times (V \cap \widetilde V) \to \GL(\bE \otimes \bF)$. By \cite[Lemma 10.6, p.~253]{Lee} $E \boxtimes F$ is a smooth vector bundle on $M \times N$ (cf.\ also \cite[Chap.~II, Problem 4, p.~84]{GHV}).

For a vector bundle $E \to M$ the space $\Gamma(M,E)$ of smooth sections of $E$ is endowed with the usual (F)-topology, the subspace $\Gamma_c(M,E)$ of compactly supported sections carries the usual (LF)-topology; similarly for $C^\infty(M)$ and its subspace $C^\infty_c(M)$ of compactly supported functions.

\subsection{Distribution theory}

For distribution theory we follow L.~Schwartz \cite{TD}. However, we will only need to consider \emph{real-valued} distributions, i.e., distributions will be continuous real-valued linear forms on corresponding spaces of real-valued test functions. In particular, $C^\infty(M)$ will be the space of smooth functions $M \to \bR$ and $\cD(M)$ the subspace of compactly supported functions for any manifold $M$ and in particular for open subsets $M \subseteq \bR^n$. These spaces are also equipped with their usual (F) and (LF)-topologies, respectively. $\delta_y \coleq \delta(.-y)$ denotes the delta distribution at $y \in \bR^n$.

For the vector valued case, given an open subset $\Omega \subseteq \bR^n$ and a finite dimensional vector space $\bE$, $\cD'(\Omega, \bE) = \Lin(\cD(\Omega), \bE)$ and $C^\infty(\Omega, \bE)$ are the spaces of distributions and smooth functions with values in $\bE$, respectively.

We will now recall some facts about distributions on manifolds (cf.~\cite{GKOS}).
For their definition we will not assume that our manifolds are oriented, hence distributions will be continuous linear forms on spaces of densities.

In the following, let $M$ be a manifold of dimension $n$. Although the objects to be integrated there naturally are densities, the notion of a locally integrable scalar function on $M$ is unambiguously defined by recourse to the Euclidean setting via local charts; the space of these locally integrable functions is denoted by $\Lint^1_{\mathrm{loc}}(M)$, where we identify functions which agree almost everywhere. More generally, a section of a vector bundle $E \to M$ is said to be \emph{locally integrable} if its components expressed in local charts are locally integrable in the usual sense. The space of (equivalence classes of) locally integrable sections of $E$ is denoted by $\Lint^1_\loc(M,E)$.

We denote by $\Vol(M)$ the line bundle of all (real) densities on $M$ \cite[Chap.~16, p.~427]{Lee}. Given a semi-Riemannian metric $g$ on $M$, its volume density is denoted by $\ud V_g$.
The space of distributions on $M$ then is defined as the strong dual $\cD'(M) \coleq (\Gamma_c(M, \Vol(M)))'$.
Every locally integrable function $f \in \Lint^1_{\mathrm{loc}}(M)$ defines a distribution $T_f \in \cD'(M)$ via the assignment
\[ \langle T_f, \omega \rangle \coleq \int f \omega \qquad (\omega \in \Gamma_c(\Vol(M))) \]
which gives a linear injection $\Lint^1_{\mathrm{loc}}(M) \to \cD'(M)$; we will usually denote $T_f$ by $f$ again. For a vector bundle $E \to M$, the space of $E$-valued distributions on $M$ is defined as
\begin{align*}
 \cD'(M,E) \coleq{} & ( \Gamma_c(M, E^* \otimes \Vol(M)))' \\
 \cong{} & \Gamma(M, E) \otimes_{C^\infty(M)} \cD'(M) \\
 \cong{} & \Hom_{C^\infty(M)} ( \Gamma(M, E^*), \cD'(M))
\end{align*}
(see.~\cite[Thm.~3.1.12, p.~239]{GKOS} for the isomorphisms), where $E^*$ is the dual bundle to $E$.
As in the scalar case, a locally integrable section $t \in \Lint^1_{\loc}(M, E)$ defines an element $T_t \in \cD'(M,E)$ via
\begin{gather*}
\langle T_t \cdot v, \omega \rangle = \langle T_t, v \otimes \omega \rangle \coleq \int (t\cdot v) \omega
\end{gather*}
for $v \in \Gamma(M, E^*)$ and $\omega \in \Gamma_c(M, \Vol(M))$, which gives a linear injection $\Lint^1_\loc(M, E) \to \cD'(M,E)$. Again, we will denote $T_t$ by $t$. Note that we use the symbol ``$\cdot$'' to denote contraction.

Let $N$ be another manifold, $F \to N$ a vector bundle and $\mu \colon E \to F$ a vector bundle isomorphism covering $\mu_0 \colon M \to N$. There is a natural pullback action $\mu^* \colon \cD'(N, F) \to \cD'(M, E)$ extending the pullback $\mu^* \colon \Lint^1_{\loc}(N,F) \to \Lint^1_{\loc}(M, E)$; it is given by
\[ \langle \mu^*(u), v \otimes \omega \rangle = \langle u, \mu_* v \otimes (\mu_0)_* \omega \rangle \qquad (v \in \Gamma(M, E^*),\ \omega \in \Gamma_c(M, \Vol(M))) \]
where $\mu_* v \in \Gamma(N,F)$ is given by $(\mu_* v)(p) \coleq \mu(v(\mu_0^{-1}(p)))$ and $(\mu_0)_*\omega \in \Gamma_c(N, \Vol(N))$ is the usual pushforward of densities.

For later use we recall that a net $(u_\e)_{\e \in (0,1]}$ in $\cD'(M,E)$ converges to $0$ if and only if for all open subsets $U \subseteq M$ where $E$ is trivial we have
\[ \langle u_\e^i, \omega \rangle \to 0 \qquad \forall i=1 \dotsc \dim E\ \forall \omega \in \Gamma_c(U, \Vol(M)) \]
where the coordinates $u_\e^i$ of $u_\e$ are taken with respect to any basis of $\Gamma(U,E)$. Clearly, it suffices to establish this for any subfamily of such sets $U$ which covers $M$.

For further background on distributional geometry we refer to \cite[Section 3.1]{GKOS}.

\subsection{Algebras of nonlinear generalized functions}

The basic problem in calculating the curvature of metrics like \eqref{equation.1.1} is that one cannot multiply distributions, in general. Our approach of regularizing them for this purpose amounts to embedding them into an algebra.

L.~Schwartz' impossibility result \cite{Schwartz} gives a clear limitation on what can be expected at all from an embedding of the space of distributions into an algebra. His original statement
\begin{quote}
 ``on ne peut, dans aucune th\'eorie, avoir \`a la fois une multiplication, une d\'erivation, et un \'el\'ement $\delta$''
\end{quote}
can be made precise as follows \cite[p.~27]{MOBook}:
\begin{quote}
There is no associative algebra $(\cA(\bR), \circ, +)$ satisfying
\begin{enumerate}[label=(\roman*)]
 \item $\cD'(\bR)$ is linearly embedded in $\cA(\bR)$, and the constant function $1$ is the unit in $\cA(\bR)$;
 \item there is a derivation operator $\pd$ on $\cA(\bR)$, that is, a linear map satisfying the Leibnitz rule;
 \item $\pd|_{\cD'(\bR)}$ coincides with the usual derivative;
 \item $\circ|_{C(\bR) \times C(\bR)}$ coincides with the usual product.
\end{enumerate}
\end{quote}

However, it is possible to construct such an algebra (which is even commutative) if one replaces condition (iv) by the more restricted condition
\begin{quote}
\begin{enumerate}
 \item[(iv')] $\circ|_{C^\infty(\bR) \times C^\infty(\bR)}$ coincides with the usual product.
\end{enumerate}
\end{quote}
First introduced by J.~F.~Colombeau \cite{ColNew,ColElem}, such algebras have found widespread applications in the study of problems simultaneously involving singularities, nonlinearities and differentiation.

The central idea of Colombeau algebras is to represent distributions by nets of smooth functions and form their product componentwise. However, this has to  be done in a specific way in order to obtain the desired consistency properties (i)--(iii) and (iv'). Instead of all nets of smooth functions one thus only considers those which have moderate growth and factors out these which have negligible growth in a certain sense; this quotient then forms the Colombeau algebra.

It is noteworthy that these nets of smooth functions contain more information than would be necessary for representing distributions. In fact, a distribution can be approximated in many ways which are equivalent from a distributional point of view but different in the Colombeau algebra. For example, consider the two approximations of the delta distribution obtained from $\varphi \in \cD(\bR)$ with $\supp \varphi \subseteq [-1,1]$ and $\int \varphi(x)\,\ud x = 1$ by setting
\[ \rho_\e(x) \coleq \frac{1}{\e} \varphi \left( \frac{x+\e}{\e} \right), \qquad \tilde \rho_\e(x) \coleq \frac{1}{\e} \varphi \left( \frac{x-\e}{\e} \right). \]
While both $\rho_\e$ and $\tilde\rho_\e$ converge to $\delta$ distributionally, their pointwise product $\rho_\e \cdot \tilde \rho_\e$ vanishes everywhere. At the same time, the product $\rho_\e \cdot \rho_\e$ is nonzero but does not converge distributionally \cite[Ex.~2.2, p.~24]{MOBook}. This ambiguity can be understood as a distinctive microscopic structure of the elements of the Colombeau algebra which is invisible from a macroscopic or purely distributional point of view. The latter enters the picture through the concept of association (see \Autoref{definition.3.53} and \Autoref{definition.5.12} below): whenever a net of smooth functions representing a Colombeau generalized function $R$ has a limit in the space of distributions, this limit is called the distribution \emph{associated} to the nonlinear generalized function $R$, or its \emph{distributional shadow}. Association defines a relation on a subspace of the Colombeau algebra which is coarser than equality; moreover, it is not compatible with multiplication and nonlinear operations in general, but only with classical operations from distribution theory. As a rule of thumb, one can say that a calculation which makes sense in the context of distributions can be repeated in the Colombeau algebra and the result will be associated to the former distributional result.

If, however, in the Colombeau algebra one performs operations on the distributions which are ill-defined classically, two things can happen: either one obtains a result which is associated to a distribution; this means that the microscopic structure of the generalized functions does not matter in the end. Or one obtains a generalized function which is not associated to a distribution; this is then a genuinely new object which can be studied in more detail, for example by prescribing a peculiar regularization procedure on physical grounds.

We will see that in the case of the curvature of the conical metric, which cannot be calculated classically, we do obtain an associated distribution in the end.

\subsection{Regularization of distributions}

A proper handling of the regularization procedure which represents a distribution by a net of smooth functions is indispensable for the understanding of nonlinear generalized functions. We will examine this in the setting of $\bR^n$ first. For this purpose, most standard textbooks on distribution theory use some form of delta net
\begin{gather*}
(\rho_\e)_\e \in \cD(\bR^n)^{(0,1]},\quad \supp \rho_\e \to \{0\}\ \textrm{ for }\e \to 0,\\
\quad \int \rho_\e(x)\,\ud x = 1,\quad \sup_{\e \in (0,1]} \norm{\rho_\e}_{L^1} < \infty.
\end{gather*}
Given $u \in \cD'(\bR^n)$, because $\rho_\e$ converges to $\delta$ the convolution $u * \rho_\e$ (which is a smooth function) converges to $u * \delta = u$, hence is an approximation of $u$ by a net of smooth functions. A similar procedure applies to distributions on open subsets $\Omega$ of $\bR^n$ by using cut-off functions.

It is helpful now to focus on the original intent of regularization, which is to represent distributions by smooth functions. If this is to be done in a linear and continuous way, we end up exactly with the space $\Lin(\cD'(\Omega), C^\infty(\Omega))$ of so-called \emph{smoothing operators} on $\Omega$. For a net $(\Phi_\e)_\e$ in this space to be able to give an accurate representation of the distribution to which it is applied, it should converge to the identity in $\Lin(\cD'(\Omega), \cD'(\Omega))$. A different viewpoint, which will be quite important for us, is obtained through L.~Schwartz' kernel theorem in the form
\begin{equation}\label{equation.2.4}
\Lin(\cD'(\Omega), C^\infty(\Omega)) \cong C^\infty(\Omega, \cD(\Omega))
\end{equation}
\cite[Th\'eor\`eme 3, p.~127]{FDVV}. Explicitly, this isomorphism says that smoothing operators $\Phi \in \Lin(\cD'(\Omega), C^\infty(\Omega))$ are in bijective correspondence with \emph{smoothing kernels} $\vec\varphi \in C^\infty(\Omega, \cD(\Omega))$ via the relations
\begin{align*}
\Phi(u)(x) &= \langle u(y), \vec\varphi(x)(y) \rangle,\\
\vec\varphi(x)(y) &= \Phi(\delta_y)(x).
\end{align*}
For $\Phi_\e(u) = u * \rho_\e$ the corresponding kernel is $\vec\varphi_\e(x)(y) \coleq \rho_\e(x-y)$, so convolution (which was used in Colombeau's original approach) is a mere special case of this situation.

Note that when we speak of smoothing operators $\Phi,\Psi$ we will often, without further mention, denote the corresponding kernels by $\vec\varphi, \vec\psi$, and vice versa; for nets $\bPhi = (\bPhi_\e)_\e, \bPsi = (\bPsi_\e)_\e$ we usually denote the corresponding nets of kernels by $\bvarphi = (\bvarphi_\e)_\e, \bpsi = (\bpsi_\e)_\e$.

The virtue of the formalism of smoothing operators -- as opposed to convolution -- is that it transfers directly to manifolds. In fact, to regularize a distribution on a manifold $M$ it is natural to employ elements of the space
\begin{equation}\label{equation.2.5}
\Lin ( \cD'(M), C^\infty(M)) \cong C^\infty(M, \Gamma_c(M, \Vol(M)))
\end{equation}
and as above, if we want nets $(\bPhi_\e)_\e$ to accurately represent the distributions they are applied to then we demand $\bPhi_\e \to \id$ in $\Lin(\cD'(M), \cD'(M))$.

In the vector valued case the reasoning is similar: distributions with values in a vector bundle $E \to M$ are regularized by applying elements of
\[ \Lin( \cD'(M, E), \Gamma(M, E)) \]
which we call \emph{vector smoothing operators}.

A cornerstone of our theory (see \cite{vecreg}) will be the isomorphism of locally convex spaces
\begin{equation}
\label{equation.2.6}
\begin{aligned}
\Lin(\cD'(&M_y, E), \Gamma(M_x, E)) \\
& \cong \Gamma(M_x\times M_y, E \boxtimes E^*) \otimes_{C^\infty(M_x \times M_y)} \Lin(\cD'(M_y), C^\infty(M_x)).
\end{aligned}
\end{equation}
Here, we added variable names to the manifold $M$ to distinguish its two instances.
The vector smoothing operator $\Theta \in \Lin(\cD'(M, E), \Gamma(M, E))$ which corresponds to a tensor product $A \otimes \Phi$ with $A \in \Gamma(M \times M, E \boxtimes E^*)$ and $\Phi \in \Lin(\cD'(M), C^\infty(M))$ is given by
\begin{equation}
 \label{equation.2.7}
(\Theta (t) \cdot v)(x) = \langle t(y), ( A(x,y) \cdot v(x) ) \otimes \vec\varphi(x)(y) \rangle
\end{equation}
for $v \in \Gamma(M, E^*)$, $t \in \cD'(M, E)$ and $x \in M$, where $\vec\varphi \in C^\infty(M, \Gamma_c(M, \Vol(M)))$ corresponds to $\Phi$ via \eqref{equation.2.5}.
To clarify this formula, note that $A(x,y) \cdot v(x)$ denotes the contraction of
\[ A \in \Gamma(M_x \times M_y, E \boxtimes E^*) \cong \Hom_{C^\infty(M_x)} ( \Gamma(M_x, E^*), C^\infty(M_x, \Gamma(M_y, E^*))) \]
with $v \in \Gamma(M_x, E^*)$, giving an element of $C^\infty(M_x, \Gamma(M_y, E^*))$. The product with $\vec\varphi \in C^\infty(M_x, \Gamma_c(M_y, \Vol(M)))$ gives an element of $C^\infty(M_x, \Gamma_c(M_y, E^* \otimes \Vol(M)))$ which, by applying $t$, is mapped to a function in $C^\infty(M_x)$.
Conversely,
given $\Theta$ we obtain an element
\[ T \in \Hom_{C^\infty(M_x \times M_y)} ( \Gamma(M_x \times M_y, E^* \boxtimes E), \Lin( \cD'(M_y), C^\infty(M_x) )) \]
by
\begin{equation}\label{equation.2.8}
 T ( v \otimes w)(t)(x) = ( \Theta ( t \otimes w) \cdot v)(x)
\end{equation}
for $v \in \Gamma(M, E^*)$, $w \in \Gamma(M, E)$, $\omega \in \Gamma(M, E)$, $t \in \cD'(M)$ and $x \in M$.

While isomorphism \eqref{equation.2.6} is reasonable to expect, two points are important to note: first, by combining arbitrary transport operators $A \in \Gamma(M \times M, E \boxtimes E^*)$ and smoothing operators $\Phi \in \Lin(\cD'(M), C^\infty(M))$ one really obtains \emph{all} vector smoothing operators on $E$. And second, this is in fact an isomorphism in the category of locally convex spaces, i.e., convergence of a net $\bTheta_\e = \bA_\e \otimes \bPhi_\e$ to the identity on $\cD'(M,E)$ is obtained by requiring $\bA_\e$ and $\bPhi_\e$ to converge correspondingly (see \cite{vecreg} and \Autoref{theorem.3.29} below).

Let us introduce the following abbreviations.
\begin{definition}\label{definition.2.1}
Let $\Omega \subseteq \bR^n$ be open, $\bE$ a finite dimensional vector space, $M$ a manifold and $E \to M$ a vector bundle. We set
 \begin{align*}
  \SK(\Omega) & \coleq C^\infty(\Omega, \cD(\Omega)), \\
  \SO(\Omega) & \coleq \Lin(\cD'(\Omega), C^\infty(\Omega)), \\
  \SK(M) & \coleq C^\infty(M, \Gamma_c(M, \Vol(M))), \\
  \SO(M) & \coleq \Lin(\cD'(M), C^\infty(M)), \\
  \TO(\Omega, \bE) & \coleq C^\infty(\Omega \times \Omega, \bE \otimes \bE^*), \\
  \TO(M, E) & \coleq \Gamma(M \times M, E \boxtimes E^*), \\
  \VSO(\Omega, \bE) & \coleq \Lin(\cD'(\Omega, \bE), C^\infty(\Omega, \bE)), \\
  \VSO(M, E) &\coleq \Lin(\cD'(M, E), \Gamma(M, E)).
 \end{align*}
 Here, $\SK$ stands for \emph{smoothing kernels}, $\SO$ for \emph{smoothing operators}, $\TO$ for \emph{transport operators} and $\VSO$ for \emph{vector smoothing operators}.
\end{definition}

In the local case $(\Omega \subseteq \bR^n$ open, $\bE$ a finite-dimensional vector space) isomorphism \eqref{equation.2.6} holds in the form
\[ \Lin(\cD'(\Omega, \bE), C^\infty(\Omega, \bE)) \cong \TO(\Omega, \bE) \otimes_{C^\infty(\Omega \times \Omega)} \SO(\Omega) \]
with \eqref{equation.2.7} and \eqref{equation.2.8} adapted accordingly.

For later use, let us mention that diffeomorphisms $\mu$ and (local) vector bundle isomorphisms $F$ act naturally on the respective spaces in \Autoref{definition.2.1}. We denote the corresponding pullback mappings by $\mu^*$ and $F^*$ and set $\mu_* \coleq (\mu^{-1})^*$, $F_* \coleq (F^{-1})^*$. Given a chart $(U,\varphi)$ on $M$ there are canonical isomorphisms
\[ \SK(\varphi(U)) \cong \SK(U), \qquad \SO(\varphi(U)) \cong \SO(U), \]
and if in addition $E$ is trivial on $U$ with trivializing fiber $\bE$ there are canonical isomorphisms
\[ \TO(\varphi(U), \bE) \cong \TO(U,E), \qquad \VSO(\varphi(U), \bE) \cong \VSO(U, E). \]

The representation of vector smoothing operators by transport operators and smoothing operators is particularly useful because it decouples, in the smoothing of vector distributions, the geometric component from the smoothing component. This allows us to use a pair $(A, \Phi)$ to regularize distributional sections also on other vector bundles obtained functorially from $E$, as will be explained in the following subsection.

\subsection{Functorial constructions}

We quickly recall how functorial constructions can be applied to vector bundles (see \cite[6.7, p.~53]{KMS}). Note that in our terminology functors may be co- or contravariant and may depend on more than one argument.

Let $\catVec$ be the category of finite dimensional real vector spaces with linear mappings as morphisms. A covariant functor $\lambda \colon \catVec \to \catVec$ is called \emph{smooth} if for all objects $V,W$ in $\catVec$ the mapping $\lambda \colon \Lin(V,W) \to \Lin(\lambda V, \lambda W)$ is smooth. 

Given a vector bundle $E$ with typical fiber $\bE$, such a functor gives rise to a vector bundle $\lambda E = \bigsqcup_x \lambda(E_x)$ with typical fiber $\lambda \bE$. If $\Phi \colon \pi^{-1}(U) \to U \times \bE$ is a local trivialization of $E$ and $\pi_{\lambda E} \colon \lambda E \to M$ the projection of $\lambda E$ onto the base manifold then a local trivialization $\lambda \Phi \colon \pi_{\lambda E}^{-1}(U) \to U \times \lambda \bE$ is obtained by setting $(\lambda \Phi)(v) \coleq (p, \lambda (\Phi|_{E_p})\cdot v)$ for $v \in \lambda(E_p)$. Here, $\Phi|_{E_p} \colon E_p \to \{ p \} \times \bE \cong \bE$ gives rise to a linear map $\lambda(\Phi|_{E_p}) \colon \lambda (E_p) \to \{ p \} \times \lambda \bE \cong \lambda \bE$. Moreover, if $\Phi_\alpha \colon \pi^{-1}(U_\alpha) \to U_\alpha \times \bE$ and $\Phi_\beta \colon \pi^{-1}(U_\beta) \to U_\beta$ are local trivializations with $U_\alpha \cap U_\beta \ne \emptyset$, the corresponding transition function $\tau_{\alpha\beta} \colon U \cap V \to \GL(\bE)$ which satisfies $(\Phi_\alpha \circ \Phi_\beta^{-1}) (p, v) = (p, \tau_{\alpha\beta}(p)\cdot v)$ gives rise to the transition function $(\lambda \tau)_{\alpha \beta} \colon U \cap V \to \GL(\lambda \bE)$ given by $(\lambda \tau)_{\alpha \beta}(x) = \lambda ( \tau_{\alpha \beta}(x))$; indeed, we then have $(\lambda \Phi_\alpha \circ \lambda \Phi_\beta^{-1})(p,v) = (p, (\lambda \tau)_{\alpha\beta}(p)\cdot v)$.

This procedure also works for contravariant functors if one applies $\lambda$ to the respective inverse mappings instead, i.e., by taking $(\lambda \Phi)(v) \coleq (p, \lambda(\Phi|_{E_p}^{-1})\cdot v)$ for the local trivialization and $(\lambda \tau)_{\alpha\beta}(x) = \lambda ( \tau_{\alpha\beta}^{-1} (x) )$ for the transition function. Analogously, one can apply functors in any number of arguments, covariant in some and contravariant in the other arguments, to obtain new vector bundles.

We will restrict our attention to functors of a particular kind. First, recall that for a functor
\[ \lambda \colon \catVec^n = \catVec \times \dotsm \times \catVec \to \catVec \]
and vector spaces $\bE_i$, $\bF_i$ ($i=1 \dotsc n$) it makes sense to talk of multilinearity of the mapping
\[ \lambda \colon \Lin(\bE_1, \bF_1) \times \dotsm \times \Lin(\bE_n, \bF_n) \to \Lin( \lambda(\bE_1, \dotsc, \bE_n), \lambda( \bF_1, \dotsc, \bF_n)) \]
(where the factor $\Lin(\bE_j, \bF_j)$ in the domain is replaced by $\Lin(\bF_j, \bE_j)$ if $\lambda$ is contravariant in the $j$th factor). If this mapping is multilinear for all choices of $\bE_i$ and $\bF_i$ we say that the functor $\lambda$ is \emph{multilinear}.

\begin{definition}\label{definition.2.2}We say that a functor $\lambda \colon \catVec^k \to \catVec$ is $(n_1, \dotsc, n_k)$-ho\-mo\-ge\-neous if there exists a multilinear functor
 \[ \sigma \colon \catVec^{n_1} \times \dotsm \times \catVec^{n_k} \to \catVec \]
 such that
\[ \lambda = \sigma \circ ( \delta_{n_1} \times \dotsm \times \delta_{n_k} ) \]
where for $n \in \bN$, $\delta_n \colon \catVec \to \catVec^n$ is the diagonal functor $\bE \mapsto \bE^n$, $f \mapsto f \times \dotsm \times f$.
\end{definition}

\begin{definition}
 Let $\bE_1, \dotsc, \bE_k, \bF$ be finite dimensional vector spaces. A mapping $f \colon \bE_1 \times \dotsm \times \bE_k \to \bF$ is called an $(n_1, \dotsc, n_k)$-homogeneous polynomial if there exists a multilinear mapping
\[ h \colon \bE_1^{n_1} \times \dotsm \times \bE_k^{n_k} \to \bF \]
such that
\[ f(x_1, \dotsc, x_k) = h(\Delta_{n_1}(x_1), \dotsc, \Delta_{n_k}(x_k)) \]
where for $n \in \bN$, $\Delta_n \colon \bE \to \bE^n$, $x \mapsto (x, \dotsc, x)$ is the diagonal mapping.

\end{definition}
Consequently, if $\lambda \colon \catVec^k \to \catVec$ is an $(n_1, \dotsc, n_k)$-homogeneous functor then for all vector spaces $\bE_i$, $\bF_i$ the mapping
 \[ \lambda \colon \Lin(\bE_1, \bF_1) \times \dotsm \times \Lin(\bE_k, \bF_k) \to \Lin(\lambda(\bE_1, \dotsc, \bE_k), \lambda(\bF_1, \dotsc, \bF_k)) \]
(where the factor $\Lin(\bE_j, \bF_j)$ in the domain is replaced by $\Lin(\bF_j, \bE_j)$ if $\lambda$ is contravariant in the $j$th factor) is an $(n_1, \dotsc, n_k)$-homogeneous polynomial.

From now on, a \emph{homogeneous functor} will always mean a functor $\catVec^k \to \catVec$ as in \Autoref{definition.2.2} (usually in one argument, i.e., with $k=1$) without mentioning the degrees of homogeneity. Similarly, a smooth functor will mean a smooth functor $\catVec^k \to \catVec$.

\begin{lemma}
\begin{enumerate}[label=(\roman*)]
 \item Every homogeneous functor $\catVec \times \dotsm \times \catVec \to \catVec$ is smooth.
 \item The composition of homogeneous functors is a homogeneous functor.
 \item The functors $\bE$, $\bE^*$, $\bE \otimes \bF$, $\bE^r_s = \underbrace{\bE \otimes \dotsm \otimes \bE}_{\textrm{$r$ factors}} \otimes \underbrace{\bE^* \otimes \dotsm \otimes \bE^*}_{\textrm{$s$ factors}}$ ($(r,s)$-tensors), $S^k \bE$ (symmetric tensors), $A^k \bE$ (antisymmetric tensors) are homogeneous.
 \end{enumerate}
\end{lemma}
\begin{proof}
(i) follows because the diagonal mappings $\Delta_n$ as well as multilinear mappings between finite dimensional vector spaces are smooth; (ii) and (iii) are immediate.
\end{proof}

A smooth functor also acts on transport operators in the following way:

\begin{definition}\label{definition.2.5}Let $\lambda \colon \catVec^k \to \catVec$ be a smooth functor in $k$ arguments and, for each $i=1 \dotsc k$, $E_i \to M$ a vector bundle over a manifold $M$ and $A_i \in \TO(M, E_i)$. We define $\lambda (A_1, \dotsc, A_k) \in \TO(M, \lambda(E_1,\dotsc,E_k))$ by
\[ (\lambda(A_1, \dotsc, A_k))(x,y) \coleq \lambda ( A_1(x,y), \dotsc, A_k(x,y)), \]
with $A_i(x,y)$ replaced by $A_i(y,x)$ if $\lambda$ is contravariant in the $i$th factor.

Similarly, if $\Omega \subseteq \bR^n$ is open, $\bE_i$ is a family of finite-dimensional vector spaces and $a_i \in \TO(\Omega, \bE_i)$ for each $i=1 \dotsc k$, we define $\lambda (a_1,\dotsc,a_k)$ as the element of $\TO(\Omega, \lambda(\bE_1, \dotsc, \bE_k))$ given by
\[ (\lambda (a_1,\dotsc,a_k))(x,y) \coleq \lambda ( a_1(x,y), \dotsc, a_k(x,y)), \]
again with $a_i(x,y)$ replaced by $a_i(y,x)$ if $\lambda$ is contravariant in the $i$th factor.
\end{definition}

The following is easily seen:
\begin{lemma}\label{lemma.2.6}The mappings $(A_1,\dotsc,A_k) \mapsto \lambda (A_1,\dotsc,A_k)$ and $(a_1,\dotsc,a_k) \mapsto \lambda (a_1,\dotsc,a_k)$ of \Autoref{definition.2.5} are smooth.
\end{lemma}

For the treatment of Lie derivatives we recall a useful notion: a \emph{vector bundle functor} is a functor $F$ which associates a vector bundle $F(M) \to M$ to each manifold $M$ and a vector bundle homomorphism $F(f) \colon F(M) \to F(N)$ to each $f \in C^\infty(M,N)$ (where $N$ is another manifold). Our prime example for this will be the tangent bundle functor $\tang M$. Moreover, $\lambda(\tang M)$ is a vector bundle functor for every homogeneous functor (in fact, even for every smooth functor) \cite[6.14, p.~56]{KMS}.

Given a vector bundle functor $F$ one can define the \emph{Lie derivative} of a smooth section $s \in \Gamma(M, F(M))$ with respect to a vector field $X \in \fX(M)$ by
\[ \Lie_X s \coleq \left.\frac{\ud}{\ud t}\right|_{t=0} ( \Fl^X_t)^* s \in \Gamma(M, F(M)) \]
where $\Fl^X_t$ is the flow of $X$ at time $t$.

If $f \colon M \to N$ is a diffeomorphism there is an induced map $f^* \colon F(N) \otimes F(N)^* \to F(M) \otimes F(M)^*$ which allows us to consider the Lie derivative of transport operators:

\begin{definition}Let $F$ be a vector bundle functor and $M$ a manifold.
 For any $X \in \fX(M)$ we define the Lie derivative of $A \in \TO(M, F(E))$ by
 \[ \Lie_X A \coleq \left. \frac{\ud}{\ud t}\right|_{t=0} (\Fl^X_t)^* A \in \TO(M,F(E)). \]
\end{definition}

Note that if $\lambda$ is a homogeneous functor then $\lambda \circ F$ is also a vector bundle functor. The notation $\Lie_X A$ actually is a shorthand for $\Lie_{X \times X} A$ and by the product rule we have $\Lie_X A = \Lie_{X \times X} A = \Lie_{X \times 0} A + \Lie_{0 \times X} A$.

\section{Generalized sections}\label{section.3}

In this section we present in detail the framework enabling us to perform nonlinear operations on distributions on manifolds and calculate the distributional curvature tensor. Our approach is based on Colombeau's space of generalized functions \cite{ColNew,ColElem}. However, we incorporate in our construction several extensions which were developed during the last couple of years: a functional analytic formulation \cite{papernew}, which led to a geometrization of the theory \cite{bigone}; a representation of smoothing operators for vector distributions \cite{vecreg}; and a study of locality conditions, which are used to obtain the sheaf property \cite{specfull}. Only through the combination of these results the scope of the theory of nonlinear generalized functions can be extended to tensor fields in sufficient generality while at the same time it is concrete enough such that calculations with singular  metrics can be performed.

\subsection{Preliminaries}

In \cite{bigone} we gave a construction of spaces of nonlinear generalized sections of vector bundles containing distributional sections such that common differential geometric operations (tensor products, Lie derivatives, covariant derivatives, contractions) are well-defined. While this settled most technical questions that will be relevant for us, it does not yet incorporate the extra structure provided by the fact that in semi-Riemannian geometry we usually only work with the family of bundles $\lambda E$ obtained functorially by applying a functor $\lambda$ to the tangent bundle $E = \tang M$.

Before we begin to adapt that construction for our purposes, let us recall its basics. Following the general principle introduced in \cite{papernew}, the basic space of nonlinear generalized sections of a vector bundle $E \to M$ consists of families of corresponding smooth sections which are parameterized by the space of smoothing operators $\VSO(M, E)$. A nonlinear generalized section then is a mapping $R \colon \VSO(M, E) \to \Gamma(M,E)$, and a distribution $u \in \cD'(M, E)$ canonically corresponds to the mapping which assigns to each $\Theta \in \VSO(M,E)$ its value $\Theta(u)$. 

Most operations on smooth sections can be extended to generalized sections simply by applying them for fixed $\Theta$. However, in case sections over different vector bundles are involved this requires the dependence of generalized sections on different smoothing operators. For example, suppose $F \to M$ is another vector bundle and a generalized section is given by $S \colon \VSO(M,F) \to \Gamma(M,F)$. Forming the tensor product of $R$ and $S$ pointwise, i.e., for fixed smoothing operators, we end up with the mapping $\VSO(E) \times \VSO(F) \to \Gamma(E \otimes F)$, $(\Theta, \Xi) \to R(\Theta) \otimes S(\Xi)$.

More generally, a basic space depending on an index set $\Delta$ of vector bundles was defined \cite[Definition 4.1, p.~191]{bigone} by
\[ \cE^\Delta(M, E) \coleq C^\infty\bigl( \prod_{G \in \Delta} \VSO(G), \Gamma(M, E) \bigr). \]
Note that smooth dependence (in the sense of \cite{KM}) on smoothing operators is necessary to define a Lie derivative which commutes with the embedding of distributional tensor fields. While the usual operations on smooth sections are easily extended to this space, the quotient construction which makes sure that one retains as much compatibility to distributional and classical calculus as possible needs more work; we refer to \cite{bigone} for further details.

In this article we will make the following amendments to the theory outlined above:
\begin{itemize}
\item We \emph{modularize smoothness} of representatives. Only this smoothness makes it possible to define a (geometric) Lie derivative of nonlinear generalized functions extending the Lie derivative of distributions via the canonical embedding. If one does not need this Lie derivative one can do a simpler theory without smoothness. We develop both theories in parallel here.
 \item We apply \emph{locality conditions} in the spirit of \cite{specfull} to obtain the sheaf property from a suitably chosen basic space.
 \item We replace dependence on smoothing operators on $\lambda E$ by dependence on smoothing operators on $E$ only; these will then be \emph{functorially extended} to smoothing operators on $\lambda E$ in the testing procedure.
 \item We develop a complete \emph{local formalism} of the theory.
 \end{itemize}

This way, we will obtain a space of nonlinear generalized sections on any vector bundle $\lambda E$ obtained functorially from $E$ such that the embedding of distributions also preserves symmetry properties. For example, if $u \otimes u \in \cD'(M, E \otimes E)$ is given, then $(A \otimes \Phi)(u \otimes u)$ is symmetric for $A = B \otimes B$ with $B \in \TO(M, E)$ but not for general $A \in \TO(M, E \otimes E)$.

\subsection{The basic space}

\newcommand{\sm}{\textbf{smooth theory}}
\newcommand{\si}{\textbf{simple theory}}

We will first develop the theory for general vector bundles $E$ and specialize this to the tangent bundle $E = \tang M$ in Section \ref{section.4}. As mentioned above, we will develop the \sm\ and the \si\ in parallel. We will use the same notation in both cases.

\begin{definition}\label{definition.3.1}
Let $M$ be a manifold, $E$ a vector bundle over $M$ and $\lambda$ a homogeneous functor in one variable.
\begin{enumerate}[label=(\roman*)]
 \item The basic space of \emph{generalized sections} of $\lambda E$ is defined in the \si\ as
\begin{align*}
&\cE( M, \lambda E) \coleq \{ R \colon \TO(M, E) \times \SO(M) \to \Gamma(M, \lambda E) \} \\
\intertext{and in the \sm\ as}
&\cE ( M, \lambda E ) \coleq C^\infty \bigl( \TO(M, E) \times \SO(M), \Gamma(M, \lambda E) \bigr).
\end{align*}
\item For $\lambda E = M \times \bR$ we call
\[ \cE(M) \coleq \cE(M, M \times \bR) 
\]
the basic space of \emph{generalized functions} on $M$.
\item The canonical embeddings of $u \in \cD'(M, \lambda E)$ and $t \in \Gamma(M, \lambda E)$ into $\cE(M, \lambda E)$ are defined as
\begin{align*}
(\iota u)(A, \Phi) & \coleq (\lambda A \otimes \Phi)(u), \\
(\sigma t) (A, \Phi) & \coleq t.
\end{align*}
\end{enumerate}
\end{definition}

Note that although it is not visible in the notation, the space of generalized functions still depends on $E$ and transport operators on it! The need for this becomes evident for example from considering the contraction of $R \in \cE(M,E)$ with $S \in \cE(M, E^*)$, which gives an element in $\cE(M)$ that clearly depends on transport operators.

For the embedding $\iota$ to be well-defined in case of the \sm\ we have to know that $\iota u$ is smooth, but this follows because the mapping $A \mapsto \lambda A$ is smooth by \Autoref{lemma.2.6} and $(B, \Phi) \mapsto B \otimes \Phi$ is linear and continuous.

\begin{remark}One could in principle construct a similar theory with $\lambda$ being a functor in more than one variable, but as we will only need tensor calculus (hence $E = TM)$ for our purposes we will not build the theory in this generality.
\end{remark}

The space $\cE(M, \lambda E)$ inherits a $C^\infty(M)$-module structure from the codomain $\Gamma(M, \lambda E)$ of its elements. By reduction to trivial bundles one can show the isomorphisms of $C^\infty(M)$-modules
\begin{equation}\label{equation.3.9}
\begin{aligned}
\cE(M, \lambda E) & \cong \Gamma(M, \lambda E)\otimes_{C^\infty(M)} \cE(M) \\
& \cong \Hom_{C^\infty(M)} ( \Gamma(M, (\lambda E)^*), \cE(M))
\end{aligned}
\end{equation}
(see \cite[Theorem 4.4, p.~192]{bigone}) explicitly given by
\[ (R \cdot v)(A, \Phi) = R(A, \Phi) \cdot v. \]
for $R \in \cE(M, \lambda E)$ and $v \in \Gamma(M, (\lambda E)^*)$. In other words, generalized sections are given by smooth sections with coefficients in $\cE(M)$.

We will also need a local version of \Autoref{definition.3.1}.

\begin{definition}\label{definition.3.3} Let $\Omega \subseteq \bR^n$ be open, $\bE$ a finite-dimensional vector space and $\lambda$ a homogeneous functor in one variable.
\begin{enumerate}[label=(\roman*)]
 \item We define the basic space of \emph{generalized functions on $\Omega$ with values in $\lambda \bE$} in the \si\ as
\begin{align*}
\cE(\Omega, \lambda \bE) \coleq \{ R \colon \TO(\Omega, \bE) \times \SO(\Omega) \to  C^\infty(\Omega, \bE) \} \\
\intertext{and in the \sm\ as}
\cE(\Omega, \lambda \bE) \coleq C^\infty \bigl( \TO(\Omega, \bE) \times \SO(\Omega), C^\infty(\Omega, \bE) \bigr).
\end{align*}
\item In case of the constant functor $\lambda \bE = \bR$ we call $\cE(\Omega) \coleq \cE(\Omega, \bR)$ the basic space of \emph{generalized functions} on $\Omega$.
\item We define the canonical embeddings $\iota \colon \cD'(\Omega, \lambda \bE) \to \cE(\Omega, \lambda \bE)$ as well as $\sigma \colon C^\infty(\Omega, \lambda \bE) \to \cE(\Omega, \lambda \bE)$ by
\begin{align*}
 (\iota u) (a, \Phi) &\coleq (\lambda a \otimes \Phi)(u), \\
 (\iota f)(a, \Phi) &= f.
\end{align*}
\end{enumerate}
\end{definition}

\Autoref{definition.3.1} and \Autoref{definition.3.3} are compatible: if $(U,\varphi)$ is a chart on $M$ over which $E$ (and hence $\lambda E$) is trivial and we consider the obvious canonical isomorphisms
then the following diagrams commute:
\[
 \xymatrix{
 \cE(U, \lambda E) \ar[r]^-\cong & \cE(\varphi(U), \lambda \bE) & \cE(\varphi(U), \lambda  \bE) \ar[r]^-\cong & \cE(U, \lambda E) \\
 \cD'(U, \lambda E) \ar[u]^\iota \ar[r]^-\cong & \ar[u]^\iota \cD'(\varphi(U), \lambda \bE) & C^\infty(U, \lambda \bE) \ar[r]^-\cong \ar[u]^\sigma & \Gamma(U, \lambda E) \ar[u]^\sigma
 }
\]
Note that a coordinatewise embedding would depend on the choice of bases and contradict the Schwartz impossibility result \cite[Section 4]{global2}.

If $R \in \cE(U, \lambda E)$ corresponds to $R_0 \in \cE(\varphi(U), \lambda \bE)$ we call $R_0$ the local expression of $R$.

To make the isomorphism $\TO(U,E) \cong \TO(\varphi(U), \bE)$ explicit, suppose we are given a trivialization $\Phi \colon \pi_E^{-1}(U) \to U \times \bE$ of $E$ over an open subset $U \subseteq M$. Then a trivialization of $E^*$ over $U$ is given by
\begin{equation}\label{equation.3.10}
\begin{aligned}
&\widetilde \Phi \colon \pi_{E^*}^{-1}(U) \to U \times \bE^*, \\
& \widetilde \Phi(\omega) = (x, (\Phi|_{E_x}^{-1})^*\omega) \qquad (\omega \in E_x^*).
\end{aligned}
\end{equation}
By \eqref{equation.2.3} a trivialization $\Psi \colon \pi^{-1}(U \times U) \to U \times U \times (\bE \otimes \bE^*)$ of $E \boxtimes E^*$ is then given by
\begin{equation}\label{equation.3.11}
\Psi(v \otimes w) = (x,y, (\pr_\bE \circ \Phi)(v) \otimes (\pr_{\bE^*} \circ \widetilde \Phi)(w)) \qquad (v \in E_x, w \in E_y^*)
\end{equation}
which induces an isomorphism
\[ \TO(U, E) \cong C^\infty(U \times U, \Lin(\bE, \bE)). \]
Given a basis $(e_i)_i$ of $\bE$ with dual basis $(\e^j)_j$, a transport operator $A \in \TO(U, E)$ can be identified with an $m \times m$-matrix $(A^i_j)_{i,j=1\dotsc m}$ of smooth functions on $U \times U$ (where $m = \dim E$), determined by
\[ A^i_j(x,y) = (\e^i \circ \pr_{\bE} \circ \Phi \circ A(x,y) \circ \Phi^{-1})(y,e_j) \qquad (x,y\in U). \]
In case $U$ is also the domain of a chart $\varphi$ we define the local expression of $A$ by
\begin{align*}
  a(x,y)v &= (\pr_{\bE} \circ \Phi \circ A(x,y) \circ \Phi^{-1})(\varphi^{-1}(y),v) \\
  &= ( \Phi_{\varphi^{-1}(x)} \circ A(x,y) \circ \Phi_{\varphi^{-1}(y)}^{-1}  )(v) \qquad (x,y\in \varphi(U),\ v\in \bR^m)
\end{align*}
and the local expression of $A^i_j$ by
\[ a^i_j(x,y) \coleq A^i_j(\varphi^{-1}(x), \varphi^{-1}(y)) \qquad (x,y \in \varphi(U)). \]

A cornerstone of our theory will be that it is intrinsically diffeomorphism invariant:

\begin{definition}\label{definition.3.4}Let $M, N$ be manifolds, $E \to M$ and $F \to N$ vector bundles and $\mu \colon E \to F$ a vector bundle isomorphism covering $\mu_0 \colon M \to N$. We then define $\mu^* \colon \cE(N, F) \to \cE(M, E)$ by
 \[ (\mu^* R)(A, \Phi) \coleq \mu^* ( R(\mu_* A, (\mu_0)_* \Phi)). \]
\end{definition}

The central result justifying this definition is that it is compatible with the pullback of distributions and smooth sections:

\begin{theorem}\label{theorem.3.5} In the situation of \Autoref{definition.3.4}, let $\lambda$ be a homogeneous functor, $u \in \cD'(N, \lambda F)$ and $t \in \Gamma(N, \lambda F)$. We have $\iota(\mu^* u) = \mu^* ( \iota u)$ and $\sigma( \mu^* t) = \mu^* (\sigma t)$.
\end{theorem}
\begin{proof}
We essentially only have to plug in all definitions; for the first claim, we have
\begin{align*}
\mu^*(\iota u)(A, \Phi)(x) &= \mu^*( (\iota u)(\mu_* A, (\mu_0)_* \Phi))(x) \\
&= \mu^* ( ( \lambda ( \mu_* A) \otimes (\mu_0)_* \Phi)(u))(x) \\
&= \mu^* ( ( \mu_* ( \lambda A) \otimes (\mu_0)_* \Phi)(u))(x) \\
&= \langle u(y), \mu_* ( \lambda A)(\mu(x), y) \cdot v(\mu(x)) \otimes (\mu_0)_*(\vec\varphi(x))(y) \rangle \\
&= \langle (\mu^* u)(y), (\lambda A)(x) \cdot (\mu^* v)(x) \otimes \vec\varphi(x)(y) \rangle \\
&= \iota ( \mu^* u)(A, \Phi)(x).
\end{align*}
The second claim is clear from
\begin{equation*}
 \mu^* ( \sigma t)(A, \Phi)(x) = (\mu^* t) (x) = \sigma(\mu^*t)(A, \Phi)(x).\qedhere
\end{equation*}

\end{proof}

\Autoref{definition.3.4} and \Autoref{theorem.3.5} hold as well in the local case; moreover, the pullback of the local expression equals the local expression of the pullback.
\subsection{Test objects}

The quotient construction, which ensures that the calculus of smooth sections is preserved in spaces of nonlinear generalized sections, is realized by evaluating elements of $\cE(M, \lambda E)$ on certain nets of vector smoothing operators. In light of \eqref{equation.2.6} these will be constructed as a combination of transport operators and scalar test objects.

\subsubsection{Transport operators}

In this section we will introduce and study so-called \emph{admissible nets of transport operators} which will be used for the regularization of vector-valued distributions on manifolds. For what follows let $M, N$ be manifolds, $E \to M$, $F \to N$ vector bundles, $\Omega, \Omega'$ open subsets of $\bR^n$ and $\bE, \bF$ finite-dimensional vector spaces. Moreover, $E_i \to M$ ($i=1 \dotsc k$) will be a family of vector bundles, $\bE_i$ ($i=1 \dotsc k$) a family of finite dimensional vector spaces and $\lambda$ a homogeneous functor in $k$ arguments.

We begin by defining the spaces of interest:

\begin{definition}\label{definition.3.6}
By $\Upsilon(M, E)$ we denote the space of all nets $\bA = (\bA_\e)_\e \in \TO(M,E)^I$ satisfying
\begin{enumerate}[label=(\roman*)]
 \item\label{definition.3.6.1} $\forall x \in M$ $\exists U \in \cU_x(M)$ $\exists \e_0>0$: $\{ \bA_\e|_{U \times U} : \e < \e_0 \}$ is bounded in $\TO(U, E)$. 
 \item\label{definition.3.6.2} $\forall p \in \csn (\Gamma(M, E \otimes E^*))$ $\forall m \in \bN$:
 \[ p(\Delta^* \bA_\e - \id_E) = O(\e^m), \]
where $\Delta \colon M \to M \times M$, $x \mapsto (x,x)$ is the diagonal mapping and $\id_E$ is the identity on each fiber of $E$.
\end{enumerate}
Elements of $\Upsilon(M,E)$ are called \emph{admissible} nets of transport operators.

By $\Upsilon_0(M,E)$ we denote the space of all nets $\bA \in \TO(M,E)^I$ satisfying \ref{definition.3.6.1} and
\begin{enumerate}
 \item[(ii')] $\forall p \in \csn (\Gamma(M, E \otimes E^*))$ $\forall m \in \bN$: $p(\Delta^* \bA_\e ) = O(\e^m)$.
\end{enumerate}
\end{definition}

The local variant of this definition is as follows:

\begin{definition}\label{definition.3.7}
$\Upsilon(\Omega, \bE)$ is the space of all nets $\ba = (\ba_\e)_\e \in \TO(\Omega, \bE)^I$ satisfying
\begin{enumerate}[label=(\roman*)]
 \item \label{definition.3.7.1} $\forall x \in \Omega$ $\exists V \in \cU_x(\Omega)$ $\exists \e_0>0$: $\{\ba_\e|_{V \times V} : \e<\e_0\}$ is bounded in $\TO(V, \bE)$.
 \item \label{definition.3.7.2} $\forall K \csub \Omega$ $\forall l \in \bN_0$ $\forall m \in \bN$:
 \[ \sup_{x \in K} \norm{ (\D^l_x (\ba_\e(x,x) - \id_\bE) } = O(\e^m). \]
\end{enumerate}
Here, $\norm{\cdot}$ denotes any norm on the corresponding space of multilinear mappings. Elements of $\Upsilon(\Omega, \bE)$ are called \emph{admissible} nets of transport operators.

By $\Upsilon_0(\Omega,\bE)$ we denote the space of all nets $\ba \in \TO(\Omega,\bE)^I$ satisfying \ref{definition.3.7.1} and

\begin{enumerate}
 \item[(ii')] \label{Item.17} $\forall K \csub \Omega$ $\forall l \in \bN_0$ $\forall m \in \bN$:
 \[ \sup_{x \in K} \norm{ (\D^l_x (\ba_\e(x,x)) } = O(\e^m). \]
\end{enumerate}
\end{definition}

The spaces $\Upsilon_0(M,E)$ and $\Upsilon_0(\Omega, \bE)$ will be needed only for the \sm.

We will study the following operations on nets of transport operators:

\begin{enumerate}
 \item {\bfseries Restriction.} Clearly, $\TO(\unterstrich,E)^I$ is a presheaf of vector spaces on $M$ with restriction mapping defined by
\begin{equation}\label{equation.3.12}
 \bA|_V \coleq (\bA_\e|_{V \times V})_\e
\end{equation}
and similarly for $\TO(\unterstrich, \bE)$.
 \item {\bfseries Expression in local charts and trivializations.} If $(U, \varphi)$ is a chart on $M$ over which $E$ is trivial with typical fiber $\bE$, then to $\bA \in \TO(U, E)^I$ there corresponds the local expression $\ba = (\ba_\e)_\e \in \TO(\varphi(U),\bE)^I$ and to the coordinates $\bA^i_j \coleq ((\bA_\e)^i_j)_\e \in C^\infty(U \times U)^I$ there correspond the local expressions $\ba^i_j \coleq ((\ba_\e)^i_j)_\e \in C^\infty(\varphi(U) \times \varphi(U))^I$.
 \item {\bfseries Application of homogeneous functors.} Similarly, the action of smooth functors given by \Autoref{definition.2.5} is extended componentwise:
 \begin{align*}
  \lambda (\bA_1, \dotsc, \bA_k) & \coleq (\lambda(\bA_{1\e}, \dotsc, \bA_{k\e}))_\e \\
  \lambda (\ba_1, \dotsc, \ba_k) & \coleq (\lambda(\ba_{1\e}, \dotsc, \ba_{k\e}))_\e.
 \end{align*}
 \item {\bfseries (Local) vector bundle morphisms.} For a vector bundle isomorphism $\mu \colon E \to F$ and a local vector bundle isomorphism $\sigma \colon \Omega \times \bE \to \Omega' \times \bF$ we have componentwise pullback operations $\mu^* \colon \TO(N, F)^I \to \TO(M, E)^I$ and $\sigma^* \colon \TO(\Omega', \bF)^I \to \TO(\Omega, \bE)^I$, respectively.
 \item {\bfseries Lie derivatives.}  On nets of transport operators $\Lie_X$ acts componentwise, and similar for partial derivatives in the local case:
 \begin{align*}
  \Lie_X \bA & \coleq (\Lie_X \bA_\e)_\e && (X \in \fX(M), \bA \in \TO(M,E)^I), \\
  \pd_x^\alpha \pd_y^\beta \ba & \coleq (\pd_x^\alpha \pd_y^\beta \ba_\e)_\e && (\alpha,\beta \in \bN_0^n, \ba \in \TO(\Omega, \bE)^I).
 \end{align*}
\end{enumerate}

First, we summarize how these operations relate to each other on the level of transport operators:

\begin{proposition}\label{proposition.3.8}
 \begin{enumerate}[label=(\roman*)]
  \item\label{alles.1} Given $A_i \in \TO(M, E_i)$ for all $i$ and an open subset $U \subseteq M$, we have $\lambda(A_1, \dotsc, A_k)|_U = \lambda ( A_1|_U, \dotsc, A_k|_U)$.
  \item\label{alles.2} For $A \in \TO(N, F)$ and $U \subseteq M$ open, $\mu^*(A)|_U = \mu^*(A|_{\mu(U)})$ (and similarly in the local case)
  \item\label{alles.3} With $A \in \TO(M, E)$ let $(U,\varphi)$ be a chart of $M$ where $E$ is trivial, $V \subseteq U$ open, $B = A|_U$ and $a$ and $b$ the local representations of $A$ and $B$, respectively. We then have $b = a|_{\varphi(U)}$.
  \item\label{alles.4} For $X \in \fX(M)$, $U \subseteq M$ open and $A \in \TO(M,E)$ we have $(\Lie_X A)|_U = \Lie_{X|_U} A|_U$.
  \item\label{alles.6} Let $(U, \varphi)$ be a chart on $M$ where each $E_i$ is trivial with typical fiber $\bE_i$. Given $A_i \in \TO(U, E_i)$ with local expression $a_i \in \TO(\varphi(U), \bE_i)$, the local expression of $\lambda(A_1, \dotsc, A_k) \in \TO(U, \lambda(E_1, \dotsc, E_k))$ is given by $\lambda(a_1, \dotsc, a_k) \in \TO(\varphi(U), \lambda(\bE_1, \dotsc, \bE_k))$.
  \item\label{alles.8} If $\mu$ is a vector bundle isomorphism and $\sigma$ its local expression, $B = \mu^*A$ for transport operators $A$ and $B$ implies $b = \sigma^* a$ for their local expressions.
\end{enumerate}
\end{proposition}
\begin{proof}
We only show \ref{alles.6}, the rest being evident.

Let $\Phi_i \colon \pi^{-1}(U) \to U \times \bE_i$ be a trivialization of $E_i$ and $\Phi \colon \pi^{-1}(U) \to U \times \lambda(\bE_1, \dotsc, \bE_k)$ the corresponding trivialization of $\lambda(E_1, \dotsc, E_k)$. Suppose that $A \coleq \lambda(A_1, \dotsc, A_k)$ has local expression $a$. Then (in case $\lambda$ is covariant in all arguments) for $x,y \in \varphi(U)$ and $v \in \lambda(\bE_1, \dotsc, \bE_k)$ we have
\begin{align*}
a(x,y) \cdot v &= (\pr_2 \circ \Phi \circ A(\varphi^{-1}(x), \varphi^{-1}(y)) \circ \Phi^{-1})(\varphi^{-1}(y),v) \\
&= ( \Phi_{\varphi^{-1}(x)} \circ A(\varphi^{-1}(x), \varphi^{-1}(y)) \circ \Phi^{-1}_{\varphi^{-1}(y)})(v) \\
&= ( \lambda( (\Phi_1)_{\varphi^{-1}(x)}, \dotsc, (\Phi_k)_{\varphi^{-1}(x)}) \\
& \quad \circ \lambda ( A_1(\varphi^{-1}(x), \varphi^{-1}(y)), \dotsc, A_k(\varphi^{-1}(x), \varphi^{-1}(y))) \\
&\quad \circ \lambda ( (\Phi_1)^{-1}_{\varphi^{-1}(y)}, \dotsc, (\Phi_k)^{-1}_{\varphi^{-1}(y)} ))(v) \\
&= \lambda ( (\Phi_1)_{\varphi^{-1}(x)} \circ A_1(\varphi^{-1}(x), \varphi^{-1}(y)) \circ (\Phi_1)^{-1}_{\varphi^{-1}(y)}, \\
&\quad \dotsc, (\Phi_k)_{\varphi^{-1}(x)} \circ A_k ( \varphi^{-1}(x), \varphi^{-1}(y)) \circ (\Phi_k)^{-1}_{\varphi^{-1}(y)} ) \cdot v \\
&= \lambda ( a_1(x,y), \dotsc, a_k(x,y))\cdot v.
\end{align*}
For contravariant arguments the claim is seen analogously.
\end{proof}

Next, we show that these operations preserve the spaces defined above:

\begin{proposition}\label{proposition.3.9} Let $\bA \in \TO(M,E)^I$. If $\bA \in \Upsilon(M, E)$ then for any open subset $U \subseteq M$, $\bA|_U \in \Upsilon(U, E)$. Conversely, if each point $x \in M$ has an open neighborhood $U \subseteq M$ such that $\bA|_U \in \Upsilon(U,E)$ then $\bA \in \Upsilon(M,E)$. An analogous statement holds for $\Upsilon_0(M, E)$ and, in the local case, for $\Upsilon(\Omega,\bE)$ and $\Upsilon_0(\Omega, \bE)$.
\end{proposition}
\begin{proof}This is clear because the conditions of \Autoref{definition.3.6} and \Autoref{definition.3.7} are local.
\end{proof}

\begin{proposition}\label{proposition.3.10}Let $(U,\varphi)$ be a chart on $M$, $\Phi \colon \pi^{-1}(U) \to U \times \bE$ a trivialization of $E$ over $U$, $\bA \in \TO(U, E)^I$ and $\ba \in \TO(\varphi(U), \bE)^I$ its local expression. Then $\bA \in \Upsilon(U, E)$ if and only if $\ba \in \Upsilon(\varphi(U), \bE)$. Similarly, $\bA \in \Upsilon_0(U, E)$ if and only if $\ba \in \Upsilon_0(\varphi(U), \bE)$.
\end{proposition}

\begin{proof}
Suppose $\bA$ is admissible. To show \Autoref{definition.3.7} \ref{definition.3.7.1}, let $x \in \varphi(U)$. By \Autoref{definition.3.6} (i) there is $U_0 \in \cU_{\varphi^{-1}(x)}(M)$ and $\e_0>0$ such that $\{\bA_\e|_{U_0 \times U_0} \colon \e < \e_0 \}$ is bounded in $\TO(U_0,E)$. Set $V \coleq \varphi(U \cap U_0)$. We know that $\{ \bA_\e|_{\varphi^{-1}(V) \times \varphi^{-1}(V)} : \e < \e_0\}$ is bounded in $\TO(\varphi^{-1}(V), E)$, which implies that $\{ \ba_\e|_{V \times V} : \e < \e_0\}$ is bounded in $\TO(V, \bE)$.

For \Autoref{definition.3.7} \ref{definition.3.7.2}, let the trivialization $\widetilde \Phi$ of $E^*$ be given by \eqref{equation.3.10} and fix $K \csub \varphi(U)$, $l \in \bN_0$ and $m \in \bN$. Then the mapping on $\Gamma(M, E \otimes E^*)$ defined by
\begin{equation}\label{equation.3.13}
 p \colon s \mapsto \sup_{x \in K} \norm{ (\D^l \tilde s)(x)}
\end{equation}
with
\[ \tilde s(x) \coleq (\Phi_{E_{\varphi^{-1}(x)}} \otimes \widetilde \Phi_{E^*_{\varphi^{-1}(x)}}) ( s(\varphi^{-1}(x))) \]
is a continuous seminorm, hence we have
\[ p( s_\e ) = O(\e^m) \]
for $s_\e = \Delta^* \bA_\e - \id_E$ by assumption. This gives the claim because then
\begin{equation}\label{equation.3.14}
 \tilde s_\e(x) = ( \Phi_{E_{\varphi^{-1}(x)}} \otimes \widetilde \Phi_{E^*_{\varphi^{-1}(x)}}) ( \bA_\e(\varphi^{-1}(x), \varphi^{-1}(x)) - \id_E) = \ba_\e(x,x) - \id_\bE.
\end{equation}

Conversely, suppose $\ba$ is admissible. Fix $x \in M$. By \Autoref{definition.3.7} \ref{definition.3.7.1} there is some $V \in \cU_{\varphi(x)}(\varphi(U))$ and $\e_0>0$ such that $\{\ba_\e|_{V \times V} : \e < \e_0 \}$ is bounded in $\TO(\varphi(V), \bE)$, which implies that $\{ \bA_\e|_{\varphi^{-1}(V) \times \varphi^{-1}(V)} : \e<\e_0 \}$ is bounded in $\TO(\varphi^{-1}(V), E)$, i.e., \ref{definition.3.6.1} of \Autoref{definition.3.6}. For \ref{definition.3.6.2}, fix $p \in \csn(\Gamma(M, E \otimes E^*))$. We may assume $p$ to be of the form \eqref{equation.3.13} for some $K \csub \varphi(U)$ and $\tilde s$ defined as above because the family of all such $p$ forms a fundamental system of continuous seminorms of $\Gamma(M, E \otimes E^*)$. Fix $m \in \bN_0$. By \Autoref{definition.3.7} \ref{definition.3.7.2} we know that
\[ \sup_{x \in K} \norm{ (\D^l_x(\ba_\e(x,x) - \id_\bE)} = O(\e^m) \]
which gives the claim because for $s_\e = \Delta^*\bA_\e - \id_E$, \eqref{equation.3.14} holds.

The statement for $\Upsilon_0$ is seen similarly by removing $\id$ from the above.
\end{proof}

The reason why we consider homogeneous functors is that they preserve admissibility.
\begin{proposition}\label{proposition.3.11}
 \begin{enumerate}[label=(\roman*)]
  \item Given $\bA_i \in \Upsilon(M, E_i)$ for all $i$ we have
  \[ \lambda(\bA_1, \dotsc, \bA_k) \in \Upsilon(M, \lambda(E_1, \dotsc, E_k)). \]
  If at least one of the $\bA_i$ is in $\Upsilon_0(M,E_i)$ instead, we have $\lambda(\bA_1, \dotsc, \bA_k) \in \Upsilon_0(M, \lambda(E_1, \dotsc, E_k))$.
 \end{enumerate}
 \begin{enumerate}[label=(\roman*),resume]
  \item\label{proposition.3.11.3} Given $\ba_i \in \Upsilon(\Omega, \bE)$ for all $i$ we have
  \[ \lambda(\ba_1, \dotsc, \ba_k) \in \Upsilon(\Omega, \lambda(\bE_1, \dotsc, \bE_k)). \]
  If at least one of the $\ba_i$ is in $\Upsilon_0(\Omega, \bE)$ instead, we have $\lambda(\ba_1, \dotsc, \ba_k) \in \Upsilon_0(\Omega, \lambda(\bE_1, \dotsc, \bE_k))$.
 \end{enumerate}
\end{proposition}

\begin{proof}
By \Autoref{proposition.3.9}, \Autoref{proposition.3.10} and \Autoref{proposition.3.8} \ref{alles.1}, \ref{alles.6} it suffices to show \ref{proposition.3.11.3}. For simplicity we assume that $\lambda$ is covariant in all arguments; the general case is obtained by replacing $\ba_i(x,y)$ by $\ba_i(y,x)$ where necessary.

First, $\ba \coleq \lambda(\ba_1, \dotsc, \ba_k)$ is given by
\begin{equation}\label{equation.3.15}
 \ba_\e(x,y) = h( \underbrace{\ba_{1\e}(x,y), \dotsc}_{n_1\textrm{ times}}, \underbrace{ \ba_{2\e}(x,y), \dotsc}_{n_2\textrm{ times}}, \dotsc, \underbrace{\ba_{k\e}(x,y), \dotsc}_{n_k\textrm{ times}})
\end{equation}
for some multilinear function
\[ h \colon \Lin(\bE_1, \bE_1)^{n_1} \times \dotsc \times \Lin(\bE_k, \bE_k)^{n_k} \to \Lin(\lambda(\bE_1, \dotsc, \bE_k),\lambda(\bE_1, \dotsc, \bE_k)). \]
To verify \Autoref{definition.3.7} \ref{definition.3.7.1} we apply the chain rule to \eqref{equation.3.15}.
By multilinearity of $h$, $\D^l\ba_\e(x,y)$ is given by terms of the form
\begin{align*}
h(& D^{l^1_1} \ba_{1\e}(x,y), \D^{l^2_1}\ba_{1\e}(x,y), \dotsc, \D^{l^1_2} \ba_{2\e}(x,y), \D^{l^2_2}\ba_{2\e}(x,y), \dotsc, \\
& \dotsc, \D^{l_k^1} \ba_{k\e}(x,y), \D^{l^2_k}\ba_{k\e}(x,y), \dotsc )
\end{align*}
for some $0 \le l^j_i \le l$ ($i = 1 \dotsc k$, $j = 1 \dotsc n_i$) and the norm of such a term may be estimated by products of derivatives of $\ba_{i\e}(x,y)$, which by assumption satisfy the desired conditions of \Autoref{definition.3.7} \ref{definition.3.7.1}.

For \Autoref{definition.3.7} \ref{definition.3.7.2} we have to estimate derivatives of
\begin{multline*}
  \ba_\e(x,x) - \id = h(\ba_{1\e}(x,x) , \dotsc, \ba_{2\e}(x,x), \dotsc, \ba_{k\e}(x,x), \dotsc) \\
  - h(\id_{\bE_1}, \dotsc, \id_{\bE_2}, \dotsc, \id_{\bE_k}, \dotsc).
\end{multline*}
By multilinearity of $h$, this can be written as
\begin{gather*}
h(\ba_{1\e}(x,x) - \id_{\bE_1}, \dotsc) + h(\dotsc, \ba_{2\e}(x,x) - \id_{\bE_2}, \dotsc) \\
+ \dotsc + h(\dotsc, \ba_{k\e}(x,x) - \id_{\bE_k}).
\end{gather*}
As above, the norm of each of these terms may be estimated by products of derivatives of $\ba_{i\e}(x,x)$ and of $\ba_{i\e}(x,x) - \id_{\bE_i}$, where the last factor has to appear at least once. By assumption, this factor satisfies the $O(\e^m)$-estimate of \Autoref{definition.3.7} \ref{definition.3.7.2} and the other factors are uniformly bounded for small $\e$, which completes the first part of the proof of \ref{proposition.3.11.3}. For the second part one puts zero instead of $\id$ and $\id_{\bE_i}$ in the above and proceeds similarly.
\end{proof}

For use in the \sm, we describe differentials of homogeneous functors applied to transport operators:

\begin{proposition}\label{proposition.3.12}
For $\bA_1, \dotsc, \bA_k$ in $\Upsilon(M,E) \cup \Upsilon_0(M,E)$ and $\bA_0 \in \Upsilon_0(M,E)$ we have
\[  (\ud_i \lambda)(\bA_1, \dotsc, \bA_k)(\bA_0) \in \Upsilon_0(M,E), \]
and similarly in the local case.
\end{proposition}
\begin{proof}
 This follows from \Autoref{proposition.3.11} using homogeneity of $\lambda$.
\end{proof}

Admissible transport operators are well-behaved under bundle isomorphisms:

\begin{proposition}\label{proposition.3.13}Let $\mu \colon E \to F$ be a smooth bundle isomorphism.
\begin{enumerate}[label=(\roman*)]
 \item\label{proposition.3.12.1} For $\bA \in \Upsilon(N, F)$, $\mu^* \bA \in \Upsilon(M, E)$.
 \item\label{proposition.3.12.2} For $\bA \in \Upsilon_0(N, F)$, $\mu^* \bA \in \Upsilon_0(M, E)$.
\end{enumerate}
In case of a local vector bundle isomorphism $\mu \colon \Omega \times \bE \to \Omega' \times \bF$, we have
\begin{enumerate}[label=(\roman*),resume]
 \item\label{proposition.3.12.3} For $\ba \in \Upsilon(\Omega', \bF$), $\mu^* \ba \in \Upsilon(\Omega, \bE)$.
 \item\label{proposition.3.12.4} For $\ba \in \Upsilon_0(\Omega', \bF$), $\mu^* \ba \in \Upsilon_0(\Omega, \bE)$.
\end{enumerate}
\end{proposition}
\begin{proof}
Again, we reduce \ref{proposition.3.12.1} and \ref{proposition.3.12.2} to the local case. We have
\begin{align*}
 \mu^* A \in \Upsilon(M, E) & \Longleftrightarrow (\mu^* \bA)|_U = \mu^*(\bA|_{\mu(U)}) \in \Upsilon(U,E) \textrm{ for all }U\\
 & \Longleftrightarrow \varphi_*(\mu^*(A|_{\mu(U)})) = \mu^*(\varphi_*(A|_{\mu_0(U)})) \in \Upsilon(\Omega, \bE)
\end{align*}
by \Autoref{proposition.3.9}, \Autoref{proposition.3.8} \ref{alles.2} and \ref{alles.8} and \Autoref{proposition.3.10}. Moreover, $\varphi_*(A|_{\mu(U)})$ is in $\Upsilon(\Omega, \bF)$ by \Autoref{proposition.3.9} and \Autoref{proposition.3.10}, so \ref{proposition.3.12.3} implies \ref{proposition.3.12.1}; similarly, \ref{proposition.3.12.4} implies \ref{proposition.3.12.2}.

For \ref{proposition.3.12.3} and \ref{proposition.3.12.4}, if $\mu$ is given by $\mu(x,v) = (\mu_1(x), \mu_2(x) v)$ then $\mu^* \ba$ is given by
\[ (\mu^* \ba)_\e(x,y) = \mu_2(x)^{-1} \circ \ba_\e(\mu_1(x), \mu_1(y)) \circ \mu_2(y) \]
and the desired conditions of \Autoref{definition.3.7} hold by the chain rule.
\end{proof}

Concerning the structure of the above spaces, $\Upsilon_0(M, E)$ and $\Upsilon_0(\Omega, \bE)$ are vector spaces and $\Upsilon(M,E)$ and $\Upsilon(\Omega,\bE)$ are affine spaces over $\Upsilon_0(M,E)$ and $\Upsilon_0(\Omega, \bE)$, respectively. More generally, we can say:

\begin{lemma}\label{lemma.3.14}
For $\chi \in C^\infty(M \times M)$ the following holds:
\begin{enumerate}[label=(\roman*)]
 \item\label{lemma.3.20.1}  $\forall \bA,\bB \in \Upsilon(M,E)$: $\chi\bA + (1-\chi)\bB \in \Upsilon(M,E)$.
 \item\label{lemma.3.20.2} $\forall \bA \in \Upsilon_0(M,E)$: $\chi \bA \in \Upsilon_0(M,E)$.
 \end{enumerate}
 Similarly, for $\chi \in C^\infty(\Omega \times \Omega)$ the following holds:
\begin{enumerate}[label=(\roman*),resume]
 \item\label{lemma.3.20.3} $\forall \ba,\bb \in \Upsilon(\Omega,\bE)$: $\chi\ba + (1-\chi)\bb \in \Upsilon(\Omega,\bE)$.
 \item\label{lemma.3.20.4} $\forall \ba \in \Upsilon_0(\Omega,\bE)$: $\chi \ba \in \Upsilon_0(\Omega,\bE)$.
 \end{enumerate}
\end{lemma}
\begin{proof}
By \Autoref{proposition.3.9} and \Autoref{proposition.3.10}, \ref{lemma.3.20.1} and \ref{lemma.3.20.2} are reduced to showing \ref{lemma.3.20.3} and \ref{lemma.3.20.4}, respectively. The latter, however, are evident from the definition.
\end{proof}

\begin{proposition}
For $\bA \in \Upsilon(M, E) \cup \Upsilon_0(M, E)$ and $X \in \fX(M)$, $\Lie_X \bA \in \Upsilon_0(M,E)$. 
\end{proposition}
\begin{proof}
The claim follows directly from the definitions using
\[
 p ( \Delta^* ( \Lie_X \bA_\e ) ) = p ( \Lie_X ( \Delta^* \bA_\e - \id_E ))
\]
for $p \in \csn ( \Gamma(M, E \otimes E^*))$, bacuse $p \circ \Lie_X$ is a continuous seminorm as well.
\end{proof}

Trivially, for any (nonvoid) open subset $\Omega \subseteq \bR^n$, $\Upsilon(\Omega, \bE)$ is non-empty as is seen by simply defining $\ba_\e(x,y) \coleq \id_\bE$. Setting $\ba_\e(x,y) \coleq 0$ gives an element of $\Upsilon_0(\Omega, \bE)$. By \Autoref{proposition.3.10} also $\Upsilon(U,E)$ and $\Upsilon_0(U,E)$ are non-empty whenever $U$ is a chart domain on $M$ where $E$ is trivial.

To obtain a sheaf structure for nets of transport operators we introduce the following notion of transport operators being locally eventually equal.

\begin{definition}\label{definition.3.16}
For any open set $W \subseteq M$ we define a relation $\sim_W$ on the union
\[ \bigcup \{ \TO(X,E)^I\ |\ X \textrm{ open},\ W \subseteq X \subseteq M \} \]
as follows: for open sets $U, V \subseteq M$ with $W \subseteq U \cap V$, $\bA \in \TO(U, E)^I$ and $\bB \in \TO(V,E)^I$ we write ``$\bA \sim_W \bB$'' or ``$\bA \sim \bB$ on $W$'' if
\[ \forall x \in W\ \exists Z \in \cU_x(U \cap V) \ \exists \e_0>0\ \forall \e<\e_0: \bA_\e|_{Z \times Z} = \bB_\e|_{Z \times  Z}. \]
\end{definition}

If $W = U \cap V$ we simply write ``$\bA \sim \bB$''.

Clearly, $\sim_U$ defines a congruence relation on $\TO(U,E)^I$, so we can form the quotient vector space $\TO(U,E)/{\sim_U}$. Moreover, $\bA \sim 0$ implies $\bA|_V \sim 0$, which turns $U \mapsto \TO(U,E)^I/{\sim_U}$ into a presheaf of vector spaces, as the restriction on the quotient does not depend on the representative.

\begin{lemma} Let $\bA, \bB \in \TO(M,E)^I$ be such that $\bA \sim \bB$. Then the following holds:
\begin{enumerate}[label=(\roman*)]
 \item $\bA \in \Upsilon(M, E)$ $\Longleftrightarrow$ $\bB \in \Upsilon(M,E)$.
 \item $\bA \in \Upsilon_0(M, E)$ $\Longleftrightarrow$ $\bB \in \Upsilon_0(M,E)$.
\end{enumerate}
\end{lemma}
\begin{proof}
This follows immediately from \Autoref{proposition.3.9}. In fact, for $\bA \in \Upsilon(M,E)$ with $\bA \sim \bB$, given any $x \in M$ there is $Z \subseteq M$ open and $\e_0>0$ such that $\bA_\e|_{Z \times Z} = \bB_\e|_{Z \times Z}$, which gives $\bB|_Z \in \Upsilon(Z,E)$ and similarly for the case of $\bA \in \Upsilon_0(M,E)$.
\end{proof}

\begin{proposition} $U \mapsto \TO(U, E)^I / {\sim_U}$ is a sheaf of vector spaces on $M$.
\end{proposition}
\begin{proof}
Let $U \subseteq M$ be open and $(U_i)_i$ an open cover of $U$. Suppose we are given $\bA \in \TO(U,E)^I$ with $\bA|_{U_i} \sim 0$ for all $i$. For showing that $\bA \sim 0$, let $x \in U$. Then $x$ is in some $U_i$ and $\exists Z \in \cU_x(U_i)$, $\e_0>0$ such that $\bA_\e|_{U_i \times U_i}|_{Z \times Z} = \bA_\e|_{Z \times Z} = 0$ for $\e<\e_0$, which gives the claim.

Next, suppose that $A_i \in \TO(U_i, E)^I$ are given with $\bA_i \sim \bA_j$ on $U_i \cap U_j$. Choose a family of functions $\chi_j \in C^\infty(U \times U)$ such that $\supp \chi_j \subseteq U_j \times U_j$, $0 \le \chi_j \le 1$, $(\supp \chi_j)_j$ is locally finite and $\sum \chi_j = 1$ in a neighborhood of the diagonal of $U \times U$.

With this, set $\bA_\e(x,y) \coleq \sum_j \chi_j(x,y) \bA_{j\e} (x,y)$ for $x,y \in U$. We claim that $\bA|_{U_i} \sim \bA_i$. Fix $x \in U_i$. Then there are a neighborhood $W$ of $x$ and a finite index set $F$ such that
\[ (W \times W) \cap \supp \chi_j \ne \emptyset \Rightarrow j \in F. \]
Moreover, we may assume that $x \in U_j$ for all $j \in F$. By assumption, $\forall j \in F$ $\exists W_j \in \cU_x(U_i \cap U_j)$ $\exists \e_j>0$ $\forall \e<\e_j$: $\bA_{j\e}|_{W_j \times W_j} = \bA_{i\e}|_{W_j \times W_j}$. Furthermore, there is a neighborhood $V$ of $x$ in $U$ such that $\sum \chi_j(x,y) = 1$ for all $x,y \in V$. Set $W = V \cap \bigcap_{j \in F}W_j$ and $\e_0 \coleq \min_j \e_j$. Then for $\e<\e_0$ and $x,y \in W$,
\begin{align*}
 A|_{U_i \times U_i}(x,y) &= \sum_{j \in F} \chi_j(x,y) \bA_{j\e}(x,y) \\
 & = \sum_{j \in F}\chi_j(x,y) \bA_{i\e}(x,y) = \bA_{i\e}(x,y). \qedhere
\end{align*}
\end{proof}

\begin{corollary}\label{corollary.3.19}$U \mapsto \Upsilon_0(U,E)/{\sim_U}$ is a sheaf of vector spaces on $M$ and $U \mapsto \Upsilon(U,E)/{\sim_U}$ is a sheaf of affine spaces over $U \mapsto \Upsilon_0(U, E)/{\sim_U}$.
\end{corollary}

The following technical result will be useful later on for the proof of the sheaf property of generalized functions.

\begin{corollary}\label{corollary.3.20} Let $U,V,W \subseteq M$ be open such that $\overline{W} \subseteq U \cap V \ne \emptyset$.
\begin{enumerate}[label=(\roman*)]
 \item For $\bA \in \Upsilon(V,E)$ there exists $\bA' \in \Upsilon(U,E)$ such that $\bA|_W \sim \bA'|_W$.
 \item For $\bA \in \Upsilon_0(V,E)$ there exists $\bA' \in \Upsilon_0(U,E)$ such that $\bA|_W \sim \bA'|_W$.
\end{enumerate}
\end{corollary}

\begin{proof}Choose an open neighborhood $X$ of $\overline{W}$ such that $\overline{X} \subseteq U \cap V$, any $\bA^0 \in \Upsilon(U,E)$ (or $\Upsilon_0(U,E)$) and $\chi \in C^\infty(U \cap V)$ with $\supp \chi \subseteq X$ and $\chi \equiv 1$ on $W$. Then
\[ \bB \coleq \chi \cdot \bA|_{U \cap V} + (1-\chi) \bA^0|_{U \cap V} \in \Upsilon(U \cap V, E) \quad \textrm{(or }\Upsilon_0(U \cap V, E)\textrm{)} \]
by \Autoref{lemma.3.14}. We have
\[ \bB|_{(U \cap V) \cap (U \setminus \overline{X})} = \bA^0|_{(U \cap V)\cap(U \setminus \overline{X})} \]
so by \Autoref{corollary.3.19} there exists $\bA' \in \Upsilon(U,E)$ (or $\Upsilon_0(U, E)$) such that
\[ \bA'|_W \sim \bB|_W = \bA|_W. \qedhere \]
\end{proof}

The above statements about the sheaf structure (\Autoref{definition.3.16} to \Autoref{corollary.3.20} are valid also in the local case, with identical proofs.

For later use we note the following immediate consequence of \Autoref{definition.3.6}.

\begin{lemma}\label{lemma.3.21}
Let $\bA \in \Upsilon(\Omega,\bE)$, $C>0$ and $K \csub \Omega$. Then with respect to any basis of $\bE$,
\[ \sup_{x \in K, y \in B_{C\e}(x)} \abso{ {\bA_\e}^i_j(x,y) - \delta^i_j } \to 0. \]
\end{lemma}
As $A^*(x,y) \coleq A(y,x)^*$ is in $\Upsilon(\Omega,\bE)$ the same result holds with $\bA^i_j(y,x)$ in place of $A^i_j(x,y)$.

\subsubsection{Scalar Test objects}

We will now study scalar test objects; their combination with admissible nets of transport operators will give vector test objects later on. In this section let $M,N$ be manifolds and $\Omega, \Omega'$ open subsets of $\bR^n$.

For the next definition, we denote by $\cP_k(M)$ the space of linear differential operators $C^\infty(M) \to C^\infty(M)$ of order $k \in \bN_0$ \cite[Chapter 6]{Kahn}. Due to $C^\infty(M, \bE) = C^\infty(M) \widehat\otimes \bE$, each $P \in \cP_k(M)$ induces a continuous linear mapping $P \colon C^\infty(M, \bE) \to C^\infty(M, \bE)$ for any complete Hausdorff locally convex space $\bE$, given by the continuous extension of the mapping $f \otimes e \mapsto P(f) \otimes e$ defined on the dense subspace $C^\infty(M) \otimes \bE$. Below, we set $\bE = \Gamma_c(M, \Vol(M))$ and write $P(\pd_x)\vec\varphi(x)$ for $\vec\varphi \in C^\infty(M, \Gamma_c(M, \Vol(M)))$ to make clear that $P$ acts on the $x$-variable of $\vec\varphi$.

Note that the absolute value $\abso{\mu}$ of a continuous density $\mu \colon M \to \Vol(M)$, defined by
\[ \abso{\mu}(p)(v_1, \dotsc, v_n) \coleq \abso{ \mu(p)(v_1, \dotsc, v_n) } \qquad (p \in M,\ v_1, \dotsc, v_n \in T^*_p M) \]
is a continuous density again.

\begin{definition}\label{definition.3.22}
Let $\bPhi = (\bPhi_\e)_\e \in \SO(M)^I$ correspond via \eqref{equation.2.5} to $\bvarphi =(\bvarphi_\e)_\e \in \SK(M)^I$. $\bPhi$ is said to be a \emph{test object} on $M$ if the following conditions hold:
\begin{enumerate}[label=(\roman*)]
 \item\label{def_scaltestobj.1} $\forall g \in \Riem(M)$ $\forall K \csub M$ $\exists C,\e_0>0$ $\forall x \in K$ $\forall \e<\e_0$: $\supp \bvarphi_\e(x) \subseteq B_{C\e}^g(x)$.
 \item\label{def_scaltestobj.2} $\bPhi_\e \to \id$ in $\Lin ( \cD'(M), \cD'(M) )$.
 \item\label{def_scaltestobj.3} $\forall p \in \csn ( \SO(M) )$ $\exists N \in \bN$: $p ( \bPhi_\e ) = O(\e^{-N})$.
 \item\label{def_scaltestobj.4} $\forall p \in \csn ( \Lin( C^\infty(M), C^\infty(M) ) )$ $\forall m \in \bN$: $p ( \bPhi_\e|_{C^\infty(M)} - \id_{C^\infty(M)} ) = O(\e^m)$.
 \item\label{def_scaltestobj.5} $\forall K \csub M$ $\forall k \in \bN_0$ $\forall P \in \cP_k(M)$: $\sup_{x \in K} \int \abso{P(\pd_x) \bvarphi_\e(x) } = O(\e^{-k})$.
 \item\label{def_scaltestobj.6} $\forall K \csub M$: $\sup_{x \in K} \int \abso{\bvarphi_\e(x)} \to 1$.
\end{enumerate}
By $S(M)$ we denote the space of test objects on $M$.

$\bPhi$ is said to be a \emph{zero test object} on $M$ if it satisfies \ref{def_scaltestobj.1}, \ref{def_scaltestobj.3}, \ref{def_scaltestobj.5} as well as
\begin{enumerate}[label=(\roman*')]
 \item[(ii')] $\bPhi_\e \to 0$ in $\Lin ( \cD'(M), \cD'(M) )$.
 \item[(iv')] $\forall p \in \csn ( \Lin( C^\infty(M), C^\infty(M) ) )$ $\forall m \in \bN$: $p ( \bPhi_\e|_{C^\infty(M)} ) = O(\e^m)$.
 \item[(vi')] $\forall K \csub M$: $\sup_{x \in K} \int \abso{\bvarphi_\e(x)} \to 0$.
\end{enumerate}
By $S_0(M)$ we denote the space of all zero test objects on $M$.

By $S_1(M)$ we denote the space of all $\bPhi \in \Lin ( \cD'(M), C^\infty(M) )^I$ satisfying \ref{def_scaltestobj.1}.
\end{definition}

Note that while conditions \ref{def_scaltestobj.1}--\ref{def_scaltestobj.4} are classical in the functional analytic approach to Colombeau algebras and constitute the core components necessary for the construction of an algebra of nonlinear generalized functions, conditions \ref{def_scaltestobj.5} and \ref{def_scaltestobj.6} are included to ensure good association properties. These latter conditions are at the root of many good properties which convolution by scaled mollifiers has in the local case (cf.~\cite{GKOS}), and here we have its equivalent on manifolds (cf.~\cite{SV:09} for the role of condition \ref{def_scaltestobj.6} and also \cite[Def.\ 7.1]{MOBook}, where a weaker version of this condition is used for defining model delta nets). 

On open subsets of $\bR^n$, the definition of test objects looks as follows:

\begin{definition}\label{definition.3.23}
Let $\bPhi = (\bPhi_\e)_\e \in \SO(\Omega)^I$ correspond to $\bvarphi =(\bvarphi_\e)_\e \in \SK(\Omega)^I$ via \eqref{equation.2.4}. $\bPhi$ is said to be a \emph{test object} on $\Omega$ if the following conditions hold:
\begin{enumerate}[label=(\roman*)]
 \item\label{def_scaltestobjloc.1} $\forall K \csub \Omega$ $\exists C,\e_0>0$ $\forall x \in K$ $\forall \e<\e_0$: $\supp \bvarphi_\e(x) \subseteq B_{C\e}(x)$.
 \item\label{def_scaltestobjloc.2} $\bPhi_\e \to \id$ in $\Lin ( \cD'(\Omega), \cD'(\Omega) )$.
 \item\label{def_scaltestobjloc.3} $\forall p \in \csn ( \SO(\Omega) )$ $\exists N \in \bN$: $p ( \bPhi_\e ) = O(\e^{-N})$.
 \item\label{def_scaltestobjloc.4} $\forall p \in \csn ( \Lin( C^\infty(\Omega), C^\infty(\Omega) ) )$ $\forall m \in \bN$: $p ( \bPhi_\e|_{C^\infty(\Omega)} - \id_{C^\infty(\Omega)} ) = O(\e^m)$.
 \item\label{def_scaltestobjloc.5} $\forall K \csub \Omega$ $\forall \alpha \in \bN_0^n$: $\sup_{x \in K}\int \abso{\pd_x^\alpha \bvarphi_\e(x) }\,\ud x = O(\e^{-\abso{\alpha}})$.
 \item\label{def_scaltestobjloc.6} $\forall K \csub \Omega$: $\sup_{x \in K} \int\abso{\bvarphi_\e(x)(y)}\,\ud y \to 1$.
\end{enumerate}
By $S(\Omega)$ we denote the space of test objects on $\Omega$.

$\bPhi$ is said to be a \emph{zero test object} on $\Omega$ if it satisfies \ref{def_scaltestobjloc.1}, \ref{def_scaltestobjloc.3}, \ref{def_scaltestobjloc.5} as well as
\begin{enumerate}[label=(\roman*')]
 \item[(ii')] $\bPhi_\e \to 0$ in $\Lin ( \cD'(\Omega), \cD'(\Omega) )$.
 \item[(iv')] $\forall p \in \csn ( \Lin( C^\infty(\Omega), C^\infty(\Omega) ) )$ $\forall m \in \bN$: $p ( \bPhi_\e|_{C^\infty(\Omega)}) = O(\e^m)$.
 \item[(vi')] $\forall K \csub \Omega$: $\sup_{x \in K} \int\abso{\bvarphi_\e(x)(y)}\,\ud y \to 0$.
\end{enumerate}
By $S_0(\Omega)$ we denote the space of all zero test objects on $\Omega$.

By $S_1(\Omega)$ we denote the space of all $\bPhi \in \Lin ( \cD'(\Omega), C^\infty(\Omega) )^I$ satisfying \ref{def_scaltestobjloc.1}.
\end{definition}

Given a chart $(U,\varphi)$ on $M$, a smoothing operator $\Phi \in \SO(U)$ has local expression $\varphi_*(\Phi) \coleq \varphi_* \circ \Phi \circ \varphi^* \in \SO(\varphi(U))$. With this, \Autoref{definition.3.22} and \Autoref{definition.3.23} are compatible:
\begin{proposition}
 Let $(U,\varphi)$ be a chart on $M$, $\bPhi \in \SO(U)^I$ and let its local expression be given by $\bPsi \in \SO(\varphi(U))^I$. Then
 \begin{enumerate}[label=(\roman*)]
  \item $\bPhi \in S_1(U)$ $\Longleftrightarrow$ $\bPsi \in S_1(\varphi(U))$,
  \item $\bPhi \in S(U)$ $\Longleftrightarrow$ $\bPsi \in S(\varphi(U))$,
  \item $\bPhi \in S_0(U)$ $\Longleftrightarrow$ $\bPsi \in S_0(\varphi(U))$.
 \end{enumerate}
\end{proposition}
The proof is immediate from the definitions (see \cite{bigone} for more details).

Defining the pullback of nets of smoothing operators componentwise, we can state diffeomorphism invariance of the space of test objects:

\begin{proposition}\label{proposition.3.25}
If $\mu \colon M \to N$ is a diffeomorphism then $\mu^* \bPhi \in S(M)$ for $\bPhi \in S(N)$, $\mu^* \bPhi \in S_1(M)$ for $\bPhi \in S_1(N)$ and $\mu^* \bPhi \in S_0(M)$ for $\bPhi \in S_0(N)$.

Similarly, if $\mu \colon \Omega \to \Omega'$ is a diffeomorphism then $\mu^* \bPhi \in S(\Omega)$ for $\bPhi \in S(\Omega')$, $\mu^* \bPhi \in S_1(\Omega)$ for $\bPhi \in S_1(\Omega')$ and $\mu^* \bPhi \in S_0(\Omega)$ for $\bPhi \in S_0(\Omega')$.
\end{proposition}

Again, this is easily seen from the definitions.

We now recall the main properties related to localization of test objects from \cite{specfull}. All proofs from there are transferred with ease to our setting (cf.~also \cite{bigone}) and skipped here.
\begin{proposition}\label{proposition.3.26}
\begin{enumerate}[label=(\roman*)]
 \item For any open set $W \subseteq M$ we define a relation $\sim_W$ on the union $\bigcup \{\, \SK(X)^I\ |\ X\textrm{ open}, W \subseteq X \subseteq M \, \}$ as follows: given two open sets $U,V \subseteq M$ with $W \subseteq U \cap V$, $\bvarphi \in \SK(U)^I$ and $\bpsi \in \SK(V)^I$ we say ``$\bvarphi \sim_W \bpsi$'' or ``$\bvarphi \sim \bpsi$ on $W$'' if $\forall x \in W\ \exists Z \in \cU_x ( U \cap V )\ \exists \e_0>0\ \forall \e<\e_0: \bvarphi_\e|_Z = \bpsi_\e|_Z$. If $W = U \cap V$ we simply write ``$\bvarphi \sim \bpsi$''.
 \item Let $U,V$ be open subsets of $M$ with $V \subseteq U$ and $\rho_V \in C^\infty(V, C^\infty_c(V))$ be equal to 1 in a neighborhood of the diagonal in $V \times V$. Then the mapping $\rho_V \colon \SK(U) \to \SK(V)$, $\vec\varphi \mapsto \rho_V \cdot \vec\varphi|_V$ is linear and continuous. Applying it componentwise to elements $\bvarphi, \bpsi$ of $S^1(U)$, we obtain:
\begin{enumerate*}[label=(\roman*)]
 \item $\rho_V \bvarphi \sim \bvarphi $,
 \item $\rho_V \bvarphi  \in S^1(V)$,
 \item $\bvarphi \sim \bpsi$ implies $\rho_V \bvarphi  \sim \rho_V \bpsi$,
 \item for any open subset $W \subseteq V$ we have $\rho_W(\rho_V \bvarphi )) \sim \rho_W \bvarphi $,
 \item given any other $\rho_V' \in C^\infty(V, C^\infty_c(V))$ which is equal to 1 on an open neighborhood of the diagonal in $V \times V$, we have $\rho_V \bvarphi  \sim \rho_V' \bvarphi $.
\end{enumerate*}
\item $U \mapsto S^1(U) / {\sim_U}$ is a sheaf of $C^\infty$-modules on $M$. Denoting the equivalence class of $\bvarphi$ in the quotient by $[\bvarphi]$, the restriction map from $U$ to $V$ is given by $S^1(U) / {\sim_U} \to S^1(V) / {\sim_V}$, $[\bvarphi] \mapsto [\rho_V \bvarphi]$ and does not depend on the choice of $\rho_V$.
\item Let $(U_i)_i$ be a family of open subsets of $M$ and $U = \bigcup_i U_i$. If an element $\bvarphi \in S^1(U)$ satisfies any of the conditions
of Definition \ref{definition.3.22} and $\bvarphi \sim \bpsi \in S^1(U)$ then $\bpsi$ also satisfies that condition. Moreover, $\bvarphi$ satisfies any of these conditions if and only if for each $i$, some representative of $[\bvarphi]|_{U_i}$ satisfies it.
\item\label{proposition.3.26.4} $U \mapsto S(U)/{\sim_U}$ and $U \mapsto S_0(U)/{\sim_U}$ are sheaves of sets on $M$.
\item\label{Item.49} Let $U,V,W$ be open sets such that $\overline{W} \subseteq V \cap U \ne \emptyset$, and let $\bvarphi \in S(V)$ (or $S_0(V)$). Then there exists $\bpsi \in S(U)$ (or $S_0(U)$) such that $[\bvarphi]|_W = [\bpsi]|_W$.
\end{enumerate}
\end{proposition}

The statements of \Autoref{proposition.3.26} are also valid in the local case.

\subsubsection{Vector test objects}

In this subsection let again $M, N$ be manifolds, $E \to M$, $F \to N$ vector bundles, $\Omega, \Omega'$ open subsets of $\bR^n$ and $\bE, \bF$ finite-dimensional vector spaces.

Convergence of a net $\bPhi \in \Lin(\cD'(M), C^\infty(M))^I$ to the identity on $\cD'(M)$ is equivalent to
\[ \langle \langle u(y), \bvarphi_\e(x)(y) \rangle, \omega(x) \rangle \to \langle u(x), \omega(x) \rangle \]
for all $u \in \cD'(M)$ and $\omega \in \Gamma_c(M, \Vol(M))$.

More generally, as is seen from \eqref{equation.2.7} we will also need to know the behavior of
\begin{equation}\label{equation.3.16}
\langle \langle u(y), f_\e(x,y) \bvarphi_\e(x)(y) \rangle, \omega(x) \rangle
\end{equation}
for $(f_\e)_\e \in C^\infty(M \times M)^I$. Supposing that $f_\e \to f$ we expect \eqref{equation.3.16} to converge to $\langle u(x) f(x,x), \omega(x)\rangle$. This will be established by the next result, which allows us to consider $\SO(M)$ as a locally convex $C^\infty(M \times M)$-module such that \eqref{equation.3.16} corresponds to applying the smoothing operator $f_\e \bullet \bPhi_\e$, where the module action $\bullet$ is defined as follows:

\begin{lemma}\label{lemma.3.27}
 Suppose we are given two locally convex spaces $\bE,\bF$ and two nuclear Fr\'echet spaces $X,Y$ together with hypocontinuous bilinear mappings\begin{align*}
X \times \bE & \xrightarrow{\cdot} \bE, \\
Y \times \bF & \xrightarrow{\cdot} \bF.
\end{align*}
 Then there exists a unique separately continuous mapping
\[ (X \widehat\otimes Y) \times \Lin ( \bF, \bE) \xrightarrow{\bullet} \Lin(\bF, \bE) \]
such that $ ( ( g \otimes h ) \bullet T ) (x) = g \cdot T ( h \cdot x)$ for all $g \in X$, $h \in Y$, $T \in \Lin(\bF, \bE)$ and $x \in \bF$. Moreover, $\bullet$ is hypocontinuous.
\end{lemma}
\begin{proof}
Consider the trilinear mapping
\[ \varphi \colon C^\infty(M) \times C^\infty(M) \times \Lin(\bF, \bE) \to \Lin(\bF, \bE) \]
defined by $\varphi(g,h,T)(x) \coleq g \cdot T(h \cdot x)$. It is partially continuous in $g$ and $h$ and
$\varphi(.,.,T)$ induces a unique linear continuous mapping $\varphi_T \colon (X \widehat\otimes Y) \to \Lin(\bF,\bE)$ such that $\varphi_T(g \otimes h) = \varphi(g,h,T)$ for all $g \in X$, $h \in Y$ and $T \in \Lin(\bF,\bE)$. We define $f \bullet T \coleq \varphi_T(f)$. To show continuity of this mapping in $T$, let $f$ be fixed. By \cite[Theorem 21.5.8, p.~495]{Jarchow} we can choose bounded sets $B_1 \subseteq X$, $B_2 \subseteq Y$ such that $f \in \cacx ( B_1 \otimes B_2)$. Let $U$ be an absolutely convex closed $0$-neighborhood in $\bE$ and $B \subseteq \bF$ be bounded. Choose a $0$-neighborhood $U'$ in $\bE$ such that $B_1 \cdot U' \subseteq U$ and set $W = \{ T \in \Lin(\bF, \bE): T(B_2 \cdot B) \subseteq U' \}$. Then
\begin{equation}\label{equation.3.17}
(( B_1 \otimes B_2) \bullet W)(B) = \varphi(B_1, B_2, W)(B) \subseteq U
\end{equation}
and hence
\begin{gather*}
(f \bullet W)(B) \in ( \cacx ( B_1 \otimes B_2) \bullet W)(B) \\
\subseteq \cacx ( ( ( B_1 \otimes B_2) \bullet W)(B)) \subseteq \cacx U = U
\end{gather*}
Uniqueness follows from \cite[\S 40.3(1), p.~162]{Koethe2} because $X \otimes Y$ is dense in $X \widehat\otimes Y$.

Because $X \widehat\otimes Y$ is a nuclear Fr\'echet space and hence barrelled we already know that $\bullet$ is hypocontinuous with respect to bounded subsets of $\Lin(\bF,\bE)$. Moreover, the above argument for continuity in $T$ is also valid for $f$ in a bounded subset $D \subseteq X \widehat\otimes Y$ because one can then find bounded sets $B_1$ and $B_2$ such that $D \subseteq \cacx(B_1 \otimes B_2)$.
\end{proof}

We will use this lemma in the case where $X=Y=C^\infty(M)$ and hence $X \widehat\otimes Y = C^\infty(M \times M)$.

\begin{lemma}\label{lemma.3.28}Let $f \in C^\infty(M \times M)$.
\begin{enumerate}[label=(\roman*)]
\item \label{corollary.3.33.1} 
With $\bullet$ the product
\[ \bullet \colon C^\infty(M_x \times M_y) \times \Lin(\cD'(M_y), C^\infty(M_x)) \to \Lin(\cD'(M_y), C^\infty(M_x)), \]
for $\Phi \in \SO(M)$, $u \in \cD'(M)$ and $\omega \in \Gamma_c(M, \Vol(M))$ we have
 \[ \langle (f \bullet \Phi)(u), \omega \rangle = \langle \langle u(y), f(x,y) \vec\varphi(x)(y) \rangle, \omega(x) \rangle. \]
 \item
With $\bullet$ the product
\[ \bullet \colon C^\infty(M_x \times M_y) \times \Lin(\cD'(M_y), \cD'(M_x)) \to \Lin(\cD'(M_y), \cD'(M_x)) \]
we have 
$(f \bullet \id_{\cD'(M)})(u) = (\Delta^* f)\cdot u$, where $(\Delta^* f)(x) = f(x,x)$.
\end{enumerate}
\end{lemma}
\begin{proof}
 The equations hold for $f = g \otimes h$ with $g,h \in C^\infty(M)$, hence by continuity also for general $f \in C^\infty(M \times M)$.\
\end{proof}

Finally, we can combine admissible nets of transport operators and scalar test objects as in \eqref{equation.2.7}. Their main properties are captured by the following result.

\begin{theorem}\label{theorem.3.29}
 If $\bA \in \Upsilon(M, E)$ and $\bPhi \in S(M)$ then $\bTheta \coleq \bA \otimes \bPhi$ satisfies the following conditions:
 \begin{enumerate}[label=(\roman*)]
 \item\label{def_testobj.1} $\forall g \in \Riem(M)$ $\forall K \csub M$ $\exists C,\e_0>0$ $\forall x \in K$ $\forall \e<\e_0$ $\forall u \in \cD'(M,E)$:
 \[ u = 0 \textrm{ on }B^g_{C\e}(x) \Longrightarrow \bTheta_\e(u)(x) = 0. \]
 \item\label{def_testobj.2} $\bTheta_\e \to \id$ in $\Lin ( \cD'(M, E), \cD'(M, E) )$.
 \item\label{def_testobj.3} $\forall p \in \csn ( \VSO(M, E))$ $\exists N \in \bN$: $p ( \bTheta_\e ) = O(\e^{-N})$.
 \item\label{def_testobj.4} $\forall p \in \csn ( \Lin ( \Gamma(M,E), \Gamma(M,E) ) )$ $\forall m \in \bN$: $p ( \bTheta_\e|_{\Gamma(M,E)} - \id_{\Gamma(M,E)} ) = O(\e^m)$.
\end{enumerate}
For $\bA \in \Upsilon_0(M,E)$ or $\bPhi \in S_0(M,E)$ instead the same statements hold with $\id$ replaced by $0$ in \ref{def_testobj.2} and \ref{def_testobj.4}.
\end{theorem}

Such objects are called ``vector test objects'' in \cite{bigone}. Note that we will not need to formulate analogons of \Autoref{definition.3.22} \ref{def_scaltestobj.5} and \ref{def_scaltestobj.6} here, as these will only be used in direct calculations (see Section \ref{section.6}).

\begin{proof}
\ref{def_testobj.1} Let $g \in \Riem(M)$ and $K \csub M$. By \Autoref{definition.3.22} \ref{def_scaltestobj.1} there are $C_0,\e_0 > 0$ such that
\[ \supp \bvarphi_\e(x) \subseteq B^g_{C\e}(x)\qquad \forall x \in K, \e < \e_0. \]
Let $x \in K$, $\e<\e_0$ and $u \in \cD'(M,E)$ with $u=0$ on $B^g_{C\e}(x)$ be given. Then for $v \in \Gamma(M, E^*)$,
\begin{align*}
( \bTheta_\e(u) \cdot v)(x) &= \langle u(y), \bA_\e(x,y) v(x) \otimes \bvarphi_\e(x)(y) \rangle \\
&= \langle u(y) \bA_\e(x,y) v(x), \bvarphi_\e(x)(y) \rangle = 0 
\end{align*}
because $u \bA_\e(x,.)v(x) = 0$ on $B^g_{C\e}$.

\ref{def_testobj.2} For preciseness, let us introduce a useful notation: for any sheaf of $C^\infty$-modules $\sF$ on $M$ and any open subset $U \subseteq M$ there is a unique ``extension by zero'' mapping
\[ \ext \colon \{ u \in \sF(U)\ |\ \supp u\textrm{ is closed in }M \} \to \sF(M) \]
such that $\ext(u)|_U = u$ and $\supp(\ext(u)) = \supp u$. We will use this in the cases where $\sF$ is one of $\Gamma(\unterstrich, E^* \otimes \Vol(M))$, $\Gamma(\unterstrich, E^*)$, $\Gamma(\unterstrich, \Vol(M))$ or $\cD'(\unterstrich)$.

Fix $x \in M$ and choose an open neighborhood $U$ of $x$ and $\e_0>0$ according to \Autoref{definition.3.6} \ref{definition.3.6.1} such that $E$ is trivial over $U$ and $\{\bA_\e|_{U \times U}\ :\ \e < \e_0 \}$ is bounded in $\TO(U,E)$. Let $v_1, \dotsc, v_m$ be a basis of $\Gamma(U,E)$ and $\alpha^1, \dotsc, \alpha^m$ the basis of $\Gamma(U, E^*)$ dual to it. We denote by
\begin{align*}
 \bTheta_\e(u)^i &\coleq \bTheta_\e(u)|_U \cdot \alpha^i \in \cD'(U), \\
 {\bA_\e}^i_j & \coleq [ (x,y) \mapsto \bA_\e(x,y) \cdot \alpha^i(x) v_j(y) ] \in C^\infty(U \times U), \\
 u^j &\coleq u|_U \cdot \alpha^j \in \cD'(U)
\end{align*}
the coordinates of $\bTheta_\e(u)$, $\bA_\e$ and $u \in \cD'(M)$ with respect to these bases. We will show that
\[ \langle \bTheta_\e(u)^i, \omega \rangle \to \langle u^i, \omega \rangle \qquad \forall u \in \cD'(M)\ \forall \omega \in \Gamma_c(U, \Vol(M)), \]
which will establish \ref{def_testobj.2}.

Given $u$ and $\omega$ let $\rho \in C^\infty_c(U)$ be equal to 1 on $\supp \omega$ and choose a compact set $K \csub U$ and $\e_0' \le \e_0$ such that $\supp \bvarphi_\e(x) \subseteq K$ for all $x \in \supp \omega$ and all $\e<\e_0'$. Let $\chi \in C^\infty_c(U)$ be equal to 1 on $K$. Then for such $\e$, we have
\begin{align}
 \nonumber \langle \bTheta_\e(u)^i, \omega \rangle & = \langle \bTheta_\e(u)|_U \cdot \alpha^i, \omega \rangle \\
\nonumber  &= \langle \bTheta_\e(u)|_U \cdot (\rho \alpha^i), \omega \rangle \\
 \nonumber & = \langle \bTheta_\e(u)|_U \cdot \ext(\rho \alpha^i)|_U, \omega \rangle \\
 \nonumber &= \langle ( \bTheta_\e(u) \cdot \ext(\rho \alpha^i))|_U, \omega \rangle \\
\label{equation.3.18} &= \langle \bTheta_\e(u) \cdot \ext(\rho \alpha^i), \ext(\omega) \rangle.
\end{align}
Applying \eqref{equation.2.7} to \eqref{equation.3.18} gives
\begin{gather}
\nonumber \langle \langle u(y), \bA_\e(x,y) \ext(\rho \alpha^i)(x) \otimes \bvarphi_\e(x)(y)\rangle, \ext(\omega)(x) \rangle \\
\label{equation.3.19} = \langle \langle u(y), \bA_\e|_{U \times M}(x,y) \rho(x) \alpha^i(x) \otimes \bvarphi_\e|_U(x)(y) \rangle, \omega(x) \rangle.
\end{gather}
For each $x \in \supp \omega$ we have $\supp \bA_\e|_{U \times M}(x,.) \alpha^i(x) \otimes \bvarphi_\e|_U (x)(.) \subseteq K$ and hence
\begin{align*}
 & \langle u(y), \bA_\e|_{U \times M}(x,y)\rho(x) \alpha^i(x) \otimes \bvarphi_\e|_U(x)(y) \rangle \\
 = {} & \langle u|_U(y), \bA_\e|_{U \times U}(x,y) \alpha^i(x) \otimes \bvarphi_\e|_U(x)|_U(y) \rangle \\
 = {} & \langle u^j(y), \bA_\e|_{U \times U}(x,y) \alpha^i(x) v_j(y) \bvarphi_\e|_U(x)|_U(y) \rangle \\
 = {} & \langle u^j(y), {\bA_\e}^i_j(x,y) \bvarphi_\e|_U(x)|_U(y) \rangle \\
 = {} & \langle \chi(y) u^j(y), \rho(x)\chi(y) {\bA_\e}^i_j(x,y) \bvarphi_\e|_U(x)|_U(y) \rangle \\
 = {} & \langle \ext ( \chi u^j)(y), \ext ( \rho(x) \chi {\bA_\e}^i_j(x,.)\bvarphi_\e|_U(x)|_U)(y) \rangle \\
 = {} & \langle \ext ( \chi u^j)(y), \ext ( \rho(x) \chi {\bA_\e}^i_j(x,.)) \bvarphi_\e|_U(x)(y) \rangle,
\end{align*}
so \eqref{equation.3.19} equals
\[ \langle \langle \ext  ( \chi u^j)(y), \ext(\rho(x) \chi {\bA_\e}^i_j(x,.))\bvarphi_\e|_U(x)(y) \rangle, \omega(x) \rangle. \]
To get rid of the restriction to $U$ in the $x$-variable we write this as
\begin{equation}\label{equation.3.20}
\langle \ext ( x \mapsto \langle \ext ( \chi u^j)(y), \rho(x) \ext ( \chi {\bA_\e}^i_j(x,.))\bvarphi_\e|_U(x)(y) \rangle ), \ext ( \omega)(x) \rangle.
 \end{equation}
We now use the fact that for $v \in \cD'(M)$ and $\vec\psi \in C^\infty(U, \Gamma_c(M, \Vol(M)))$ with $\supp \vec\psi$ closed in $M$,
\[ \ext ( x \mapsto \langle v(y), \vec\psi(x)(y) \rangle )(x) = \langle u(y), \ext ( \vec\psi )(x)(y) \rangle, \]
so \eqref{equation.3.20} becomes
\begin{equation}\label{equation.3.21}
\langle \langle \ext ( \chi u^j)(y), \ext ( x \mapsto \rho(x) \ext ( \chi {\bA_\e}^i_j(x,.)))(x) \bvarphi_\e(x)(y) \rangle, \ext(\omega)(x) \rangle
\end{equation}
and because
\[ \ext ( x \mapsto \rho(x) \ext ( \chi {\bA_\e}^i_j(x,.))) = \ext ( (\rho \otimes \chi)  {\bA_\e}^i_j) \]
\eqref{equation.3.21} equals, by \Autoref{lemma.3.28} \ref{corollary.3.33.1},
\begin{align*}
 & \langle \langle \ext ( \chi u^j)(y), \ext ( ( \rho \otimes \chi) {\bA_\e}^i_j)(x,y) \bvarphi_\e(x)(y) \rangle, \ext(\omega)(x) \rangle \\
 = {} & \langle ( \ext ( ( \rho \otimes \chi) {\bA_\e}^i_j) \bullet \bPhi_\e) ( \ext(\chi u^j)), \ext(\omega) \rangle.
\end{align*}

On the other hand,
\begin{align*}
 \langle u^i, \omega \rangle & = \langle \chi u^i, \omega \rangle = \langle \chi u^j, \delta^i_j \omega \rangle \\
 & = \langle \ext ( \chi u^j), \ext ( \delta^i_j \omega) \rangle = \langle \ext ( \chi u^j) \delta^i_j, \ext(\omega) \rangle \\
 &= \langle ( \delta^i_j \bullet \id_{\cD'(M)} ) ( \ext ( \chi u^j)), \ext(\omega) \rangle.
\end{align*}
Hence, we can write $\langle \bTheta_\e (u)^i, \omega \rangle - \langle u^i, \omega \rangle$ as
\begin{align}
 \nonumber & \langle ( \ext ( ( \rho \otimes \chi) {\bA_\e}^i_j) \bullet \bPhi_\e - \delta^i_j \bullet \id_{\cD'(M)}) ( \ext(\chi u^j)), \ext(\omega) \rangle \\
 \label{equation.3.22} = {} & \langle ( \ext ( ( \rho \otimes \chi ) {\bA_\e}^i_j) \bullet ( \bPhi_\e - \id_{\cD'(M)}))(\ext(\chi u^j)), \ext(\omega) \rangle \\
 \label{equation.3.23} + {} & \langle ( ( \ext ( ( \rho \otimes \chi ) {\bA_\e}^i_j) - \delta^i_j) \bullet \id_{\cD'(M)} ) ( \ext(\chi u^j)), \ext(\omega) \rangle.
\end{align}
Because $\{ \ext ( ( \rho \otimes \chi) {\bA_\e}^i_j) : \e < \e_0 \}$ is bounded in $C^\infty(M \times M)$ and $\bPhi_\e \to \id_{\cD'(M)}$ in $\Lin(\cD'(M), \cD'(M))$, hypocontinuity of $\bullet$ implies that \eqref{equation.3.22} converges to zero. To see that also \eqref{equation.3.23} converges to zero, we rewrite it as
\begin{align*}
 & \langle ( \ext ( (  \rho \otimes \chi) {\bA_\e}^i_j )(x,x) - \delta^i_j) \ext ( \chi u^j), \ext(\omega) \rangle \\
 = {} & \langle ( {\bA_\e}^i_j(x,x) - \delta^i_j) u^j(x), \omega(x) \rangle
\end{align*}
and use the fact that $\Delta^* {\bA_\e}^i_j - \delta^i_j \to 0$ in $C^\infty(U)$ by \Autoref{definition.3.6} \ref{definition.3.6.2}.

In case $\bA \in \Upsilon_0(M,E)$ or $\bPhi \in S_0(M)$ one similarly sees that
\[ \langle ( \ext (( \rho \otimes \chi) {\bA_\e}^i_j ) \bullet \bPhi_\e) ( \ext ( \chi u^j)), \ext(\omega) \rangle \]
converges to $0$.

\ref{def_testobj.3} Let $p \in \csn(\VSO(M,E))$ be given. We can assume that $p$ has the form
\[ p(\Theta) = \sup_{u \in B, i=1 \dotsc m} p'(\Theta(u)^i) \]
for some bounded subset $B \subseteq \cD'(M,E)$ and $p' \in \csn(C^\infty(U))$, where $U$ is an open subset of $M$ where $E$ is trivial. We furthermore may assume $p'$ to be of the form
\begin{equation}\label{equation.3.24}
p'(f) = \sup_{x \in K} \abso{ P(\pd_x) f(x)}
\end{equation}
where $K \csub U$ and $P$ is a linear differential operator on $C^\infty(U)$. Let $\rho \in C^\infty_c(U)$ be $1$ on a neighborhood $V$ of $K$. We have, for $x \in V$,
\begin{gather}
 \bTheta_\e(u)^i(x)
 \label{equation.3.25} = \langle u(y), \bA_\e(x,y) \ext(\rho \alpha^i)(x) \otimes \bvarphi_\e(x)(y) \rangle.
\end{gather}
Choose $\e_0>0$ and a compact neighborhood $L$ of $V$ with $V \subseteq U$ such that for $\e< \e_0$, $\supp \bvarphi_\e (x) \subseteq L$. Let $\chi \in C^\infty_c(U)$ be equal to 1 on $L$. Then as in \ref{def_testobj.2}, \eqref{equation.3.25} equals
\begin{gather*}
 \langle u|_U(y), \bA_\e|_{U \times U} (x,y) \rho(x) \alpha^i(x) \otimes \bvarphi_\e|_U (x)|_U(y) \rangle \\
 = \langle u^j(y), {\bA_\e}^i_j (x,y) \bvarphi_\e|_U(x)|_U(y) \rangle \\
 = ( \ext ( ( \rho \otimes \chi) {\bA_\e}^i_j ) \bullet \bPhi_\e) ( \ext(\chi u^j))(x).
\end{gather*}
Hence, we have
\[ P(\pd_x) ( \bTheta_\e(u)^i(x)) = P(\pd_x) ( \langle u^j(y), {\bA_\e}^i_j(x,y) \bvarphi_\e|_U(x)|_U(y) \rangle ) \]
which gives
\begin{gather*}
 p(\bTheta_\e) \le C \sup_{u \in B, i,j = 1 \dotsc m} p' ( ( \ext ( ( \rho \otimes \chi) {\bA_\e}^i_j) \bullet \bPhi_\e)(\ext(\chi u^j)))
\end{gather*}
for some $C>0$ where we know that $\{ \ext( (\rho \otimes \chi) {\bA_\e}^i_j )\ |\ \e < \e_0 \}$ is bounded; by hypocontinuity the last expression can be estimated by
\[ \sup_{u \in B} p'' ( \bPhi_\e(\ext(\chi u^j))) \]
for some $p'' \in \csn ( \SO(M) )$, and as $\{ \ext ( \chi u^j)\ |\ u \in B, j=1 \dotsc m \}$ is bounded in $\cD'(M)$, this grows like $O(\e^{-N})$ for some $N \in \bN$ by \Autoref{definition.3.22} \ref{def_scaltestobj.3}.

For \ref{def_testobj.4} one proceeds similarly as above. In fact, suppose we are given a seminorm $p \in \csn ( \Lin ( \Gamma(M, E), \Gamma(M, E)))$ and $m \in \bN$. We may assume $p$ to be of the form
\[ p(\Theta) = \sup_{u \in B, i=1 \dotsc m} p'(\Theta(u)^i) \]
with $B \subseteq \Gamma(M,E)$ bounded and $p' \in \csn(C^\infty(U))$ for some open subset $U \subseteq M$ where $E$ is trivial. Moreover, $p'$ may be again assumed to be of the form \eqref{equation.3.24} with $K \csub U$ and a linear differential operator $P$. With $\rho$, $\chi$ and $\e<\e_0$ as above, we have for $x$ in a neighborhood of $K$ and $u \in B$
\begin{gather*}
 \bTheta_\e(u)^i (x) - u^i(x) = ( \ext ( ( \rho \otimes \chi){\bA_\e}^i_j) \bullet ( \bPhi_\e - \id_{C^\infty(M)} )) ( \ext ( \chi u^j))(x) \\
 + (( \ext ( ( \rho \otimes \chi) {\bA_\e}^i_j) - \delta^i_j) \bullet \id_{C^\infty(M)} ) ( \ext ( \chi u^j))(x)
\end{gather*}
so $p(\Theta_\e - \id_{C^\infty(M)})$ can be estimated by
\begin{gather*}
 \sup_{u,i} p' ( ( \ext ( ( \rho \otimes \chi){\bA_\e}^i_j) \bullet ( \bPhi_\e - \id_{C^\infty(M)}))(\ext(\chi u^j))) \\
 + \sup_{u,i} p' ( ( (\ext ( ( \rho \otimes \chi){\bA_\e}^i_j) - \delta^i_j) \bullet \id_{C^\infty(M)} ) ( \ext ( \chi u^j)))
\end{gather*}
which satisfies the desired negligibility estimates, and similarly for $\bA \in \Upsilon_0(M, E)$ or $\bPhi \in S_0(M,E)$.
\end{proof}

Similarly, in the local case we have

\begin{proposition}\label{proposition.3.30}
 If $\bA \in \Upsilon(\Omega, \bE)$ and $\bPhi \in S(\Omega)$ then $\bTheta \coleq \bA \otimes \bPhi$ satisfies the following conditions:
 \begin{enumerate}[label=(\roman*)]
 \item\label{def_testobj.1a} $\forall K \csub \Omega$ $\exists C,\e_0>0$ $\forall x \in K$ $\forall \e<\e_0$ $\forall u \in \cD'(\Omega,\bE)$:
 \[ (u = 0 \textrm{ on }B_{C\e}(x) \Longrightarrow \bTheta_\e(u)(x) = 0). \]
 \item\label{def_testobj.2a} $\bTheta_\e \to \id$ in $\Lin ( \cD'(\Omega, \bE), \cD'(\Omega, \bE) )$.
 \item\label{def_testobj.3a} $\forall p \in \csn ( \VSO(\Omega, \bE))$ $\exists N \in \bN$: $p ( \bTheta_\e ) = O(\e^{-N})$.
 \item\label{def_testobj.4a} $\forall p \in \csn ( \Lin ( C^\infty(\Omega,\bE), C^\infty(\Omega, \bE) ) )$ $\forall m \in \bN$: $p ( \bTheta_\e|_{C^\infty(\Omega,\bE)} - \id_{C^\infty(\Omega, \bE)} ) = O(\e^m)$.
\end{enumerate}
\end{proposition}

The proof is similar to that of \Autoref{theorem.3.29}.

We will also need a combination of \Autoref{corollary.3.20} and \Autoref{proposition.3.26} \ref{proposition.3.26.4} for pairs $(\bA, \bvarphi)$.

\begin{lemma}\label{lemma.3.31}Let $U,V,W \subseteq M$ be open such that $\overline{W} \subseteq U \cap V \ne \emptyset$ is compact. Given $\bA \in \Upsilon(V,E)$, $\bA_1, \dotsc, \bA_k \in \Upsilon_0(V,E)$ and $\bvarphi \in S(V)$, $\bvarphi_1, \dotsc, \bvarphi_l \in S_0(U)$ there exist $\bA' \in \Upsilon(U,E)$, $\bA'_1, \dotsc, \bA'_k \in \Upsilon_0(U,E)$ and $\bvarphi' \in S(U)$, $\bvarphi_1', \dotsc, \bvarphi_l' \in S_0(U)$ as well as $\e_0 \in I$ such that for all $\e<\e_0$ we have
 \[ \bvarphi_\e = \bvarphi_\e', \bvarphi_{i\e} = \bvarphi'_{i\e} \textrm{ on }W\ \forall i=1 \dotsc l \]
and
\[ \bA_\e(x,y) = \bA_\e'(x,y), \bA_{j\e}(x,y) = \bA'_{j\e}(x,y)\ \forall j=1 \dotsc k \]
for all $x \in W$ and $y \in \carr \bvarphi_\e(x) \cup \carr \bvarphi_{1\e}(x) \cup \dotsc \cup \carr \bvarphi_{l\e}(x)$.
\end{lemma}
\begin{proof}
Let $X$ be an open neighborhood of $\overline{W}$ whose closure is compact in $U \cap V$. By \Autoref{proposition.3.26} \ref{Item.49} there are $\bvarphi' \in S(U)$, $\bvarphi_i' \in S_0(U)$ such that $\bvarphi \sim_X \bvarphi'$, $\bvarphi_i \sim_X \bvarphi_i'$ $\forall i$. This implies that
\[ \forall p \in X\ \exists W_p \in \cU_p(X)\ \exists \e_p>0\ \forall \e<\e_p: \bvarphi_\e = \bvarphi'_\e, \bvarphi_{i\e} = \bvarphi_{i\e}'\ \forall i \textrm{ on }W_p. \]
Similarly, by \Autoref{corollary.3.20} there are $\bA' \in \Upsilon(U,E)$,$\bA_j' \in \Upsilon_0(U,E)$ such that
\begin{multline*}
\forall p \in X\ \exists W_p' \in \cU_p(X)\ \exists \e_p'>0\ \forall \e<\e_p':\\
\bA_\e|_{W_p' \times W_p'} = \bA'_\e|_{W_p' \times W_p'}, \bA_{j\e}|_{W_p' \times W_p'} = \bA'_{j\e}|_{W_p' \times W_p'}\ \forall j.
\end{multline*}
By Definition \ref{definition.3.22} \ref{def_scaltestobj.1} we have
\begin{multline*}
\forall p \in V\ \exists W_p'' \in \cU_p(W_p \cap W_p')\ \exists \e_p''>0\ \forall \e<\e_p''\ \forall x \in W_p'':\\
\supp \bvarphi_\e(x) \subseteq W_p',\ \supp \bvarphi_{i\e}(x) \subseteq W_p'\ \forall i.
\end{multline*}

Now a finite family $W_{p_1}'', \dotsc, W_{p_k}''$ covers $W$; let $\e_0$ be the minimum of the $\e_{p_1}$, $\dotsc$, $\e_{p_k}$, $\e_{p_1'}$, $\dotsc$, $\e_{p_r'}$ and take $\e<\e_0$ and $x \in W$. Then $x$ is contained in some $W_{p_l}''$ and $\bvarphi_\e(x) = \bvarphi_\e'(x)$, $\bvarphi_{i\e} = \bvarphi'_{i\e}(x)$ $\forall i$. Now $\supp \bvarphi_\e(x) \subseteq W_{p_l}'$, $\supp \bvarphi_{i\e}(x) \subseteq W_{p_l}'$, hence for $y \in \carr \bvarphi_\e(x) \cup \carr \bvarphi_{1\e}(x) \cup \dotsc \cup \bvarphi_{l\e}(x)$ we have $\bA_\e(x,y) = \bA_\e'(x,y)$, $\bA_{j\e}(x,y) = \bA_{j\e}'(x,y)$ because $x \in W''_{p_l} \subseteq W'_{p_l}$ and $y \in W_{p_l}'$.
\end{proof}

\Autoref{lemma.3.31} is also valid in the local case.

\subsection{The quotient construction}

In this subsection let $M$ be a manifold, $E \to M$ a vector bundle, $\Omega$ an open subset of $\bR^n$, $\bE$ a finite dimensional vector space and $\lambda$ a homogeneous functor.

\begin{definition}\label{definition.3.32}We call $R \in \cE(M, \lambda E)$ \emph{moderate} if, in the \si, 
\begin{gather*}
\forall p \in \csn ( \Gamma(M, \lambda E))\ \forall \bA \in \Upsilon(M,E)\\
\forall \bPhi \in S(M)\ \exists N \in \bN: p ( R ( \bA_\e, \bPhi_\e ) )) = O(\e^{-N})
\end{gather*}
or, in the \sm,
\begin{gather*}
\forall k,l \in \bN_0\ \forall p \in \csn ( \Gamma(M, \lambda E))\\
\forall \bA \in \Upsilon(M,E), \bA_1, \dotsc, \bA_k \in \Upsilon_0(M,E)\\\
\forall \bPhi \in S(M), \bPhi_1, \dotsc, \bPhi_l \in S_0(M)\ \exists N \in \bN:\\
p ( (\ud_1^k \ud_2^l R) ( \bA_\e, \bPhi_\e )(\bA_{1\e}, \dotsc, \bA_{k\e}, \bPhi_{1\e}, \dotsc, \bPhi_{l\e}) ) = O(\e^{-N}).
\end{gather*}
The subspace of moderate elements of $\cE(M, \lambda E)$ is denoted by $\cE_M(M, \lambda E)$ and we set $\cE_M(M) \coleq \cE_M(M \times \bR)$.

We call $R \in \cE(M, \lambda E)$ \emph{negligible} if, in the \si,
\begin{gather*}
\forall p \in \csn ( \Gamma(M, \lambda E))\ \forall \bA \in \Upsilon(M, E)\\
\forall \bPhi \in S(M)\ \forall m \in \bN: p ( R ( \bA_\e, \bPhi_\e ) ) = O(\e^m))
\end{gather*}
or, in the \sm,
\begin{gather*}
 \forall k,l \in \bN_0\ \forall p \in \csn ( \Gamma(M, \lambda E))\\
 \forall \bA \in \Upsilon(M,E), \bA_1, \dotsc, \bA_k \in \Upsilon_0(M,E) \\\
\forall \bPhi \in S(M), \bPhi_1, \dotsc, \bPhi_l \in S_0(M)\ \forall m \in \bN:\\
p ( (\ud_1^k \ud_2^l R) ( \bA_\e, \bPhi_\e )(\bA_{1\e}, \dotsc, \bA_{k\e}, \bPhi_{1\e}, \dotsc, \bPhi_{l\e}) ) = O(\e^m).
\end{gather*}
The subspace of negligible elements of $\cE(M, \lambda E)$ is denoted by $\cN(M, \lambda E)$ and we set $\cN(M) \coleq \cN(M \times \bR)$.
\end{definition}

The main properties of the embeddings $\iota$ and $\sigma$ are as follows.

\begin{theorem}\label{thm.3.33} Let $u \in \cD'(M, \lambda E)$ and $t \in \Gamma(M, \lambda E)$. Then the following holds:
 \begin{enumerate}[label=(\roman*)]
  \item\label{thm.3.33.1} $\iota u$ is moderate.
  \item\label{thm.3.33.2} $\sigma t$ is moderate.
  \item\label{thm.3.33.3} $(\iota - \sigma)(t)$ is negligible.
  \item\label{thm.3.33.4} If $\iota(u)$ is negligible then $u=0$.
 \end{enumerate}
\end{theorem}

These conditions form the backbone of every Colombeau-type algebra of generalized functions: \ref{thm.3.33.1} and \ref{thm.3.33.2} ensure that embedded distributions and smooth functions are moderate; \ref{thm.3.33.3} implies that in the quotient, the embeddings $\iota$ and $\sigma$ agree on smooth objects, which is key to preserving the product of smooth functions. \ref{thm.3.33.4} finally gives injectivity of the embeddings.

\begin{proof}
(i) Let $u \in \cD'(M, \lambda E)$. Fix $p \in \csn ( \Gamma(M, \lambda E))$, $\bA \in \Upsilon(M,E)$ and $\bPhi \in S(M)$. As $\Theta \mapsto p(\Theta(u))$ is a continuous seminorm on $\Lin ( \cD'(M, \lambda E), \Gamma(M, \lambda E))$, by \Autoref{theorem.3.29} \ref{def_testobj.3} there exists $N \in \bN$ such that
\[ p ( (\lambda \bA_\e \otimes \bPhi_\e)(u)) = p ( (\iota u)(\bA_\e, \bPhi_\e)) = O(\e^{-N}). \]

For the \sm\ we also have to consider differentials of $\iota u$. We have
\begin{align*}
 (\ud_1 (\iota u))(A, \Phi)(A_0) &= ((\ud \lambda)(A)(A_0) \otimes \Phi)(u) \\
 (\ud_2 ( \iota u))(A,\Phi)(\Phi_0) & = (\lambda A \otimes \Phi_0)(u)
\end{align*}
and similarly for higher derivatives. Using \Autoref{proposition.3.12} and \Autoref{theorem.3.29} the claim follows.

(ii) is clear, as $\sigma t$ does not depend on $A$ or $\Phi$.

(iii) Fix $p \in \csn ( \Gamma(M, \lambda E))$, $\bA$, $\bPhi$ and $m$. Because $\Theta \mapsto p ( \Theta (t))$ is a continuous seminorm on $\Lin ( \Gamma(M, \lambda E), \Gamma(M, \lambda E))$ we have
\[ p ( (\iota t)(\bA_\e, \bPhi_\e) - ( \sigma t)(\bA_\e, \bPhi_\e)) = p ( ( \lambda \bA_\e \otimes \bPhi_\e - \id)(t)) = O(\e^m) \]
by \Autoref{theorem.3.29} \ref{def_testobj.4}

(iv) Fix $\bA$ and $\bPhi$. Then by negligibility, for $u \in \cD'(M, \lambda E)$ we have
\[ (\iota u)(\bA_\e, \bPhi_\e) = ( \lambda \bA_\e \otimes \bPhi_\e)(u) \to 0 \]
in $\Gamma(M, \lambda E)$ and hence in $\cD'(M, \lambda E)$. By \Autoref{theorem.3.29} \ref{def_testobj.2} we also have
\[ ( \lambda \bA_\e \otimes \bPhi_\e)(u) \to u \quad \textrm{in }\cD'(M, \lambda E) \]
and hence $u=0$.
\end{proof}

As is typical in Colombeau spaces of generalized functions, if a function is known to be moderate its negligibility can be tested for without examining the derivatives:

\begin{proposition}\label{proposition.3.34}
 Suppose $R \in \cE(M, \lambda E)$ is moderate; then it is negligible if and only if $\forall K \csub M$ $\forall \bA \in \Upsilon(M, E)$ $\forall \bPhi \in S(M)$:
 \[ \sup_{p \in K} \norm{ R ( \bA_\e, \bPhi_\e )(p) } = O(\e^m) \]
 where $\norm{.}$ is the norm on $\lambda E$ induced by any Riemannian metric on $\lambda E$. 
\end{proposition}

This is proven by similar methods as in \cite[Theorem 6.3, p.~204]{bigone}. Clearly, $\cE(M, \lambda E)$ is a $C^\infty(M)$-module by setting $(f \cdot R)(A, \Phi) \coleq f \cdot R(A,\Phi)$ for $f \in C^\infty(M)$. Furthermore, $\cE_M(M, \lambda E)$ is a submodule of $\cE(M, \lambda E)$ and $\cN(M, \lambda E)$ is a submodule of $\cE_M(M, \lambda E)$, hence we can define the quotient module:

\begin{definition}$\cG(M, \lambda E) \coleq \cE_M(M, \lambda E) / \cN ( M, \lambda E)$ is called the space of nonlinear generalized sections of $\lambda E$.
\end{definition}

As in \cite[Theorem 6.5]{bigone} one sees that the isomorphisms \eqref{equation.3.9} preserve moderateness and negligibility:
\begin{lemma}\label{lemma.3.36}We have the following $C^\infty(M)$-module isomorphisms:
\begin{align*}
 \cE_M(M, \lambda E) & \cong \Gamma(M, \lambda E)\otimes_{C^\infty(M)} \cE_M(M) \\
 & \cong \Hom_{C^\infty(M)} ( \Gamma(M, (\lambda E)^*), \cE_M(M)), \\
 \cN(M, \lambda E) & \cong \Gamma(M, \lambda E)\otimes_{C^\infty(M)} \cN(M) \\
 & \cong \Hom_{C^\infty(M)} ( \Gamma(M, (\lambda E)^*), \cN(M)), \\
 \cG(M, \lambda E) & \cong \Gamma(M, \lambda E)\otimes_{C^\infty(M)} \cG(M) \\
 & \cong \Hom_{C^\infty(M)} ( \Gamma(M, (\lambda E)^*), \cG(M)).
 \end{align*}
\end{lemma}

In case $\lambda E$ is trivial this allows us to talk of the \emph{coordinates} of the elements of $\cE(M, \lambda E)$. Moreover, in this case a generalized section is moderate (negligible) if all its coordinates are moderate (negligible).

Similarly to \Autoref{definition.3.32}, we define the quotient in the local theory as follows:

\begin{definition}We call $R \in \cE(\Omega, \lambda \bE)$ \emph{moderate} if, in the \si,
\begin{gather*}
\forall p \in \csn ( C^\infty(\Omega, \lambda \bE)\ \forall \bA \in \Upsilon(\Omega,\bE)\\
\forall \bPhi \in S(\Omega)\ \exists N \in \bN: p ( R ( \bA_\e, \bPhi_\e ) ) = O(\e^{-N}).
\end{gather*}
or, in the \sm,
\begin{gather*}
  \forall k,l \in \bN_0\ \forall p \in \csn ( C^\infty(\Omega, \lambda \bE))\ \forall \bA \in \Upsilon(\Omega,\bE), \bA_1, \dotsc, \bA_k \in \Upsilon_0(\Omega, \bE) \\\
\forall \bPhi \in S(\Omega), \bPhi_1, \dotsc, \bPhi_l \in S_0(\Omega)\ \forall m \in \bN:\\
p ( (\ud_1^k \ud_2^l R) ( \bA_\e, \bPhi_\e )(\bA_{1\e}, \dotsc, \bA_{k\e}, \bPhi_{1\e}, \dotsc, \bPhi_{l\e}) ) = O(\e^m).
\end{gather*}
The subspace of moderate elements of $\cE(\Omega, \lambda \bE)$ is denoted by $\cE_M(\Omega, \lambda \bE)$ and we set $\cE_M(\Omega) \coleq \cE_M(\Omega \times \bR)$.

We call $R \in \cE(\Omega, \lambda \bE)$ \emph{negligible} if, in the \si,
\begin{gather*}
\forall p \in \csn ( C^\infty(\Omega, \lambda \bE))\ \forall \bA \in \Upsilon(\Omega, \bE)\\
\forall \bPhi \in S(\Omega)\ \forall m \in \bN: p ( R ( \bA_\e, \bPhi_\e ) ) = O(\e^m)).
\end{gather*}
or, in the \sm,
\begin{gather*}
   \forall k,l \in \bN_0\ \forall p \in \csn ( \Gamma(\Omega, \lambda \bE))\ \forall \bA \in \Upsilon(\Omega,\bE), \bA_1, \dotsc, \bA_k \in \Upsilon_0(\Omega, \bE) \\\
\forall \bPhi \in S(\Omega), \bPhi_1, \dotsc, \bPhi_l \in S_0(\Omega)\ \forall m \in \bN:\\
p ( (\ud_1^k \ud_2^l R) ( \bA_\e, \bPhi_\e )(\bA_{1\e}, \dotsc, \bA_{k\e}, \bPhi_{1\e}, \dotsc, \bPhi_{l\e}) ) = O(\e^m).
 \end{gather*}

The subspace of negligible elements of $\cE(\Omega, \lambda \bE)$ is denoted by $\cN(\Omega, \lambda \bE)$ and we set $\cN(\Omega) \coleq \cN(\Omega \times \bR)$.
\end{definition}

As before, we have (with an identical proof):

\begin{theorem}Let $u \in \cD'(\Omega, \lambda \bE)$ and $t \in C^\infty(\Omega, \lambda \bE)$. Then the following holds:
 \begin{enumerate}[label=(\roman*)]
  \item $\iota u$ is moderate.
  \item $\sigma t$ is moderate.
  \item $(\iota - \sigma)(t)$ is negligible.
  \item If $\iota(u)$ is negligible then $u=0$.
 \end{enumerate}
\end{theorem}

\begin{proposition}\label{proposition.3.39}
 Suppose $R \in \cE(\Omega, \lambda \bE)$ is moderate. Then it is negligible if and only if $\forall K \csub M$ $\forall \bA \in \Upsilon(\Omega, \bE)$ $\forall \bPhi \in S(\Omega)$:
 \[ \sup_{p \in K} \norm{ R ( \bA_\e, \bPhi_\e )(p) } = O(\e^m) \]
 where $\norm{\cdot}$ is any norm on $\lambda \bE$.
 \end{proposition}

Again,  $\cE_M(\Omega, \lambda \bE)$ is a submodule of $\cE(\Omega, \lambda \bE)$ and $\cN(\Omega, \lambda \bE)$ is a submodule of $\cE_M(\Omega, \lambda \bE)$, so we can define the quotient module:

\begin{definition}$\cG(\Omega, \lambda \bE) \coleq \cE_M(\Omega, \lambda \bE) / \cN ( \Omega, \lambda \bE)$ is called the space of nonlinear generalized functions on $\Omega$ with values in $\lambda \bE$.
\end{definition}

Moreover, the local and global settings are compatible in the following sense:

\begin{theorem}\label{theorem.3.41}
 Let $(U,\varphi)$ be a chart on $M$ where $E$ is trivial, $R \in \cE(U, \lambda E)$ and $R' \in \cE(\varphi(U), \bE)$ its local expression. Then $R$ is moderate (negligible) if and only if $R'$ is moderate (negligible).
\end{theorem}
\begin{proof}
This is clear from the definitions.
\end{proof}

It follows that we have an isomorphism $\cG(U, \lambda E) \cong \cG(\varphi(U), \lambda \bE)$.

\begin{theorem}The pullback $\mu^* \colon \cE(N, F) \to \cE(M, E)$ along a vector bundle isomorphism $\mu \colon E \to F$ preserves moderateness and negligibility and hence is also well-defined as a mapping $\cG(N, F) \to \cG(M, E)$.
\end{theorem}

\begin{proof}This follows directly from \Autoref{proposition.3.13} and \Autoref{proposition.3.25}.
\end{proof}

\begin{theorem}Pullbacks commute with the embeddings, i.e., $\mu^* \circ \iota = \iota \circ \mu^*$ and $\mu^* \circ \sigma = \sigma \circ \mu^*$.
\end{theorem}

\begin{proof}
 As in \Autoref{theorem.3.5}, this is an elementary consequence of the definitions.
\end{proof}

\subsection{The sheaf property}\label{subsection.3.5}

The key to obtaining the sheaf property of the quotient will be to employ so-called \emph{locality conditions} (first introduced in \cite{papernew} and studied in more detail in \cite{specfull}). Again, in this whole subsection we fix a manifold $M$, a vector bundle $E \to M$ and a homogeneous functor $\lambda$.

\begin{definition}Let a mapping
 \[ \ell \colon \TO(M, E) \times \SO(M) \times M \to X_\ell \]
be given, where $X_\ell$ is some set (depending on $\ell$). An element $R \in \cE(M, \lambda E)$ is called \emph{$\ell$-local} if for all tuples $(A, \Phi, x)$ and $(B, \Psi, y)$ in the domain of $\ell$ the following implication holds:
\begin{gather*}
\ell ( A, \Phi, x ) = \ell ( B, \Psi, y ) \Longrightarrow R ( A, \Phi) (x) = R ( B, \Psi) (y).
\end{gather*}
We denote the set of all $\ell$-local elements of $\cE(M,\lambda E)$ by \emph{$\cE[\ell](M,\lambda E)$}.
\end{definition}

To explicitly display the mapping $\ell$ when we speak of an $\ell$-local generalized section, instead of $\ell$ we may write its value on $(A, \Phi, x)$ or $(A, \vec\varphi, x)$ with $\vec\varphi$ the kernel of $\Phi$ given by \eqref{equation.2.5} instead, as is exemplified in the following proposition:

\begin{proposition}\label{proposition.3.45}
 Let $u \in \cD'(M, \lambda E)$ and $t \in \Gamma(M, \lambda E)$. Then
\begin{enumerate}[label=(\roman*)]
 \item\label{proposition.3.50.1} $\iota u$ is $(\lambda A)(x,.)\otimes\vec\varphi(x)$-local,
 \item\label{proposition.3.50.2} $\sigma t$ is $x$-local.
\end{enumerate}
\end{proposition}

\begin{proof}
 This is clear from the definitions.
\end{proof}

Note also that it is more convenient to write $\vec\varphi(x)$ than $\ev_x \circ \Phi$ (where $\ev_x$ is evaluation at $x$).

The purpose of locality conditions is to obtain a subspace of $\cE(M, \lambda E)$ containing $\cD'(M, \lambda E)$ and $\Gamma(M, \lambda E)$ where the usual constructions involved in obtaining the sheaf property work out. As we will see later, a natural choice for $\ell$ in our context is
\begin{equation}\label{equation.3.26}
\fell ( A, \vec\varphi, x) = \varinjlim_{U \in \cU_x} \bigl( A|_{\bigcup_{x' \in U} ( \{x'\} \times \carr \vec\varphi(x') ) }, \vec\varphi|_U, x \bigr ).
\end{equation}

Reformulating \eqref{equation.3.26}, an element $R\in \cE(M, \lambda E)$ is $\fell$-local if the following holds: given $A,B,\Phi,\Psi$, whenever there is an open set $U$ with $\vec\varphi|_U = \vec\psi|_U$ and such that for all $x \in U$ and $y \in \carr \vec\varphi(x)$, $A(x,y) = B(x,y)$ holds, we have $R(A, \Phi)|_U = R(B, \Psi)|_U$. Both locality types mentioned in \Autoref{proposition.3.45} imply this locality type, hence the resulting basic space will be large enough to accommodate distributional and smooth tensor fields.

\begin{proposition}If $R \in \cE(M, \lambda E)$ is $(\lambda A)(x,.)\otimes\vec\varphi(x)$-local or $x$-local then it is $\fell$-local.
\end{proposition}
\begin{proof}
Suppose that $\fell(A, \vec\varphi, x) = \fell (B, \vec \psi, y)$. Then in the first case, $x=y$, $\vec\varphi(x) = \vec\psi(x)$, and there is $U \in \cU_x$ such that
\begin{align*}
 & A(x',y') = B(x', y') & & \forall x' \in U, y' \in \carr \vec\varphi(x') \\
 \Rightarrow & (\lambda A)(x', y') = (\lambda B)(x', y') & & \forall x' \in U, y' \in \carr \vec\varphi(x') \\
 \Rightarrow & (\lambda A)(x, y') = (\lambda B)(x, y') & & \forall y' \in \carr \vec\varphi(x) \\
 \Rightarrow & (\lambda A)(x, .) \otimes \vec\varphi(x) = (\lambda B)(y, .) \otimes \vec\psi(y).
\end{align*}
The second claim is clear from $x=y$.
\end{proof}

So, although smoothing operators $\Theta = A \otimes \Phi$ have been decoupled and appear as pairs $(A, \Phi)$ in the definition of the basic space, the locality condition \eqref{equation.3.26} recouples their supports and hence makes it possible to obtain sheaf properties. 

With this we can define the following restriction mapping.

\begin{proposition}Let $U \subseteq M$ be open and $R \in \cE[\fell](U, \lambda E)$. Then for each open subset $V \subseteq U$ there exists a unique mapping $R|_V \in \cE[\fell](V,\lambda E)$, called the \emph{restriction} of $R$ to $V$, such that for all $A \in \TO(V,E)$, $\vec \varphi \in \SK(V)$, $B \in \TO(U,E)$, $\vec\psi \in \SK(U)$ and each open subset $W \subseteq V$ the conditions
\begin{equation}\label{equation.3.27}
 \begin{aligned}
  \vec\varphi|_W &= \vec\psi|_W \\
  A(x,y) &= B(x,y)\quad \forall x \in W\ \forall y \in \carr \vec\varphi(x)
 \end{aligned}
\end{equation}
imply that
\[ R|_V(A, \vec\varphi)|_W = R(B, \vec\psi)|_W. \]
If in addition we have $A_1, \dotsc, A_k \in \TO(V,E)$, $\vec\varphi_1, \dotsc, \vec\varphi_l \in \SK(V)$ as well as $B_1, \dotsc, B_k \in \TO(U,E)$, $\vec\psi_1, \dotsc, \vec\psi_l \in \SK(U)$ such that $\vec\varphi_i|_W = \vec\psi_i|_W$ and $A_j(x,y) = B_j(x,y)$ $\forall x \in W$ $\forall y \in \carr \vec\varphi(x) \cup \carr \vec\varphi_1 \cup \dotsc \cup \carr \vec\varphi_l(x)$ then
\begin{gather*}
(\ud_1^k \ud_2^l R|_V)(A,\vec\varphi)(A_1, \dotsc, A_k, \vec\varphi_1, \dotsc, \vec\varphi_l)|_W \\
= (\ud_1^k \ud_2^l R)(B, \vec\psi)(B_1, \dotsc, B_k, \vec\psi_1, \dotsc, \vec\psi_l)|_W.
\end{gather*}
The map $R \mapsto R|_V$ is $C^\infty$-linear, i.e., $(fR)|_V = f|_V \cdot R|_V$ for $f \in C^\infty(U)$. Moreover, if $V' \subseteq V$ is open then $(R|_V)|_{V'} = R|_{V'}$.
\end{proposition}

\begin{proof}
The proof is similar to that of \cite[Theorem 9]{specfull} or \cite[Theorem 8]{papernew}, but we give it explicitly because some details differ.

Fix $A \in \TO(V,E)$ and $\vec\varphi \in \SK(V)$. We are going to define $R|_V(A,\vec\varphi)$ using the sheaf property of $\Gamma(M, \lambda E)$.

Let $X$ be an open subset of $V$ such that $\overline{X} \subseteq V$ is compact. Choose $\rho_X \in \cD(V)$ with $\rho_X = 1$ on $X$; then $\rho_X \vec\varphi \in C^\infty(U, \Gamma_c(U, \Vol(M)))$. As $\supp \rho_X$ is compact in $V$ we can find a compact set $L_X$ in $V$ such that $\supp \vec\varphi(x) \subseteq L_X$ for all $x \in \supp \rho_X$. Now choose a function $\chi_X \in \cD(V)$ with $\chi_X = 1$ on $L_X$. We then have $(\rho_X \otimes \chi_X)A \in \TO(U,E)$ and define $f_X \coleq R((\rho_X \otimes \chi_X)A, \rho_X \vec\varphi)|_X \in \Gamma(X, \lambda E)$.

We now use that for any open sets $X_1, X_2 \subseteq V$ with compact closures in $V$ and $X_1 \cap X_2 \ne \emptyset$ we have $f_{X_1}|_{X_1 \cap X_2} = f_{X_2}|_{X_1 \cap X_2}$; in fact,
\begin{align*}
 (\rho_{X_1} \vec \varphi)(x) &= \vec\varphi(x) = (\rho_{X_2} \vec\varphi)(x), \\
 ((\rho_{X_1} \otimes \chi_{X_1})A)(x,y) &= ((\rho_{X_2} \otimes \chi_{X_2})A)(x,y)
\end{align*}
for all $x \in X_1 \cap X_2$ and $y \in \carr \vec\varphi(x)$ implies
\begin{align*}
 f_{X_1}|_{X_1 \cap X_2} &= R((\rho_{X_1} \otimes \chi_{X_1})A, \rho_{X_1} \vec\varphi)|_{X_1 \cap X_2} \\
 &= R((\rho_{X_2} \otimes \chi_{X_2})A, \rho_{X_2} \vec\varphi))|_{X_1 \cap X_2} = f_{X_2}|_{X_1 \cap X_2}
\end{align*}
by $\fell$-locality of $R$. Consequently, the family $(f_X)_X$ defines a unique element $f \in \Gamma(V, \lambda E)$ satisfying $f|_X = f_X$ for all $X$. We define $R|_V (A, \vec\varphi) \coleq f$.

For the \sm, we note that $R$ is smooth because this can be tested locally.

To show that $R|_V$ has the claimed property suppose that \eqref{equation.3.27} holds for some $A,B,\vec\varphi,\vec\psi,W$ as stated. Let $X \subseteq W$ be such that $\overline{X}$ is compact in $W$ and hence also in $V$; then
\[ R|_V(A,\vec\varphi)|_X = R((\rho_X \otimes \chi_X)A, \rho_X\vec\varphi)|_X = R(B, \vec\psi)|_X \]
because
\begin{align*}
 (\rho_X\vec\varphi)(x) &= \vec\varphi(x) = \vec\psi(x) && \forall x \in X,\\
 ((\rho_X \otimes \chi_X)A)(x,y) &= A(x,y) = B(x,y) && \forall x \in X\ \forall y \in \carr \vec\varphi(x).
\end{align*}
As we can cover $W$ by such sets $X$, this shows the claim.

To obtain the second claim, we need to differentiate
\[ R|_V(A + t_1 A_1 + \dotsc + t_k A_k, \vec\varphi + s_1 \vec\varphi_1 + \dotsc + s_l \vec\varphi_l)|_W \]
with respect to each variable $t_1, \dotsc, t_k, s_1, \dotsc, s_l$ at $0$. This expression equals
\[ R(B + t_1 B_1 + \dotsc + t_k B_k, \vec\psi + s_1 \vec\psi_1 + \dotsc + s_l \vec\psi_l)|_W \]
if
\[ \vec\varphi + s_1 \vec\varphi_1 + \dotsc + s_l \vec\varphi_l = \vec\psi + s_1 \vec\psi_1 + \dotsc + s_l \vec\psi_l \]
on $W$ and
\[ (A + t_1 A_1 + \dotsc + t_k A_k)(x,y) = (B + t_1 B_1 + \dotsc + t_k B_k)(x,y) \]
for all $x \in W$ and $y \in \carr ( \vec\varphi + s_1 \vec\varphi_1 + \dotsc + s_l \vec\varphi_l)(x)$, which is implied by the assumptions and gives the second claim.

To see that $R|_V$ is $\fell$-local let $\widetilde U \subseteq V$ be open such that
\begin{align*}
 \vec\varphi(x) &= \vec\psi(x) &&  \forall x \in \widetilde U, \\
 A(x,y) &= B(x,y) && \forall x \in \widetilde U\ \forall y \in \carr \vec\varphi(x)
\end{align*}
for some $\vec\varphi,\vec\psi \in \SK(V)$ and $A,B \in \TO(U,E)$. Let an open set $X \subseteq \widetilde U$ have compact closure in $\widetilde U$ and take $\rho_X$, $L_X$ and $\chi_X$ from above.
Then
\begin{gather*}
 \vec\varphi(x) = (\rho_X \vec\varphi)(x) = (\rho_X \vec\psi)(x) = \vec\psi(x) \\
 A(x,y) = ((\rho_X \otimes \chi_X)A)(x,y) = ((\rho_X \otimes \chi_X)B)(x,y) = B(x,y)
\end{gather*}
for $x \in X$ and $y \in \carr \vec\varphi(x)$. Hence,
\begin{align*}
 R|_V(A,\vec\varphi)|_X & = R((\rho_X \otimes \chi_X)A, \rho_X\vec\varphi)|_X \\
 &= R((\rho_X \otimes \chi_X)B, \rho_X \vec\psi)|_X = R|_V(B, \vec\psi)|_X,
\end{align*}
which is $\fell$-locality of $R|_V$.

For uniqueness, suppose we are given $S \in \cE[\fell](V, \lambda E)$ satisfying the stated property. Given $A \in \TO(V,E)$ and $\vec\varphi \in \SK(V)$, choose an open set $X \subseteq V$ with $\overline{X} \subseteq V$ compact and $\rho_X, \chi_X \in \cD(V)$ as above. Then by $\fell$-locality we have
\[ R|_V (A, \vec\varphi)|_X = R((\rho_X \otimes \chi_X)A, \rho_X\vec\varphi)|_X = S(A, \vec\varphi)|_X \]
which shows uniqueness of $R|_V$.

For $C^\infty$-linearity of $R \mapsto R|_V$ suppose that $R_1,R_2 \in \cE[\fell](U, \lambda E)$ are given. Then with $X, \rho_X, \chi_X$ as above,
\begin{gather*}
 (R_1+R_2)|_V(A, \vec\varphi)|_X = (R_1 + R_2)((\rho_X \otimes \chi_X)A, \rho_X\vec\varphi)|_X \\
= R_1 ( ( \rho_X \otimes \chi_X)A, \rho_X\vec\varphi)|_X + R_2 ( (\rho_X \otimes \chi_X)A, \rho_X \vec\varphi)|_X \\
= R_1|_V ( A, \vec\varphi)|_W + R_2|_V ( A, \vec\varphi)|_X \\
= (R_1|_V + R_2|_V ) ( A, \vec\varphi)|_X.
\end{gather*}
and for $f \in C^\infty(U)$,
\begin{gather*}
 (f R)|_V ( A, \vec\varphi)|_X = (fR)((\rho_X \otimes \chi_X)A, \rho_X \vec\varphi)|_X \\
= f|_X R ( (\rho_X \otimes \chi_X)A, \rho_X \vec\varphi)|_X = (f|_V R|_V)(A, \vec\varphi)|_X.
\end{gather*}
Finally, for transitivity let $W \subseteq V'$. Then for any open $X \subseteq W$ with $\overline{X} \subseteq W$ compact and suitable $\rho_X \in \cD(V')$, $\chi_X \in \cD(V')$ as above,
\begin{gather*}
 ((R|_V)|_{V'})(A, \vec\varphi)|_X = R|_V ( ( \rho_X \otimes \chi_X ) A, \rho_X \vec\varphi)|_X \\
= R ( ( \rho_X \otimes \chi_X)A, \rho_X \vec\varphi)|_X = R|_{V'} ( A, \vec\varphi)|_X
\end{gather*}
gives the claim.
\end{proof}

We can now go about proving the sheaf property for $\cG[\fell](M, \lambda E)$. The first step is to show that moderateness and negligibility localize:

\begin{theorem}\label{theorem.3.48}
Let an open subset $U \subseteq M$, an open cover $(U_i)_i$ of $U$ and $R \in \cE[\fell](U, \lambda E)$ be given. Then $R$ is moderate or negligible if and only if all $R|_{U_i}$ are.
\end{theorem}
\begin{proof}
In order to test $R|_{U_i}$ according to \Autoref{definition.3.22} (\si) we have to estimate derivatives of the components of
\[ R|_{U_i} (\bA_\e,  \bvarphi_\e )_\e(x) \]
for $\bA \in \Upsilon(U_i, E)$, $\bvarphi \in S(U_i)$ and $x$ in a relatively compact open subset $W \subseteq U_i$. By Lemma \ref{lemma.3.31} we can choose $\bA' \in \Upsilon(U,E)$, $\bvarphi' \in S(U)$ and $\e_0>0$ such that $\forall \e<\e_0$:
\begin{equation}\label{equation.3.28}
 \begin{aligned}
\bvarphi_\e & = \bvarphi_\e' \textrm{ on }W \\
\bA_\e(x,.) & = \bA_\e'(x,.)\textrm{ on }\carr \bvarphi_\e(x)\ \forall x \in W.
 \end{aligned}
\end{equation}

Hence, the expression to be estimated equals $R(\bA'_\e, \bvarphi'_\e)(x)$ on this neighborhood for small $\e$, so moderateness and negligibility of $R|_{U_i}$ are implied by the same property of $R$.

For the \sm\ we also have to estimate
\begin{equation}\label{equation.3.29}
 (\ud_1^k \ud_2^l R|_{U_i})(\bA_\e, \bvarphi_\e)(\bA_{1\e}, \dotsc, \bA_{k\e}, \bPhi_{1\e}, \dotsc, \bPhi_{l\e})(x).
\end{equation}
As above we can choose $\bA' \in \Upsilon(U,E)$, $\bA_1', \dotsc, \bA_k' \in \Upsilon_0(U,E)$ and $\bvarphi' \in S(U)$, $\bvarphi_1', \dotsc, \bvarphi_l' \in S_0(U)$ such that for small $\e$, the statement analogous to \eqref{equation.3.28} holds such that \eqref{equation.3.29} equals
\[ (\ud_1^k \ud_2^l R)(\bA_\e', \bvarphi_\e')(\bA_{1\e}', \dotsc, \bA_{k\e}', \bPhi'_{1\e}, \dotsc, \bPhi'_{l\e})(x) \]
which gives the claim.

Conversely, in order to estimate $R(\bA_\e, \bvarphi_\e)(x)$ for $\bA \in \Upsilon(U,E)$, $\bvarphi \in S(U)$ and $x$ in a compact subset $K \subseteq U$, we can assume without limitation of generality that $K \subseteq U_i$ for some $i$. Then, the expression to be estimated equals $R|_{U_i} ( \bA_\e|_{U_i \times U_i}, \rho_{U_i} \bvarphi_\e)(x)$ for $x$ in an open neighborhood of $K$ and small $\e$, where $\rho_{U_i} \in C^\infty(U_i, \cD(U_i))$ equals $1$ on an arbitrary open neighborhood of the diagonal in $U_i \times U_i$, whence moderateness and negligibility of $R$ follows if all $R|_{U_i}$ have that property.
\end{proof}

\begin{theorem}
 $\cG[\fell](\unterstrich, \lambda E)$ is a sheaf of $C^\infty$-modules on $M$.
\end{theorem}

\begin{proof}
Let $U \subseteq M$ be open and $(U_i)_i$ a covering of $U$ by open sets. Suppose that we are given functions $R_i \in \cE_M[\fell](U_i, \lambda E)$ such that $R_i|_{U_i \cap U_j} - R_j|_{U_i \cap U_j} \in \cN[\fell] ( U_i \cap U_j, \lambda E)$ for all $i,j$ with $U_i \cap U_j \ne \emptyset$.

Choose a partition of unity $(\chi_i)_i$ on $U$ subordinate to $(U_i)_i$ and for each $i$ a function $\rho_i \in C^\infty(U_i, \cD(U_i))$ which is equal to 1 on the diagonal in $U_i \times U_i$. We define the mapping $R \colon \TO(U, E) \times \SK(U) \to \Gamma(U, \lambda E)$ by
\begin{equation}\label{equation.3.30}
R(A, \vec\varphi)(x) \coleq \sum_i \chi_i (x) \cdot R_i ( A|_{U_i \times U_i}, \rho_i \vec\varphi)(x).
\end{equation}

Clearly $R$ is smooth in case of the \sm, as each $R_i$ is smooth and the sum is locally finite.

It is easily verified that $R$ is $\fell$-local: suppose $\vec\varphi = \vec\psi$ on an open set $W \subseteq U$ and $A(x.,) = B(x,.)$ on $\carr \vec\varphi(x)$ $\forall x \in W$ and fix $x$; if $\chi_i(x) \ne 0$ then $x \in U_i$ and hence
$\rho_i \vec\varphi = \rho_i\vec\psi$ on $W \cap U_i$ and $A|_{U_i \times U_i}(x',.) = B|_{U_i \times U_i} (x',.)$ on $\carr \rho_i \vec\varphi_\e(x)$ for all $x' \in W \cap U_i$.

For moderateness of $R$ we estimate derivatives of (components of) $R(\bA_\e, \bvarphi_\e)(x)$ for $\bA \in \Upsilon(U,E)$, $\bvarphi \in S(U)$ and $x$ in a compact subset $K \subseteq U$. There exists a finite index set such that for $x \in K$, the sum in \eqref{equation.3.30} only has to be taken for $i$ in this index set. By the Leibniz rule it then suffices to estimate derivatives of $R_i (\bA_\e|_{U_i \times U_i}, \rho_i \bvarphi_\e ) (x)$ for $x$ in $\supp \chi_i \cap K$. But this expression has moderate growth by assumption because $\bA|_{U_i \times U_i}$ is admissible and $\rho_i \bvarphi$ is a test object on $U_i$; the case of the \sm\ is seen similarly.

Next, we show that $R|_{U_j} - R_j$ is negligible for all $j$. Fix $\bA \in \Upsilon(U_j, E)$ and $\bvarphi \in S(U_j)$ for testing and a compact set $K$ in $U_j$. Using Lemma \ref{lemma.3.31} choose $\bB \in \TO(U, E)^I$ and $\bpsi \in S(U)$ such that for small $\e$, $\bvarphi_\e = \bpsi_\e$ on a neighborhood $W$ of $K$ and $\bA(x,.) = \bB(x,.)$ on $\carr \bvarphi_\e(x)$ for $x \in W$.

As above, there is a finite index set such that for $x \in K$ the sum in 
\[ (R|_{U_j} - R_j)(\bA_\e, \bvarphi_\e)(x) = \sum_i \chi_i(x) \cdot ( R_i ( \bA_\e|_{U_i \times U_i}, \rho_i\bvarphi_\e) - R_j(\bA_\e, \bvarphi_\e))(x) \]
runs only over this index set. Hence, by \Autoref{proposition.3.34} it suffices to estimate $R_i ( \bA_\e|_{U_i \times U_i}, \rho_i \bvarphi_\e)(x) - R_j(\bA_\e, \bvarphi_\e)(x)$ for $x$ in the relatively compact subset $\supp \chi_i \cap K \subseteq U_i \cap U_j$. Let $\rho \in C^\infty(U_i \cap U_j, \cD(U_i \cap U_j))$ be equal to 1 on a neighborhood of the diagonal in $(U_i \cap U_j) \times (U_i \cap U_j)$. Then for $p \in \supp \chi_i \cap K$ and small $\e$,
\begin{gather*}
 R_i ( \bA_\e|_{U_i \times U_i}, \rho_i \bvarphi_\e)(p) - R_j ( \bA_\e, \bvarphi_\e)(p) \\
= R_i|_{U_i \cap U_j} ( \bB_\e|_{(U_i \cap U_j) \times (U_i \cap U_j)}, \rho \rho_i \bpsi_\e)(p) - R_j|_{U_i \cap U_j} ( \bA_\e|_{U_i \cap U_j}, \rho \bvarphi_\e)(p) \\
= R_i|_{U_i \cap U_j} ( \bB_\e|_{(U_i \cap U_j) \times (U_i \cap U_j)}, \rho\bpsi_\e)(p) \\
\qquad - R_j|_{U_i \cap U_j} ( \bB_\e|_{(U_i \cap U_j)\times (U_i \cap U_j)}, \rho\bpsi_\e)(p).
\end{gather*}
For the first equality we use that for the first terms, $\rho_i \bvarphi_\e = \rho \rho_i \bpsi_\e$ on an open neighborhood $W'$ of $p$ in $W$; moreover, for $x \in W'$ and $y \in \carr ( \rho_i \bvarphi_\e)(x) \subseteq \carr \bvarphi_\e(x)$, $\bA_\e(x,y) = \bB_\e(x,y)$. The equality of the second terms is clear.

For the second equality we use that $\rho \rho_i \bpsi_\e = \rho \bpsi_\e$ on a neighborhood of $p$ for small $\e$, which is clear; $\rho \bvarphi_\e = \rho \bpsi_\e$ on a neighborhood $W'$ of $p$ for small $\e$; and for $x \in W'$ and $y \in \carr \rho \bvarphi_\e$, $\bA_\e(x,y) = \bB_\e(x,y)$.

Finally, $R$ is $\fell$-local: suppose $p \in M$, $\e \in I$, $W \in \cU_p(M)$, and
\begin{gather}
 \label{equation.3.31} \bvarphi_\e = \bpsi_\e \textrm{ on }W, \\
\label{equation.3.32} \forall x \in W: \bA_\e(x,.) = \bB_\e(x,.) \textrm{ on }\carr \bvarphi_\e(x).
\end{gather}
Then we want that for each $i$ such that $p \in \carr \chi_i$,
\[ R_i ( \bA|_{U_i \times U_i}, \rho_i \bvarphi)_\e(p) = R_i ( \bB|_{U_i \times U_i}, \rho_i \bpsi)_\e(p). \]
But \eqref{equation.3.31} implies $\rho_i \bvarphi = \rho_i \bpsi$ on $W \cap U_i$, and \Autoref{equation.3.32} implies that $\forall x \in W \cap U_i: \bA_\e(x,.) = \bB_\e(x,.)$ on $\carr (\rho_i \bvarphi_\e)(x)$, which gives the claim.
\end{proof}

\begin{theorem}\label{theorem.3.50}
 The embeddings $\iota$ and $\sigma$ are sheaf morphisms into $\cG[\fell](\unterstrich, \lambda E)$.
\end{theorem}
\begin{proof}
 We refer to \cite[Theorem 8.11, p.~211]{bigone} for the proof, which transfers without problems to our setting.
\end{proof}

One can also show that the isomorphism of \Autoref{lemma.3.36} preserves locality:
\begin{proposition}\label{proposition.3.51}
There is a sheaf isomorphism
\begin{align*}
 \cG[\fell](\unterstrich, \lambda E) & \cong \Gamma(\unterstrich, \lambda E)\otimes_{C^\infty(\unterstrich)} \cG[\fell](\unterstrich) \\
& \cong \Hom_{C^\infty(\unterstrich)} ( \Gamma(\unterstrich, (\lambda E)^*), \cG[\fell](\unterstrich)).
\end{align*}
\end{proposition}
We refer to \cite[Theorem 8.15, p.~212]{bigone} for the proof.

\begin{proposition}
An element $R \in \cE[\ell](M, \lambda E)$ is moderate (negligible) if and only if all its coordinates with respect to a covering of $M$ by trivializations of $E$ are moderate (negligible).
\end{proposition}
\begin{proof}
This follows from \Autoref{theorem.3.48} and \Autoref{lemma.3.36}.
\end{proof}

All results of this subsection hold, mutatis mutandis, also in the local setting.

\subsection{Association}

To reconcile the theory of nonlinear generalized functions with the theory of distributions, the concept of association is useful. We cover only its mere basics here and refer to \cite{MOBook,GKOS} for a more in-depth discussion.

\begin{definition}\label{definition.3.53}
An element $R \in \cE(M, \lambda E)$ is said to be \emph{associated to zero} if for all $\bA \in \Upsilon(M, E)$, $\Phi \in \SO(M)$, $v \in \Gamma(M, (\lambda E)^*)$ and $\omega \in \Gamma_c(M, \Vol(M))$ we have
\[ \int (R(\bA_\e, \bPhi_\e) \cdot v) \omega \to 0 \qquad (\e \to 0). \]
In this case we write $R \approx 0$. For $R,S \in \cE(M,\lambda E)$ we say that $R$ is associated to $S$, written $R \approx S$, if $R-S \approx 0$.

An element $R \in \cE(M, \lambda E)$ is said to have $u \in \cD'(M, \lambda E)$ as an associated distribution if $R \approx \iota(u)$.
\end{definition}

Clearly, association defines a congruence relation.

Because every $R \in \cN(M, \lambda E)$ is associated to zero, association is well-defined on the quotient and two elements $\widehat R, \widehat S \in \cG(M, \lambda E)$ are called associated if any two (hence all) of their representatives are associated.

The definition in the local theory is similar:

\begin{definition}
An element $R \in \cG(\Omega, \lambda \bE)$ is said to be \emph{associated to zero} if for all $\bA \in \Upsilon(\Omega, \bE)$, $\Phi \in \SO(\Omega)$ and $\omega \in \cD(\Omega)$ we have
\[ \int R(\bA_\e, \bPhi_\e)(x) \omega(x)\,\ud x \to 0 \quad \textrm{in }\bE \qquad (\e \to 0) \]
In this case we write $R \approx 0$, and $R \approx S$ means that $R - S \approx 0$. Moreover, an element $R \in \cG(\Omega, \lambda \bE)$ is said to have $u \in \cD'(\Omega,\bE)$ as an associated distribution if $R \approx \iota(u)$.
\end{definition}

\begin{proposition}
 Let $(U,\varphi)$ be a chart on $M$ where $E$ is trivial, $R \in \cE(U, \lambda E)$ and $\widetilde R \in \cE(\varphi(U), \lambda \bE)$ its local expression. Then $R \approx 0$ if and only if $\widetilde R \approx 0$.
\end{proposition}

Further down we will modify this concept and talk of \emph{metric association}.

\begin{proposition}
A generalized section $R \in \cE[\fell](M, \lambda E)$ is associated to zero if and only if all its coordinates with respect to a covering of $M$ by trivializations of $E$ are zero.
\end{proposition}
\begin{proof}
This follows immediately from the definitions.
\end{proof}
 
\subsection{Invertibility}

In this section we will discuss multiplicative invertibility of generalized functions which are bounded away from zero in a suitable way. The central result (which, as a comparison with the notion of a \emph{strictly nonzero} generalized function in \cite{GKOS} suggests, is not the most general one, but sufficient for our purposes) is the following.

\begin{theorem}\label{theorem.3.57}Suppose $R \in \cE_M(M)$ satisfies
\begin{equation}\label{equation.3.33}
\begin{aligned}
&\forall K \csub M\ \exists C_K>0\ \forall \bA \in \Upsilon(M,E)\ \forall \bPhi \in S(M)\\
&\exists \e_0=\e_0(K)>0\ \forall \e<\e_0: \inf_{x \in K} R(\bA_\e, \bPhi_\e)(x) \ge C_K.
\end{aligned}
\end{equation}
Then there exists $\widetilde R \in \cE_M(M)$ such that
\begin{enumerate}[label=(\roman*)]
 \item\label{gnugger.1} $\forall A \in \TO(M, E)$ $\forall \Phi \in \SO(M)$ $\forall x \in M$: $\widetilde R(A, \Phi)(x)>0$,
 \item\label{gnugger.2} $\forall K \csub M$ $\forall \bA \in \Upsilon(M,E)$ $\forall \bPhi \in S(M)$ $\exists \e_0>0$\\
 $\forall \e<\e_0$ $\forall x \in K$: $\widetilde R(\bA_\e, \bPhi_\e)(x) = R(\bA_\e, \bPhi_\e)(x)$.
 \item\label{gnugger.3} $\widetilde R - R$ is negligible.
 \item\label{gnugger.4} $S(A,\Phi) \coleq 1 / ( \widetilde R(A, \Phi) )$ defines an element $S \in \cE_M(M)$.
 \item\label{gnugger.5} $R S - 1$ is negligible.
 \item\label{gnugger.6} If $R$ is $\fell$-local then $S$ is $\fell$-local.
\end{enumerate}
\end{theorem}
\begin{proof}
 Without limitation of generality we can assume that $C_K \ge C_L$ in \eqref{equation.3.33} if $K,L \subseteq M$ are compact sets with $K \subseteq L$. To see this, take a compact exhaustion $(K_n)_{n \in \bN}$ of $M$ and fix $C_n \coleq C_{K_n}$. Let $C_n' \coleq \min_{1 \le i \le n} C_i$ for all $n \in \bN$. For $K = K_n$, property \eqref{equation.3.33} holds with $C_K =C_n'$; in fact, fixing $\bA$, $\bPhi$ and taking $\e_0 < \min_{1 \le i \le n} \e_0(K_i)$ we have $\inf_{x \in K_n} R(\bA_\e, \bPhi_\e)(x) \ge C_n \ge C_n'$. Now we have $K \subseteq K_n$ and $L \subseteq K_l$ for some $n \le l$, so $C_k' \ge C_l'$ justifies our assumption.
 
Fix any Riemannian metric $g$ on $M$ and choose a continuous function $\rho \colon M \to (0, \infty)$ such that $B^g_{\rho(x)}(x)$ is compact in $M$ for all $x \in M$. Define $h(x) \coleq C_{B^g_{\rho(x)}(x)}$. We can find a continuous function $h_0 \colon M \to \bR$ with $0 < h_0 < h$.
Define the disjoint sets
 \begin{align*}
  A & \coleq \{ (x,t) \in M \times \bR\ |\ t \le h_0(x)/2 \}, \\
  B & \coleq \{ (x,t) \in M \times \bR\ |\ t \ge h_0(x) \}.
 \end{align*}
 As these are closed by continuity of $h_0$, we can choose a bump function $\chi \in C^\infty(M \times \bR, [0,1])$ such that $\chi|_A = 0$ and $\chi|_B = 1$. Define $\widetilde R \colon  \TO(M, E) \times \SO(M) \to  C^\infty(M)$ by
 \[ \widetilde R(A, \Phi)(x) \coleq R(A, \Phi)(x) \chi(x, R(A, \Phi)(x)) + 1 - \chi(x, R(A, \Phi)(x)).
 \]
In the \sm, it is clear that $\widetilde R$ is smooth if $R$ is so.

We now claim that
\begin{gather*}
\forall K \csub M\ \forall \bA \in \Upsilon(M,E)\ \forall \bPhi \in S(M)\ \exists \e_0>0\\
\ \forall \e<\e_0\ \forall x \in K: \chi(x, R(\bA_\e, \bPhi_\e)(x)) = 1.
\end{gather*}
It suffices to prove this for $x_0 \in M$ and $K \coleq B^g_\delta(x_0)$ where $\delta>0$ is chosen as follows: first, take $\delta_1>0$ such that $\rho(x_0)/2 \le \rho(x)$ for $d_g(x,x_0) < \delta_1$ (note for this that $\rho$ is continuous); then, set $\delta \coleq \min ( \rho(x_0)/4, \delta_1 )$. Choose $\e_0>0$ such that for $\e < \e_0$ and $x \in K$, $R(\bA_\e, \bPhi_\e)(x) \ge C_K$. For $x,y \in K$ we then have $d_g(x,y) \le 2 \delta \le \rho(x_0)/2 \le \rho(x)$, so $K \subseteq B^g_{\rho(x)}(x)$ and $C_K \ge C_{B^g_{\rho(x)}(x)} = h(x) \ge h_0(x)$, which implies $\chi(x, R(\bA_\e, \bPhi_\e)(x)) = 1$.

For \ref{gnugger.1}, if $\chi(x, R(A,\Phi)(x)) \ne 0$ then $R(A,\Phi)(x) > h_0(x)/2 > 0$ so $\widetilde R(A,\Phi)(x) \ne 0$.
\ref{gnugger.2} -- \ref{gnugger.3} are clear from the above.

For \ref{gnugger.4} smoothness of $S$ is clear and moderateness follows by applying the chain rule. Finally, \ref{gnugger.5} and \ref{gnugger.6} are evident.
\end{proof}

Note that if \eqref{equation.3.33} holds for one representative of an element of $\cG(M)$ then it also holds for every other representative. Moreover, the analogous statement of \Autoref{theorem.3.57} in the local case is valid:

\begin{theorem}\label{theorem.3.58}Suppose $R \in \cE_M(\Omega)$ satisfies
\begin{equation}\label{equation.3.34}
\begin{aligned}
&\forall K \csub \Omega\ \exists C_K>0\ \forall \bA \in \Upsilon(\Omega,\bE)\ \forall \bPhi \in S(\Omega)\\
&\exists \e_0=\e_0(K)>0\ \forall \e<\e_0: \inf_{x \in K} R(\bA_\e, \bPhi_\e)(x) \ge C_K.
\end{aligned}
\end{equation}
Then there exists $\widetilde R \in \cE_M(\Omega)$ such that
\begin{enumerate}[label=(\roman*)]
 \item\label{gnugger.1x} $\forall A \in \TO(\Omega, \bE)$ $\forall \Phi \in \SO(\Omega)$ $\forall x \in M$: $\widetilde R(A, \Phi)(x)>0$,
 \item\label{gnugger.2x} $\forall K \csub \Omega$ $\forall \bA \in \Upsilon(\Omega,\bE)$ $\forall \bPhi \in S(\Omega)$ $\exists \e_0>0$\\
 $\forall \e<\e_0$ $\forall x \in K$: $\widetilde R(\bA_\e, \bPhi_\e)(x) = R(\bA_\e, \bPhi_\e)(x)$.
 \item\label{gnugger.3x} $\widetilde R - R$ is negligible.
 \item\label{gnugger.4x} $S(A,\Phi) \coleq 1 / ( \widetilde R(A, \Phi) )$ defines an element $S \in \cE_M(M)$.
 \item\label{gnugger.5x} $R S - 1$ is negligible.
 \item\label{gnugger.6x} If $R$ is $\fell$-local then $S$ is $(\fell, x)$-local.
\end{enumerate}
\end{theorem}

If $(U,\varphi)$ is a chart then $R \in \cG(U)$ is invertible if and only if its local expression $\varphi_*R \in \cG(\varphi(U))$ is invertible; moreover, conditions \eqref{equation.3.33} and \eqref{equation.3.34} obviously correspond to each other:
\begin{lemma}\label{lemma.3.59}
Let $R \in \cE_M(U)$ and its local expression $R_0 \in \cE_M(\varphi(U))$ be given. Then the following are equivalent:
\begin{enumerate}[label=(\roman*)]
 \item\label{gla.1}  $\forall K \csub U$ $\exists C_K>0$ $\forall \bA \in \Upsilon(U, E)$ $\forall \bPhi \in S(U)$ $\exists \e_0>0$ $\forall \e<\e_0:$
 \[ \inf_{x \in K} R(\bA_\e, \bPhi_\e)(x) \ge C_K. \]
 \item\label{gla.2}  $\forall K \csub \varphi(U)$ $\exists C_K>0$ $\forall \bA \in \Upsilon(\varphi(U), \bE)$ $\forall \bPhi \in S(\varphi(U))$ $\exists \e_0>0$ $\forall \e<\e_0$:
 \[ \inf_{x \in K} R_0(\bA_\e, \bvarphi_\e)(x) \ge C_K. \]
\end{enumerate}
\end{lemma}

\subsection{Operations on generalized sections}

For the following definition, note that $E \mapsto \lambda_1 E \otimes \lambda_2 E$ is a homogeneous functor if $\lambda_1$ and $\lambda_2$ are so.

\begin{definition}\label{definition.3.60}
Let $\lambda_1, \lambda_2$ be homogeneous functors. The tensor product of $R \in \cE(M, \lambda_1 E)$ and $S \in \cE(M, \lambda_2 E)$ is defined as the element $R \otimes S \in \cE(M, \lambda_1 E \otimes \lambda_2 E)$ given by
\[ (R \otimes S)(A, \Phi) \coleq R(A, \Phi) \otimes S(A, \Phi). \]
\end{definition}

One quickly verifies the following properties:

\begin{lemma}
\begin{enumerate}[label=(\roman*)]
 \item The tensor product of moderate nonlinear generalized sections is moderate. If at least one factor is negligible, the tensor product is negligible. Consequently, this operation is well-defined on the quotient spaces.
 \item If $R$ and $S$ are $\fell$-local then $R \otimes S$ is $\fell$-local and $(R \otimes S)|_U = R|_U \otimes S|_U$ in $\cE[\fell](U, \lambda_1 E \otimes \lambda_2 E)$.
\end{enumerate}
\end{lemma}
Due to the vector bundle isomorphism $\lambda E \otimes (M \times \bR) \cong \lambda E$, \Autoref{definition.3.60} induces an $\cE(M)$-module structure on $\cE(M, \lambda E)$ and a $\cG(M)$-module structure on $\cG(M, \lambda E)$, and similarly for the spaces with $\fell$-locality.

We now introduce Lie derivatives of generalized sections. There are two natural choices for this, cf.~\cite{papernew,bigone} for a detailed discussion.

\begin{definition}
If $\lambda E$ be is a natural bundle, Lie derivatives with respect to $X \in \fX(M)$ are defined by
\begin{align*}
 (\widetilde\Lie_X R)(A, \Phi) \coleq \Lie_X ( R ( A, \Phi) ) \qquad (X \in \fX(M))
\end{align*}
and in the \sm\ additionally by
\[
 (\widehat \Lie_X R)(A,\Phi) \coleq \Lie_X ( R ( A, \Phi) ) - (\ud_1 R)(A, \Phi)(\Lie_X A) - (\ud_2 R)(A,\Phi)(\Lie_X \Phi).
\]
\end{definition}

We can extend $\widetilde \Lie_X$ even to generalized vector fields $X$:
\begin{definition}
 For $R \in \cE(M, \lambda E)$ and $X \in \cE(M, TM)$ we set
 \[ (\widetilde \Lie_X R)(A, \Phi) \coleq \Lie_{X(A, \Phi)} R(A, \Phi). \]
\end{definition}

The Lie derivatives $\widehat\Lie_X$ and $\widetilde\Lie_X$ have the following properties:
\begin{enumerate}[label={(\roman*)}]
\item \label{prop_derivprop.1} They satisfy the product rule
\begin{align*}
 \widehat\Lie_X (R \otimes S) &= \widehat\Lie_X R \otimes S + R \otimes \widehat\Lie_X S, \\
 \widetilde\Lie_X (R \otimes S) &= \widetilde\Lie_X R \otimes S + R \otimes \widetilde\Lie_X S.
\end{align*}
\item \label{prop_derivprop.2}$\widehat\Lie$ and $\widetilde\Lie$ are $\bR$-bilinear as maps $\fX(M) \times \cE(M, \lambda E) \to \cE(M, \lambda E)$.
 For $R \in \cE(M, \lambda E)$, $\widetilde\Lie_X F$ is $\cE(M)$-linear in $X$.
\item\label{5.3} $\widehat\Lie_X$ commutes with $\iota$, i.e., $\widehat\Lie_X \circ \iota = \iota \circ \Lie_X$.
\item\label{5.4} On $\Gamma(M, \lambda E)$, $\widehat\Lie_X$ and $\widetilde\Lie_X$ coincide with the classical Lie derivative $\Lie_X$ of smooth sections.
\item $\widetilde\Lie_X$ and $\widehat \Lie_X$ preserve moderateness and negligibility, hence are well-defined on the quotient.
\item $\widetilde\Lie_X$ commutes with $\iota$ on the level of association, i.e., $\widetilde \Lie_X(\iota u) \approx \iota(\Lie_X u)$.
\end{enumerate}

\section{Generalized tensor fields}\label{section.4}

We will now specialize the construction of Section \ref{section.3} to the tangent bundle $TM$ and the mixed tensor functors $\lambda = \tang^r_s$ ($r,s \ge 0$); moreover, we only consider $\fell$-local generalized tensor fields. We furthermore set $\TO(M) \coleq \TO(M, \tang M)$, $\Upsilon(U) \coleq \Upsilon(U, \tang M)$ and $\Upsilon_0(U) \coleq \Upsilon_0(U, \tang M)$.

\begin{definition}\label{definition.4.1}For $r,s \ge 0$ we set
 \begin{align*}
  \cE^r_s(M) &\coleq \cE[\fell](M, T^r_s M), &  (\cE_M)^r_s(M) &\coleq \cE_M[\fell](M, T^r_s M), \\
  \cN^r_s(M) &\coleq \cN[\fell](M, T^r_s M), &  \cG^r_s(M) &\coleq \cG[\fell](M, T^r_s M).
 \end{align*}
\end{definition}

Classically, an $(r,s)$-tensor field $t \in \cT^r_s M$ can be identified with a $C^\infty(M)$-multilinear map $\cT^0_1(M)^r \times \cT^1_0(M)^s \to C^\infty(M)$. The same holds true for generalized tensor fields:

\begin{theorem}
\begin{enumerate}[label=(\roman*)]
 \item \label{asdfasd} There is a canonical isomorphism
\begin{equation}\label{equation.4.35}
\cE^r_s(M) \cong \Hom_{\cE(M)} ( \cE^0_1(M)^r \times \cE^1_0(M)^s, \cE(M))
\end{equation}
given by
\begin{gather*}
R(\Theta^1, \dotsc, \Theta^r, X_1, \dotsc, X_s)(A, \Phi) = \\
R(A, \Phi) ( \Theta^1(A, \Phi), \dotsc, \Theta^r(A, \Phi), X_1(A, \Phi), \dotsc, X_s(A, \Phi)).
\end{gather*}
\item\label{asdfasd.2} For any open subset $V \subseteq M$, $R(\Theta_1, \dotsc, X_s)|_V = R|_V(\Theta^1|_V, \dotsc, X_s|_V)$.
\item\label{asdfasd.3} If $R$ and all $\Theta^i, X_j$ are moderate, so is $R(\Theta^1, \dotsc, X_s)$ and if one of them is negligible, so is $R(\Theta^1, \dotsc, X_s)$.
\item\label{asdfasd.4} The isomorphism \eqref{equation.4.35} induces an isomorphism
\begin{equation}\label{equation.4.36}
\cG^r_s(M) \cong \Hom_{\cG(M)} ( \cG^0_1(M)^r \times \cG^1_0(M)^s, \cG(M))
\end{equation}
given by 
\[ [R] ( [\Theta^1], \dotsc, [X_s] ) = [ R ( \Theta^1, \dotsc, X_s) ]. \]
\end{enumerate}
\end{theorem}
\begin{proof}
\ref{asdfasd} Denote by $\lambda_U$ for $U \subseteq M$ the canonical isomorphism
\[ \cT^r_s(U) \to \Hom_{C^\infty(U)} ( \cT^0_1(U)^r \times \cT^1_0(U)^s, C^\infty(U)). \]
We define $\widetilde \lambda_M \colon \cE^r_s(M) \to \Hom_{\cE(M)} ( \cE^0_1(M)^r \times \cE^1_0(M)^s, \cE(M))$ by
\begin{multline*}
(\widetilde \lambda_M R)(\Theta^1, \dotsc, \Theta^r, X_1, \dotsc, X_s)(A,\Phi) \coleq \\
\lambda_M ( R(A,\Phi))(\Theta^1(A,\Phi), \dotsc, \Theta^r(A, \Phi), X_1(A,\Phi), \dotsc, X_s(A, \Phi)).
\end{multline*}
It is easily verified that $\widetilde \lambda_M R$ indeed maps into $\cE(M)$ and is $\cE(M)$-multilinear. To see that $\widetilde \lambda_M$ is injective, suppose $\widetilde \lambda_M R = 0$. Then for all $A \in \TO(M)$, $\Phi \in \SO(M)$ and $\Theta^1, \dotsc, \Theta^r \in \cT^0_1(M)$, $X_1, \dotsc, X_s \in \cT^1_0(M)$ we have
\[ \lambda_M(R(A,\Phi))(\Theta^1, \dotsc, X_s) = (\widetilde \lambda_M R)(\Theta^1, \dotsc, X_s)(A, \Phi) = 0, \] which implies $\lambda_M ( R(A,\Phi)) = 0$ and hence $R(A,\Phi) = 0$. For surjectivity, assume $L$ in the range of $\widetilde \lambda_M$ is given; the inverse of $\lambda_M$ is then determined by
\[ \lambda_M(\widetilde \lambda_M^{-1}(L)(A, \Phi))(\Theta^1, \dotsc, X_s) \coleq L(\Theta^1, \dotsc, X_s)(A, \Phi). \]

\ref{asdfasd.2} To simplify notation we only show that $(\widetilde \lambda_M R)(S)|_V = (\widetilde \lambda_V R|_V)(S|_V)$.
We have, for suitable $A, \Phi, B, \Psi$ and $W$ as before,
\begin{align*}
(\widetilde \lambda_V R|_V)(S|_V)(A,\Phi)|_W &= \lambda_V ( R|_V(A,\Phi))(S|_V(A,\Phi))|_W \\
&= \lambda_W ( R|_V(A,\Phi)|_W)(S|_V(A,\Phi)|_W) \\
&= \lambda_W ( R(B,\Psi)|_W)(S(B,\Psi)|_W) \\
&= \lambda_U ( R(B,\Psi))(S(B,\Psi))|_W = (\widetilde \lambda_M R)(S)|_W.
\end{align*}
This holds for all $W$, so we are done.
\end{proof}

\begin{definition}For $X,Y \in \cE^1_0(M)$ we define the Lie Bracket $[X,Y] \in \cE^1_0(M)$ as
\[ [X,Y](A, \Phi) \coleq [X(A, \Phi), Y(A,\Phi)]. \]
 \end{definition}
We have $\widetilde \Lie_X \circ \widetilde \Lie_Y - \widetilde \Lie_Y \circ \widetilde \Lie_X = \widetilde \Lie_{[X,Y]}$ and $[X,Y] = \widetilde \Lie_X Y$ for all $X,Y \in \cE^1_0(M)$. The Lie bracket is $\bR$-bilinear, antisymmetric, satisfies the Jacobi identity and
\[ [FX, GY] = F G [X,Y] + F (\widetilde \Lie_X G) Y - G (\widetilde \Lie_Y F) X. \]

\subsection{Symmetric covariant generalized tensor fields}

A smooth tensor field $t \in \cT^0_s(M)$ is symmetric by definition if $\lambda_M t$ is a symmetric multilinear mapping. Similarly, we can define symmetry of generalized tensor fields on the level of representatives:

\begin{definition}\label{definition.4.4}
 An element $R \in \cE^0_s(M)$ is called \emph{symmetric} if the corresponding multilinear mapping $\widetilde \lambda_M(R)$ is symmetric; analogously, $\widehat R \in \cG^0_s(M)$ is called symmetric if the multilinear mapping given by \eqref{equation.4.36} is symmetric.
\end{definition}
By $\Sym$ we denote the symmetrizer
\[ (\Sym R)(X_1, \dotsc, X_s) \coleq \frac{1}{s!} \sum_{\sigma \in S_s} R(X_{\sigma(1)}, \dotsc, X_{\sigma(s)}) \]
acting either on $R \in \Hom_{C^\infty(U)} ( \cT^1_0(U) \times \dotsc \times \cT^1_0(U), C^\infty(U))$ or on $R \in \Hom_{\cE(U)} ( \cE^1_0(U) \times \dotsc \times \cE^1_0(U), \cE(U))$, where $S_s$ is the permutation group of $\{1, \dotsc, s\}$. The following result shows that $R$ is symmetric if and only if it has symmetric values.

\begin{lemma}\label{lemma.4.5}$R \in \cE^0_s(M)$ is symmetric if and only if $R(A, \Phi)$ is symmetric for all $A \in \TO(M)$ and $\Phi \in \SO(M)$.
\end{lemma}
\begin{proof}
Symmetry of $R$ is equivalent to $\widetilde \lambda_M R = \Sym ( \widetilde \lambda_M R)$, or in turn to
\[ R(A,\Phi) = \widetilde \lambda_M^{-1} ( \Sym ( \widetilde \lambda_M R)(A,\Phi)) \]
for all $A,\Phi$. The right hand side can be calculated as $\lambda_M^{-1}(\Sym(\lambda_M ( R(A,\Phi))))$ which proves the claim.
\end{proof}

One can then consider the submodules of moderate and negligible elements and form the quotient. Due to the inclusion $S^s T^*M \subseteq \tang^0_sM$ (where $S^k$ is the symmetric tensor functor), symmetric tensor fields fit in our functorial framework:

\begin{corollary}\label{corollary.4.6}
We have
\[ \{ R \in \cE^0_s(M)\ |\ R \textrm{ is symmetric} \} = \cE(M, S^s T^*M). \]
\end{corollary}

Henceforth, we will conveniently identify symmetric elements of $\cE^0_s(M)$ and elements of $\cE(M, S^s T^*M)$. Generalized symmetric covariant tensor fields thus are defined as follows:

\begin{definition}\label{definition.4.7}
For $s \in \bN_0$ we set
 \begin{align*}
  \cE^0_{s,\sym}(M) &\coleq \cE[\fell](M, S^s T^*M), &   (\cE_M)^0_{s,\sym}(M) &\coleq \cE_M[\fell](M, S^s T^*M), \\
  \cN^0_{s,\sym}(M) &\coleq \cN[\fell](M, S^s T^*M), &   \cG^0_{s,\sym}(M) &\coleq \cG[\fell](M, S^s T^*M).
 \end{align*}
\end{definition}

\begin{lemma}\label{lemma.4.8}
 $\cG^0_{s,\sym}(M) = \{ R \in \cG^0_s(M)\ |\ R \textrm{ is symmetric}\}$.
\end{lemma}
\begin{proof}
 The mapping $\cG^0_{s,\sym}(M) \to \cG^0_s(M)$ is clear, and also that its image consists of symmetric elements. Conversely, using the symmetrizer $\Sym \circ R$ (or rather $\lambda_M^{-1} \circ \Sym \circ \lambda_M \circ R$) one sees that
 \[ R \mapsto [(A,\Phi) \mapsto \Sym ( R (A,\Phi))] \]
 is a left inverse to $\cG^0_{s,\sym}(M) \to \cG^0_s(M)$.
\end{proof}

It follows that any symmetric element of $\cG^0_s(M)$ has a symmetric representative.

Recall that there also is an isomorphism
\[ \cD'^r_s(M) \cong \Hom_{C^\infty(M)} ( \cT^1_0(M) \times \dotsc \times \cT^1_0(M), \cD'(M)) \]
which allows one to speak of symmetric distributional tensor fields.
\begin{lemma}\label{lemma.4.9}
 If $u \in \cD'^0_s(M)$ is symmetric then $\iota u \in \cE^0_s(M)$ is symmetric as well. If $t \in \cT^0_s(M)$ is symmetric then $\sigma t \in \cE^r_s(M)$ is symmetric as well.
\end{lemma}
\begin{proof}
 For $u \in \cD'^r_s(M)$ we have, for any permutation $\sigma \in S_s$,
 \begin{align*}
  (\iota u)&(A,\varphi) \cdot (v_1 \otimes \dotsc \otimes v_s)(x) \\
  &= \langle u(y), (\lambda A)(x,y) \cdot (v_1 \otimes \dotsc \otimes v_s)(x) \otimes \vec\varphi(x)(y) \rangle \\
  &= \langle u(y), A^*(x,y) \cdot v_1(x) \otimes \dotsc \otimes A^*(x,y) \cdot v_j(x) \otimes \vec\varphi(x)(y) \rangle \\
  &= \langle u(y), A^*(x,y) \cdot v_{\sigma(1)}(x) \otimes \dotsc \otimes A^*(x,y) \cdot v_{\sigma(s)} (x) \otimes \vec\varphi(x)(y) \rangle \\
  &= (\iota u)(A, \varphi)(v_{\sigma(1)}, \dotsc, v_{\sigma(s)})(x).
 \end{align*}
The case of $\sigma t$ is clear.
\end{proof}

We will now examine how symmetry behaves with respect to restriction.

\begin{lemma}\label{lemma.4.10}If $R \in \cE^0_s(M)$ is symmetric then $R|_U \in \cE^r_s(U)$ is symmetric for all open subsets $U \subseteq M$.

Conversely, if $\widetilde R \in \cG^0_s(M)$ is such that for every point $x \in M$ there is a neighborhood $U$ such that $\widetilde R|_U$ is symmetric then $\widetilde R$ is symmetric.
\end{lemma}
\begin{proof}
For the first part we have to verify that $R|_U(A,\Phi)$ is symmetric for all $A,\Phi$. Because for smooth tensor fields we know symmetry is a local property it is enough to verify $R|_U(A,\Phi)|_V$ to be symmetric for compact subsets $V \csub U$; take $\rho \in \cD(U)$ which is 1 on $V$ and set $\vec\psi \coleq \rho \vec\varphi \in C^\infty(M, \Gamma(M, \Vol(M)))$. We know that there is a compact set $L \csub U$ such that $\vec\psi(x)$ has support in $L$ for all $x \in V$; let $\chi \in \cD(U)$ be 1 on $L$ and set $B \coleq (\rho \otimes \chi) A$. Then by construction of the restriction mapping we have
\[ R|_U ( A, \Phi)|_V = R(B, \Psi)|_V \]
and because $R$ is symmetric the claim follows.

For the second part we note that
\begin{multline*}
 \widetilde R(X_1, \dotsc, X_s)|_{U} = \widetilde R|_{U}(X_1|_{U}, \dotsc, X_s|_{U}) \\
 = \widetilde R|_{U}(X_{\sigma(1)}|_{U}, \dotsc, X_{\sigma(s)})|_{U} = \widetilde R(X_{\sigma(1)}, \dotsc, X_{\sigma(s)})|_U
\end{multline*}
for all $X_i \in \cG^1_0(M)$ and permutations $\sigma \in S_s$, which implies the claim.
\end{proof}

For the local setting, analogues of \Autoref{definition.4.4} to \Autoref{lemma.4.10} hold by the same arguments.

Finally, symmetry is preserved by the isomorphisms $\cE^r_s(U) \cong \cE^r_s(\varphi(U))$ and $\cG^r_s(U) \cong \cG^r_s(\varphi(U))$.

\section{Generalized semi-Riemannian geometry}\label{section.5}

\begin{definition}\label{definition.5.1}
A symmetric generalized tensor field $g \in \cG^0_{2,\sym}(M)$ is called \emph{nondegenerate} if the canonical homomorphism $\cG^1_0(M) \to \cG^0_1(M)$ it induces via \eqref{equation.4.36} is an isomorphism. Nondegenerateness of $g \in \cG^0_{2,\sym}(\Omega)$ is defined analogously.
\end{definition}
Note that a symmetric generalized tensor field $g \in \cG^0_{2,\sym}(U)$, where $(U,\varphi)$ is a chart, is nondegenerate if and only if its local expression in $\cG^0_{2,\sym}(\varphi(U))$ is nondegenerate.

\begin{definition}\label{definition.5.2}A \emph{generalized semi-Riemannian metric} on $M$ is a generalized tensor field $g \in \cG^0_2(M)$ which is symmetric and nondegenerate.
\end{definition}

\begin{lemma}\label{lemma.5.3}The homomorphism $\cG^1_0(M) \to \cG^0_1(M)$ induced by $g$ is an isomorphism if and only if it is so locally, i.e., if any $x \in M$ has an open neighborhood $U$ such that the homomorphism $\cG^1_0(U) \to \cG^0_1(U)$ induced by $g|_U$ is an isomorphism.
\end{lemma}
\begin{proof}
This is seen by the usual sheaf-theoretic methods and thus skipped here.
\end{proof}

Suppose $(U,\varphi)$ is a chart on $M$. As $\cG^1_0(U)$ and $\cG^0_1(U)$ are free $\cG(U)$-modules of the same dimension we can determine surjectivity of a module homomorphism $\cG^1_0(U) \to \cG^0_1(U)$ by considering the determinant of its matrix representation \cite[Prop.~4.18, p.~519]{LangAlgebra}. If this determinant is invertible we have an isomorphism and the inverse is given by the cofactor formula. Combining this with the characterization of invertibility in \Autoref{theorem.3.57} we obtain the following analytic criterion for nondegenerateness:

\begin{corollary}\label{corollary.5.4}Let $g \in \cG^0_{2,\sym}(M)$. Suppose we can cover $M$ by coordinate neighborhoods $U$ for which there is a representative $\tilde g$ of $g|_U$ that satisfies
\begin{equation}\label{equation.5.37}
\begin{aligned}
&\forall K\csub U\ \exists C_K>0\ \forall \bA \in \Upsilon(U)\ \forall \bPhi \in \S(U)\ \exists \e_0>0\ \forall \e <\e_0:\\
&\qquad \inf_{x \in K} (\det \tilde g_{ab})(\bA_\e, \bPhi_\e)(x) \ge C_K.
\end{aligned}
\end{equation}
 Then $g$ is a generalized semi-Riemannian metric.
\end{corollary}

The same holds in the local setting.

\subsection{Covariant derivatives}

\begin{definition}\label{definition.5.5}A \emph{generalized covariant derivative} on $M$ is defined to be a mapping $\nabla \colon \cE^1_0(M) \times \cE^1_0(M) \to \cE^1_0(M)$ such that for $X, Y, R, S \in \cE^1_0(M)$ and $f \in \cE(M)$,
 \begin{enumerate}[label=(\roman*)]
  \item $\nabla_{X+Y}R = \nabla_X R + \nabla_Y R$,
  \item $\nabla_{f X} R = f \nabla_X R$,
  \item $\nabla_X (R+S) = \nabla_X R + \nabla_X S$,
  \item $\nabla_X (fR) = (\widetilde\Lie_X f) R + f \nabla_X R$.
 \end{enumerate}
\end{definition}

$\nabla$ extends in a unique way to a derivation (i.e., a linear map satisfying the Leibniz rule) on $\bigoplus_{r,s\ge 0} \cE^r_s(M)$ such that $\nabla_X R = \widetilde \Lie_X R$ for $R \in \cE(M)$.  

\begin{lemma}A smooth covariant derivative $\nabla$ on $M$ extends to a generalized covariant derivative on $M$ by defining, for $X,R \in \cE^1_0(M)$,
\[ (\nabla_X R)(A,\Phi) \coleq \nabla_{X(A,\Phi)} ( R(A,\Phi) ) \qquad (A \in \TO(M,E), \Phi \in \SO(M)). \]
\end{lemma}
 
\begin{definition}
 A generalized covariant derivative $\nabla$ is called \emph{moderate} if for all $X,R \in (\cE_M)^1_0(M)$, $\nabla_X R \in (\cE_M)^1_0(M)$ .
\end{definition}

The following is easily seen:

\begin{proposition}
 \begin{enumerate}[label=(\roman*)]
  \item\label{modcd.1} Every smooth covariant derivative is moderate.
  \item\label{modcd.4} If $\nabla$ is moderate and $X,R \in (\cE_M)^1_0(M)$ then $\nabla_X R \in \cN^1_0(M)$ if $X$ or $R$ is negligible.
 \end{enumerate}
\end{proposition}

\begin{definition}\label{definition.5.9}Given a moderate generalized covariant derivative $\nabla$ on $M$, its action on $\cG^1_0(M)$ is defined as $\nabla_{[X]} [R] \coleq [\nabla_X R]$, where $[X] \in \cG^1_0(M)$ and $[R] \in \cG^1_0(M)$; this is independent of the representatives used.
\end{definition}

Hence, we also call a mapping $\nabla \colon \cG^1_0(M) \times \cG^1_0(M) \to \cG^1_0(M)$ which satisfies the identities of \Autoref{definition.5.5} a generalized covariant derivative.

\begin{theorem}\label{theorem.5.10}
 Let a manifold $M$ be given, endowed with a generalized semi-Riemannian metric $g \in \cG^0_2(M)$. There is a unique generalized covariant derivative $\nabla \colon \cG^1_0(M) \times \cG^1_0(M) \to \cG^1_0(M)$ on $M$ such that
\begin{itemize}
 \item $[X,Y] = \nabla_X Y - \nabla_Y X$
 \item $\widetilde \Lie_Z \langle X, Y \rangle = \langle \nabla_Z X, Y \rangle + \langle X, \nabla_Z Y \rangle$
\end{itemize}
for all $X,Y,Z \in (\cE^1_0)_M(M)$. $\nabla$ is called the Levi-Civit\`a derivative of $g$ and is characterized by the usual Koszul formula.
\end{theorem}
\begin{proof}
The proof is exactly as in the classical setting (cf.~\cite[Chapter 3, Theorem 11, p.~61]{oneill}).
\end{proof}

In a coordinate chart, the Levi-Civit\`a derivative can (as would be expected) be expressed by Christoffel symbols as in $\nabla_{\pd_i} \pd_j = \Gamma^k_{ij} \pd_k$ with
\[ \Gamma^i_{kl} = \frac{1}{2} g^{im} \left( g_{mk,l} + g_{ml, k} - g_{kl, m} \right). \]

\begin{definition}\label{definition.5.11}Let $\nabla \colon \cE^1_0(M) \times \cE^1_0(M) \times \cE^1_0(M)$ be a generalized covariant derivative on $M$.
The \emph{Riemannian curvature tensor} $\curv \in \cE[\fell](M, \Hom(\Lambda^2TM \otimes TM, TM))$ of 
$\nabla$ is defined by
\begin{equation}\label{equation.5.38}
\curv(X,Y)Z \coleq \nabla_X \nabla_Y Z - \nabla_Y\nabla_X Z + \nabla_{[X,Y]} Z \qquad (X,Y,Z \in \cE^1_0(M)).
\end{equation}
\end{definition}

We remark that if $\nabla$ is moderate then the curvature tensor is moderate, and in fact defines an element in $\cG(M, \Hom(\Lambda^2 TM \otimes TM, TM))$.

As $\Hom(\Lambda^2 TM \otimes TM, TM) \subseteq T^1_3 M$ we can regard $\curv$ as a generalized $(1,3)$-tensor field. If $\nabla$ is the Levi-Civit\`a derivative of a generalized semi-Riemannian metric then for its coordinates $\curv^l_{ijk}$ defined locally by
\[ \curv(\pd_i, \pd_j)\pd_k = \curv^l_{ijk} \pd_l \]
we have the usual formula
\[ \curv^l_{ijk} = \pd_{x_j} \Gamma^l_{ik} - \pd_k \Gamma^l_{ij} + \Gamma^l_{js} \Gamma^s_{ik} - \Gamma^l_{ks} \Gamma^s_{ij}. \]
and the lowered tensor $\curv_{lijk} = g_{ls}\curv^s_{ijk}$ is given by
\begin{equation}\label{curvaturecomponents}
\curv_{iklm} = \frac{1}{2} \left( g_{im,kl} + g_{kl,im} - g_{il, km} - g_{km, il} \right) + g_{np} ( \Gamma^n_{kl} \Gamma^p_{im} - \Gamma^n_{km} \Gamma^p_{il}).
\end{equation}

More generally, for a generalized covariant derivative on the quotient, i.e., a mapping $\cG^1_0(M) \times \cG^1_0(M) \to \cG^1_0(M)$, we define the curvature tensor also by \eqref{equation.5.38}.

\subsection{Metric association}

The concept of association as in \Autoref{definition.3.53} is sufficient for a purely analytic setting, but for our geometric setting we need a slightly different notion. We first recall that in the case of a smooth (semi-)Riemannian metric $g$ the volume density $\ud V_g$ determines an isomorphism between densities and functions, whence one can also view distributions as linear functionals on $C^\infty_c(M)$, the space of compactly supported smooth \emph{functions}; the volume density then is used for embedding regular distributions $f$ as in
\[ \langle f, \varphi \rangle \coleq \int f \varphi \,\ud V_g \qquad (\varphi \in C^\infty_c(M)). \]
If a manifold is endowed with a generalized metric $g$ then for the notion of association it is also reasonable to view $(\iota u)(\bA_\e, \bPhi_\e)$ for $u \in \cD'(M)$ as a linear functional on $C^\infty_c(M)$, but varying $\e$ changes the geometry of the underlying manifold. This leads to the concept of metric association:

\begin{definition}\label{definition.5.12}
For $R \in \cE(M)$ and $g \in \cE^0_{2,\sym}(M)$ we say that $R$ is \emph{metrically associated to $u \in (C_c^\infty(M))'$ with respect to $g$}, written as $R \approx_g u$, if for all $\bA \in \Upsilon(M)$, $\bPhi \in S(M)$ and $\omega \in C^\infty_c(M)$ we have
\[ \int R_\e\,\omega \,\ud V_{g_\e} \to \langle u, \omega \rangle \qquad (\e \to 0) \]
with $R_\e \coleq R(\bA_\e, \bPhi_\e)$ and $g_\e \coleq g(\bA_\e, \bPhi_\e)$.
\end{definition}

The following is easily seen in local coordinates.

\begin{lemma}
If $R_1 - R_2 \in \cN(M)$ and $g$ is moderate then $R_1 \approx_g u$ implies $R_2 \approx_g u$.

If $R$ is moderate and $g_1 - g_2 \in \cN^0_{2,\sym}(M)$ then $R \approx_{g_1} u$ implies $R \approx_{g_2} u$.
\end{lemma}

It follows that the notion of metric association is also well-defined for elements of the quotients, i.e., for $R \in \cG(M)$ and $g \in \cG^0_{2,\sym}(M)$.

It will be crucial for us to localize the concept of metric association as follows.

\begin{proposition}\label{proposition.5.14}
Let $R \in \cG(M)$ and $g \in \cG^0_{2,\sym}(M)$.
For any open cover $(U_i)_i$ of $M$, $R \approx_g u$ if and only if $R|_{U_i} \approx_{g|_{U_i}} u|_{U_i}$ for all $i$. 
\end{proposition}

Finally, this can be expressed in local charts.

\begin{proposition}\label{proposition.5.15}
 Let $(U, \varphi)$ be a chart on $M$, $R \in \cG(U)$ and $u \in \cD'(U)$. Then $R \approx_g \iota(u)$ if and only if $\varphi_* R \approx_{\varphi_*g} \varphi_* u$.
\end{proposition}

\section{The conical metric}\label{section.6}

As an application of the theory developed in the previous sections we will calculate the curvature of the conical metric \eqref{equation.1.1} as an associated distribution.
In Euclidean coordinates this metric takes the form
\begingroup
\renewcommand{\arraystretch}{2.5}
\begin{equation}\label{equation.6.39} g_{ab} = 
\frac{1}{2} (1+\alpha^2) \delta_{ab} + \frac{1}{2} (1-\alpha^2) h_{ab}
\end{equation}
\endgroup
where $\delta_{ab}$ is the identity matrix and
\[ h_{ab}(x_1,x_2) = \frac{1}{x_1^2+x_2^2}
 \left( \begin{matrix}
  x_1^2-x_2^2 & 2x_1y_1 \\
  2x_1x_2 & x_2^2-x_1^2
 \end{matrix} \right) \qquad ((x_1,y_1) \in \bR^2).
\]

This metric is smooth everywhere except at the origin and discontinuous there; moreover, it falls out of the class of gt-regular metrics introduced by Geroch and Traschen \cite{gerochtraschen}. This class gives conditions on the metric to guarantee that the curvature tensor exists as a distribution; for this one would require that $g$ is in $\Lint^\infty_{\loc} \cap H^1_{\loc}$, but it is only in $(\Lint^\infty_\loc \cap W^{1,1}_{\loc} ) \setminus H^1_{\loc}$ \cite[Lemma 1]{zbMATH06749306}.

A central feature of our theory is its global, diffeomorphism invariant nature. We assume the following situation: we are given a manifold $M$ with a chart $(U, \varphi)$ around a point $p_0$ such that $\varphi(U) = B_\mu(0)$ and $\varphi(p_0) = 0$ and $M$ is equipped with a singular metric tensor $G \in \cD'^0_{2,\sym}(M) \cap \cT^0_{2,\sym}(M \setminus \{p_0\})$. Locally on $U$ the metric $g \coleq \varphi_*(G|_U)$ is given by
\[ \ud s^2 = \ud r^2 + \alpha^2 r^2 \,\ud \phi^2 \]
with $\abso{\alpha} < 1$ (the case $\alpha=1$ corresponds to flat Euclidean space).
Let $\delta_{p_0} \in (C^\infty_c(M))'$ be the delta distribution at $p_0$, i.e.,
\[ \langle \delta_{p_0}, f \rangle = f(p_0) \qquad (f \in C^\infty_c(M)), \]
and $\hat G = [ \iota G] \in \cG^0_2(M)$ the conical metric as a generalized tensor field. Our main result then reads as follows.

\begin{theorem}\label{theorem.6.1}On $U$, the scalar curvature $\widehat R$ of the conical metric $\widehat G$ has $4\pi(1-\alpha) \cdot \delta_{p_0}$ as a metrically associated distribution:
\[ \widehat R \approx_{\widehat G} 4\pi(1-\alpha) \cdot \delta_{p_0} \qquad \textrm{on }U. \]
\end{theorem}

We shall not be concerned with $\widehat G$ outside of $U$ as it is smooth there.

The first task is to verify that $\widehat G|_U$ indeed is a generalized semi-Riemannian metric in the sense of \Autoref{definition.5.2}. Symmetry is clear by \Autoref{lemma.4.9} because $g$ is symmetric. For nondegeneracy we have to verify that the determinant $\det T_{\widehat G}$ of the mapping
\[ T_{\widehat G} \colon \cG^1_0(U) \to \cG^0_1(U) \]
associated to $\widehat G|_U$ via \eqref{equation.4.36} and the isomorphism
\[ \Hom_{\cG(U)} ( \cG^1_0(U) \times \cG^1_0(U), \cG(U)) \cong \Hom_{\cG(M)} ( \cG^1_0(M), \cG^0_1(M)) \]
is invertible. This is the case if $\varphi_*(\det T_{\widehat G}) = \det T_{\hat g}$ is invertible, where $T_{\hat g}$ is the mapping
\[ T_0 \colon \cG^1_0(\varphi(U)) \to \cG^0_1(\varphi(U)) \]
associated to
\begin{align*}
 \varphi_*(\widehat G|_U) &= \varphi_* ( [ \iota G]|_U) = \varphi_* ( [ ( \iota G)|_U]) = \varphi_* ( [ \iota(G|_U)]) \\
&= [ \varphi_*(\iota(G|_U))] = [ \iota ( \varphi_* ( G|_U))] = [\iota g] \eqcol \hat g.
\end{align*}

Talking of the determinant we will from now on identify $\hat g$ and $T_{\hat g}$ and just write $\det \hat g$.
By \Autoref{theorem.3.58} it suffices to find a representative in $\cE_M(\varphi(U))$ of $\det \hat g$ satisfying \eqref{equation.3.34}. For this, because for any representative $\tilde f$ of $\hat g$ $\det \tilde f$ is a representative of $\det \hat g$ we only have to choose $\tilde f$ appropriately. 

We will abbreviate $\Omega \coleq \varphi(\Omega)$ and $\Upsilon(\Omega) \coleq \Upsilon(\Omega, \bR^2)$. Moreover, we set $\tilde f_\e \coleq \tilde f(\bA_\e, \bPhi_\e)$ and similarly for $\tilde h_\e$, its coordinates $\tilde h_{ab\e}$ etc.
\begin{lemma}\label{lemma.6.2}
There exists a representative $\tilde f$ of $\hat g$ such that for any $K \csub \Omega$ and $0 < \kappa < \alpha^2$, $\bA \in \Upsilon(\Omega)$ and $\bPhi \in \SO(\Omega)$ there exists $\e_0>0$ such that $\kappa^2 \le \det \tilde f_\e \le (1 + \kappa)^2$ uniformly on $K$ for $\e<\e_0$.
\end{lemma}
\begin{proof}
A representative $\tilde f$ of $\hat g$ is given in Euclidean coordinates by
\[ \tilde f_{ab} = \frac{1+\alpha^2}{2} \delta_{ab} + \frac{1-\alpha^2}{2} \tilde h_{ab} \]
with
\[ \tilde h_{ab}(A, \vec\varphi)(x) = \int h_{ij}(y) A^i_a(y,x) A^j_b(y,x) \vec\varphi(x)(y) \,\ud y. \]
Here, we use the Einstein summation convention and the fact that for $\lambda \bE = \bE^*$, $(\lambda A)(x,y) = A^*(x,y) = A(y,x)^*$.
Note that $\tilde f_{ab}$ and $\tilde h_{ab}$ are symmetric.

To bound the eigenvalues of $\tilde f_{ab\e}$ and hence also its determinant we will use the well-known fact that for a real symmetric matrix $Q$ an upper or lower bound on $\sup \{ v^t \cdot Q \cdot v\ |\ \norm{v} = 1 \}$, where $v^t$ is the transpose of $v$, gives a corresponding bound on the eigenvalues of $Q$.

Via the characteristic polynomial the eigenvalues $\tilde \lambda_\e^+$ and $\tilde \lambda_\e^-$ of $\tilde h_{ab\e}$ are given by
\[ \tilde \lambda^\pm_\e = \frac{\tilde h_{11\e} + \tilde h_{22\e}}{2} \pm \sqrt{ \frac{ (\tilde h_{11\e} + \tilde h_{22\e})^2}{4} + \tilde h_{12\e}^2 - \tilde h_{11\e} \tilde h_{22\e} }. \]

Fix $K \csub \Omega$, $0 < \kappa < \alpha^2$, $\bA \in \Upsilon(\Omega)$ and $\bPhi \in S(\Omega)$. We claim that
\[ \limsup_{\e \to 0} \tilde \lambda^+_\e = 1,\qquad  \liminf_{\e \to 0} \tilde \lambda^-_\e = -1 \qquad\textrm{uniformly on }K. \]

Choose some $L \csub U$ with $K \csub L$ and take $C>0$ such that $\supp \bvarphi_\e(x) \subseteq B_{C\e}(x) \subseteq L$ for all $x \in K$ and small $\e$.
First, because $h_{11} = -h_{22}$ and $h_{12} = h_{21}$ we have
\begin{align*}
 \tilde h_{11\e}(x) & + \tilde h_{22\e}(x) = \int h_{ij}(y) ({\bA_\e}^i_1 {\bA_\e}^j_1 + {\bA_\e}^i_2 {\bA_\e}^j_2)(y,x)\bvarphi_\e(x)(y)\,\ud y \\
& = \int \bigl( h_{11}(y) ( {\bA_\e}^1_1 {\bA_\e}^1_1 + {\bA_\e}_2^1 {\bA_\e}_2^1 - {\bA_\e}_1^2 {\bA_\e}_1^2 - {\bA_\e}^2_2 {\bA_\e}^2_2)(y,x) \\
& \quad + 2 h_{12}(y) ( {\bA_\e}^1_1 {\bA_\e}^2_1 + {\bA_\e}^1_2 {\bA_\e}^2_2)(y,x) \bigr) \bvarphi_\e(x)(y)\,\ud y
\end{align*}
and hence
\begin{align*}
 | \tilde h_{11\e}(x) & + \tilde h_{22\e}(x) | \\
 & \le \norm{h_{11}}_{\infty} \cdot \sup_{y \in B_{C\e}(x)} \abso{ ({\bA_\e}^1_1 {\bA_\e}^1_1 - {\bA_\e}^2_2 {\bA_\e}^2_2)(y,x) } \cdot \norm{\bvarphi_\e(x)}_{L^1}\\
& + \norm{ h_{11}}_{\infty} \sup_{y \in B_{C\e(x)}} ( \abso{{\bA_\e}_2^1(y,x)}^2 + \abso{{\bA_\e}_1^2(y,x)}^2 ) \norm{\bvarphi_\e(x)}_{L^1} \\
& + 2 \norm{h_{12}}_\infty \cdot \sup_{y \in B_{C\e}(x)} \abso{{\bA_\e}^1_1(y,x)} \cdot \sup_{y \in B_{C\e}(x)} \abso{ {\bA_\e}^2_1(y,x)} \norm{ \bvarphi_\e(x)}_{L_1} \\
& + 2 \norm{h_{12}}_\infty \cdot \sup_{y \in B_{C\e}(x)} \abso{{\bA_\e}^1_2(y,x)} \cdot \sup_{y \in B_{C\e}(x)} \abso{ {\bA_\e}^2_2(y,x)} \norm{ \bvarphi_\e(x)}_{L_1}.
\end{align*}

By \Autoref{lemma.3.21} and boundedness of $\norm{\bvarphi_\e(x)}_{\Lint^1}$ (\Autoref{definition.3.23} \ref{def_scaltestobjloc.5}) these expressions converge to zero uniformly for $x \in K$. It hence remains to estimate the square root of
\[ \tilde h_{12\e}^2 - \tilde h_{11\e} \tilde h_{22\e} = \tilde h_{12\e}^2 + \tilde h_{11\e}^2 - \tilde h_{11\e} ( \tilde h_{22\e} + \tilde h_{11\e} ) \]
where the last term converges to zero by the above and by boundedness of $\tilde h_{11\e}$. We then have
\[ \sqrt{ \tilde h^2_{12\e}(x) + \tilde h^2_{11\e}(x) } \le \int \sqrt{h_{11}^2(y) + h_{12}^2(y) } \abso{ \bvarphi_\e(x)(y) } \,\ud y = \norm{\bvarphi_\e(x)}_{L^1} \to 1 \]
by \Autoref{definition.3.23} \ref{def_scaltestobjloc.6}, which gives the claim.

It follows that $\forall \delta > 0$ $\exists \e_0>0$ $\forall \e<\e_0$: $-(1+\delta) \le \tilde \lambda^-_\e \le \tilde \lambda^+_\e \le 1 + \delta$. We then need to have
\[ \frac{1+\alpha^2}{2} + \frac{1-\alpha^2}{2} \tilde\lambda^-_\e \ge \kappa, \qquad \frac{1+\alpha^2}{2} + \frac{1 - \alpha^2}{2} \tilde \lambda_\e^+ \le 1 + \kappa \]
for the statement of the theorem to hold, which is the case for $\e$ small enough.
\end{proof}

We hence have established that $\widehat G|_U$ is a generalized semi-Riemannian metric which by \Autoref{theorem.5.10} has an associated Levi-Civit\`a derivative whose curvature tensor $\widehat {\boldsymbol R} \in \cG^1_3(U)$ is given by \eqref{equation.5.38}. As in the classical situation, $\widehat {\boldsymbol R}$ is completely determined by $\widehat{\curv}_{1212} = \widehat G_{1s} \widehat{\curv}^s_{212}$ which is related to the scalar curvature
\[ \widehat R = \widehat G^{ij} \widehat G^{lm} \widehat {\boldsymbol R}_{iljm} \]
by
\[ \widehat R = \frac{2}{\det \widehat G} \widehat {\boldsymbol R}_{1212}. \]
We will now show that 
\[ \widehat R \approx_{\widehat G} 4\pi(1-A) \delta_{p_0}. \]
Although we have a nice expression for a representative of $\hat g$ to work with -- namely $g$ itself -- it will be necessary to have a representative $\widetilde R$ of the local expression $\varphi_*(\widehat R)$ of the scalar curvature $\widehat R$ such that $\widetilde R_\e$ is in fact the scalar curvature of the metric $\tilde g_\e$ on $\Omega$, at least for small $\e$. For this purpose we note that \Autoref{theorem.3.58} applies not only to $\det \tilde f$ but also to $\det \tilde g$, which means that there exists $\tilde s \in \cE_M(\Omega)$ such that
\begin{equation}\label{equation.6.40}
 \begin{aligned}
\forall K \csub \Omega\ \forall \bA \in \Upsilon(\Omega)\ \forall \bPhi \in S(\Omega)\ \exists \e_0>0 \\
\forall \e < \e_0\ \forall x \in K:\ \tilde s_\e(x) \cdot \det \tilde g_\e(x) = 1.
 \end{aligned}
\end{equation}
As a representative of the inverse metric of $\hat g$ we define $\tilde h \in (\cE_M)^2_0(\Omega)$ with the cofactor formula but using $\tilde s$ in place of $(\det \tilde g)^{-1}$:
\[ \tilde h^{ab} \coleq \tilde s \cdot C_{ba} \]
where $C_{ba}$ is the $(b,a)$-cofactor of $\tilde g$. Fixing any $K,\bA,\bPhi,\e_0$ and $\e<\e_0$ as in \eqref{equation.6.40} we then have
\begin{equation}\label{equation.6.41}
\tilde h^{ab}_\e = (\tilde g_\e)^{ab}
\end{equation}
on $K$, i.e., $\tilde h$ is an inverse of $\tilde g$ not only in the quotient but also on the level of representative for small $\e$.
However, $\tilde g_\e$ still is invertible only for small $\e$; let us denote its curvature tensor for such $\e$ by $k_\e$. This curvature tensor is a polynomial in the components of the metric and their derivatives of up to second order as well as the components of the inverse metric. If we replace in the corresponding formula for $\varphi_*(\widehat R)$ all occurrences of the inverse metric by the components of $\tilde h$ we also obtain an expression for a representative $\widetilde R$ thereof, given by the corresponding formula for the curvature of $\tilde g$ but with $\tilde h$ in place of the inverse metric. Now because of \eqref{equation.6.41} we know that given $K$, $\bA$ and $\bPhi$, $\widetilde R_\e = \widetilde R(\bA_\e, \bPhi_\e)$ equals the scalar curvature of $\tilde g_\e = \tilde g(\bA_\e, \bPhi_\e)$ on $K$ for small $\e$.

By \Autoref{proposition.5.15} we now have to show that
\[ \int \widetilde R_{\e}(x) \omega(x) \sqrt{\abso{\tilde g_\e(x)}}\,\ud x \to \langle 4\pi(1-A) \delta, \omega \rangle \qquad \forall \omega \in C^\infty_c(\Omega). \]
Fix $\omega$, $\bA \in \TO(\Omega)$ and $\bPhi \in \SO(\Omega)$ and let $0 < \lambda < \mu$ be such that $\supp \omega \subseteq B_\lambda(0)$. From the above discussion we can assume that for $x \in B_\lambda(0)$ and small $\e$,
$\widetilde R_{\e}$ in fact equals the scalar curvature of $\tilde g_\e$.

Following an idea of \cite{JVClarke}, we then write
\[ \int_{B_\lambda(0)} \widetilde R_\e(x) \omega(x) \sqrt{ \abso{ \tilde g_\e(x) } } \,\ud x = \omega(0,0) I_1 + I_2\]
where we expand $\omega$ at $0$ as in
\[ \omega(x) = \omega(0) + \int_0^1 (\D \omega)(tx) \cdot x \,\ud t \qquad (x \in B_\lambda(0)) \]
and the integrals $I_1$ and $I_2$ are given by
\begin{align*}
 I_1 &= \int_{B_\lambda(0)} \widetilde R_\e (x) \sqrt{\abso{\tilde g_\e(x)}}\,\ud x \\
 I_2 &= \int_{B_\lambda(0)} \int_0^1 \widetilde R_\e(x)  (\D \omega)(t x) \cdot x \sqrt{ \abso{ \tilde g_\e(x)}} \,\ud t \,\ud x.
\end{align*}

We can calculate $I_1$ by help of the Gauss-Bonnet theorem:
\begin{theorem}[{\cite[Theorem 9.3]{Lee}}]
 Let $\gamma$ be a smooth curve on an oriented Riemannian manifold $(M,g)$ of dimension 2 which is positively oriented as the boundary of an open set $\Omega$ with compact closure. Then
 \[ \int_\Omega K \,\ud V_g + \int_\gamma \kappa_g \,\ud s = 2\pi \]
 where $K$ is the Gaussian curvature of $g$, $\ud V_g$ its Riemannian volume element and $\kappa_g$ the geodesic curvature of $\gamma$.
\end{theorem}

Orienting $\Omega$ by its standard orientation induced from $\bR^2$ and noting that the scalar curvature $\widetilde R_\e$ is twice the Gaussian curvature $\widetilde K_\e$ of $\tilde g_\e$ we can calculate $I_1$ as
\[
 I_1 = 2 \int_{B_{\lambda}(0)} K_\e(x) \,\abso{\sqrt{\tilde g_\e(x)}}\,\ud x = 2 \left( 2\pi - \int_{\pd B_{\mu}(0)} \kappa_{\tilde g_\e}\,\ud s \right)
\]
where $\kappa_{\tilde g_\e}$ is the geodesic curvature of $\pd B_{\lambda}(0)$. As the metric is smooth in a neighborhood of this boundary we have, due to \Autoref{definition.3.23} \ref{def_scaltestobjloc.4} and continuity of all operations making up the integrand,
\[ \int_{\pd B_\lambda(0)} \kappa_{\tilde g_\e} \,\ud s \to \int_{\pd B_\lambda(0)} \kappa_g\,\ud s, \]
i.e., it suffices to integrate the non-regularized geodesic curvature. A short calculation then gives the value of this integral as $2\pi\alpha$, which results in
\[ I_1 \to 2(2\pi - 2 \pi \alpha) =  4\pi(1-\alpha). \]

For the second integral our approach is different and more general than the one from \cite{JVClarke}. We will employ the following result, which describes how a function which is smooth outside a singularity is approximated there.

\begin{lemma}\label{lemma.6.4}Let $f \in C^\infty(\Omega \setminus \{0\}) \cap \Lint^1_{\loc}(\Omega)$, $\mu_1,\mu_2>0$ such that $\mu_1 + \mu_2 \le \mu$, $\bPhi \in S(\Omega)$ and take $C>0$, $\e_0>0$ such that $C \e_0 < \min(\mu_2, \mu_1/2)$ and $\supp \varphi_\e(x) \subseteq B_{C\e}(x)$ $\forall \e<\e_0$ $\forall x \in B_{\mu_1}(0)$. Then for $\alpha \in \bN_0^n$ and $q \in \bN$ with $q > \abso{\alpha}$ there is $L>0$ such that for $\e<\e_0$ and $2 C \e \le \abso{x} \le \mu_1$, with
\[ \tilde f_\e(x) \coleq \int f(y) \bvarphi_\e(x)(y)\,\ud y \]
we have
\[ \abso{ (\pd^\alpha \tilde f_\e)(x) - (\pd^\alpha f)(x) } \le L \e^{q - \abso{\alpha}} \sup_{\substack{\abso{\beta} \le q \\ y \in B_{C\e}(x)}} \abso{ \pd^\beta f(y) }. \]
\end{lemma}
\begin{proof}
By Taylor's theorem we have for $2C\e \le \abso{x} \le \mu_1$
\begin{gather*}
 (\pd^\alpha \tilde f_\e)(x) - (\pd^\alpha f)(x) = \int f(y) (\pd_x^\alpha \bvarphi_\e)(x)(y) \,\ud y - (\pd^\alpha f)(x) \\
= \sum_{\abso{\beta} < q} (\pd^\beta f)(x)  \int \frac{(y-x)^\beta}{\beta!} (\pd_x^\alpha \bvarphi_\e)(x)(y)\,\ud y - (\pd^\alpha f)(x) \\
+ \sum_{\abso{\beta} = q} \frac{q}{\beta!} \int (y-x)^\beta \int_0^1 (1-t)^{q-1} \cdot \\
\qquad (\pd^\beta f)(x+t(y-x))\,\ud t \, (\pd_x^\alpha \bvarphi_\e)(x)(y)\,\ud y.
\end{gather*}
First, we note that
\[ \int \frac{(y-x)^\beta}{\beta!} (\pd_x^\alpha \bvarphi_\e)(x)(y)\,\ud y - \delta_{\beta\alpha} = O(\e^{q-\abso{\alpha}}) \]
uniformly for $\abso{x} \le \mu_1$ by \Autoref{definition.3.23} \ref{def_scaltestobjloc.4}. Because $q > \abso{\alpha}$ this means that
\[ \sum_{\abso{\beta} < q} (\pd^\beta f)(x)  \int \frac{(y-x)^\beta}{\beta!} (\pd_x^\alpha \bvarphi_\e)(x)(y)\,\ud y - (\pd^\alpha f)(x) = O(\e^{q-\abso{\alpha}}) \cdot \sup_{\abso{\beta} < q} \abso{\pd^\beta f(x)} \]
The remainder terms are estimated as in
\begin{gather*}
 \abso{ \int_{B_{C\e}(x)} (y-x)^\beta \int_0^1 (1-t)^{q-1} (\pd^\beta f)(x + t(y-x)) \,\ud t\, (\pd_x^\alpha \bvarphi_\e)(x)(y)\,\ud y } \\
\le C' \e^q \sup_{y \in B_{C\e}(x)} \abso{ \pd^\beta f(y) } \sup_{\abso{x} \le \mu_1} \norm{ \pd_x^\alpha \bvarphi_\e(x) }_1 \\
\le C'' \e^{q-\abso{\alpha}} \sup_{y \in B_{C\e}(x)} \abso{\pd^\beta f(y)}.
\end{gather*}
which proves the claim.
\end{proof}
A similar result holds if we incorporate transport operators:

\begin{lemma}\label{lemma.6.5}
Let $f \in C^\infty(\Omega \setminus\{0\}, S^2 \bR^2) \cap L^1_{\loc}(\Omega, S^2 \bR^2)$ and $\mu_1,\mu_2>0$ such that $\mu_1 + \mu_2 \le \mu$. Let $\bA \in \TO(\Omega)$, $\bPhi \in S(\Omega)$ and take $C>0$, $\e_0>0$ such that $C\e_0 < \min(\mu_2, \mu_1/2)$ and $\supp \varphi_\e(x) \subseteq B_{C\e}(x)$ $\forall \e<\e_0$ $\forall x \in B_{\mu_1}(0)$. Then for $\alpha \in \bN_0^n$ and $q \in \bN$ with $q > \abso{\alpha}$ there is $L>0$ such that for $\e<\e_0$ and $2 C \e \le \abso{x} \le \mu_1$ we have
\[ \abso{ (\pd^\alpha \tilde f_{ij\e})(x) - (\pd^\alpha f_{ij})(x) } \le L \e^{q - \abso{\alpha}} \sup_{\substack{\abso{\beta} \le q \\ y \in B_{C\e}(x) \\ k,l}} \abso{ \pd^\beta f_{kl}(y) }. \]
\end{lemma}
The proof is similar to that of \Autoref{lemma.6.4} (albeit a bit more technical due to the presence of transport operators) and thus omitted. To calculate $I_2$ we split the integral into two regions. Fix $C$ such that
\[ \supp \bvarphi_\e(x) \subseteq B_{C \e}(x) \]
for $x \in B_\lambda(0)$ and small $\e$. With
\[ M = \sup_{x \in B_\lambda(0)} \norm{(\D \omega)(x)} < \infty \]
we have
\begin{equation}\label{equation.6.42}
 \begin{aligned}
 \abso{I_2} \le \quad & M \int_{\abso{x} < 2 C \e} \abso{ \widetilde R_\e(x)}  \sqrt{ \abso{\tilde g_\e(x)}} \abso{x}\,\ud x\\
+ & M \int_{2 C \e < \abso{x} < \lambda } \abso{ \widetilde R_\e(x)}  \sqrt{ \abso{\tilde g_\e(x)}} \abso{x} \,\ud x.
 \end{aligned}
\end{equation}

For the first integral in \eqref{equation.6.42} we have to be able to bound $\widetilde R_\e$ by $O(\e^{-2})$. From the coordinate formula for $\widetilde R_{1212\e}$ (see \eqref{curvaturecomponents}) we see that it consists of a sum of either second derivatives of the metric or the products of two first derivatives, eventually multiplied by the metric which is bounded. Without transport operators, the relevant estimate is obtained from
\begin{gather*}
 \abso{ \pd^\alpha \tilde f_\e(x)} = \abso{ \int f(y) \pd_x^\alpha \bvarphi_\e(x)(y) } \le \norm{f}_{\infty} \norm{\pd_x^\alpha \bvarphi_\e(x)}_{L^1} = O(\e^{-\abso{\alpha}}).
\end{gather*}
With transport operators, some of the derivatives fall on $\bA_\e$ instead of $\bvarphi_\e$, which does not add any growth due to \Autoref{definition.3.7} \ref{definition.3.7.1} and we still have the same estimate: $\abso{\pd^\alpha \tilde f_{ab\e} (x) }$ can be estimated by terms of the form (with $\alpha_1 + \alpha_2 + \alpha_3 = \alpha$)
\begin{gather*}
\abso{ \int f_{ij}(y) \pd_x^{\alpha_1} ( {\bA_\e}^i_a(y,x)) \pd_x^{\alpha_2} ( {\bA_\e}^j_b(y,x)) \pd_x^{\alpha_3} \bvarphi_\e(x)(y) } \\
\le \sup_{i,j} \norm{f_{ij}}_{\infty} \cdot \norm{ \pd_x^{\alpha_1} ( {\bA_\e}^i_a(.,x))}_{\infty} \norm{\pd_x^{\alpha_2} ( {\bA_\e}^j_b)(.,x))}_{\infty} \norm{\pd^{\alpha_3}_x \bvarphi_\e(x)}_{L^1} \\
= O(\e^{-\abso{\alpha_3}}) = O(\e^{-\abso{\alpha}}).
\end{gather*}
Hence, $\widetilde R_\e = O(\e^{-2})$ holds in a neighborhood of the origin and the first integral in \eqref{equation.6.42} can be estimated by
\[ \e^{-2} \int_{\abso{x} < 2C \e} \abso{x} \, \ud x = \e^{-2} \int_0^{2 C\e} r^2 \,\ud r \to 0 \qquad (\e \to 0). \]

For the second integral in \eqref{equation.6.42} we claim that for some $q \in \bN$ we have
\begin{equation}\label{equation.6.43}
\abso{\tilde R_\e(x)} = O\left(\frac{\e^q}{\abso{x}^{q+2}}\right) \qquad\textrm{uniformly on }2C\e < \abso{x} < \mu.
\end{equation}

For $\abso{x} \ne 0$ we know that the curvature $R$ of $g$ vanishes at $x$; thus, we can write (suppressing $\e$ and $x$ in the notation from now on)
\[ \widetilde R = \frac{2}{\det \tilde g} (\tilde R_{1212} - R_{1212}). \]
As $\det \tilde g$ and its inverse are bounded by \Autoref{lemma.6.2} we can safely ignore it for the estimate. We rewrite the difference as
\begin{multline*}
\widetilde R_{iklm} - R_{iklm} = \\
\frac{1}{2} \left( \tilde g_{im,kl} + \tilde g_{kl,im} - \tilde g_{il,km} - \tilde g_{km,il} \right) + \tilde g_{np} ( \tilde \Gamma^n_{kl} \tilde \Gamma^p_{im} - \tilde \Gamma^n_{km} \tilde \Gamma^p_{il}) \\
- \frac{1}{2} \left( g_{im,kl} + g_{kl,im} - g_{il,km} - g_{km,il} \right) - g_{np} ( \Gamma^n_{kl} \Gamma^p_{im} - \Gamma^n_{km} \Gamma^p_{il})
\end{multline*}
The first terms to estimate here are of the form $\tilde g_{im,kl} - g_{im,kl}$. By \Autoref{lemma.6.5} we have
\[ \abso{\tilde g_{im,kl} - g_{im,kl}} \le L \e^{q} \sup_{\substack{\abso{\beta} \le q+2\\y \in B_{C\e}(x)\\k,l}} \abso{\pd^\beta g_{kl}(y)} \]
and as $g_{kl}$ is homogeneous of degree $0$ this is $O(\e^q \abso{x}^{-(q+2)})$ as desired because for $2C\e \le \abso{x} \le \lambda$ and $\abso{y-x} \le C\e$ we have $\abso{y} \ge \abso{x}/2$ and hence $\abso{\pd^\beta g_{kl}(y)} \sim \abso{y}^{-\abso{\beta}} \lesssim \abso{x}^{-\abso{\beta}} \le \abso{x}^{-(q+2)}$.
The other terms are of the form
\[ \tilde g_{np} \tilde \Gamma^n_{kl} \tilde \Gamma^p_{im} - g_{np} \Gamma^n_{kl} \Gamma^p_{im}. \]
We first write
\[ g_{np} \Gamma^n_{kl} \Gamma^p_{im} = \frac{1}{4} g_{np}g^{ns}g^{pr}(g_{sk,l} + g_{sl,k} - g_{kl,s})(g_{ri,m} + g_{rm,i} - g_{im,r}) \]
and similarly for the regularized version; expanding the parentheses we hence end up with terms of the form
\[ g_{ab}g^{cd}g^{ef}g_{gh,i}g_{jk,l}. \]
The expressions to estimate then are obtained from
\begin{align*}
\nonumber \tilde g_{ab} \tilde  g^{cd} & \tilde g^{ef}\tilde g_{gh,i}\tilde g_{jk,l}
 - g_{ab}g^{cd}g^{ef}g_{gh,i}g_{jk,l} = \\
& ( \tilde g_{ab} - g_{ab} ) \tilde  g^{cd}\tilde g^{ef}\tilde g_{gh,i}\tilde g_{jk,l} \\
+\, & g_{ab} (\tilde g^{cd} - g^{cd}) \tilde g^{ef}\tilde g_{gh,i}\tilde g_{jk,l} \\
+\, & g_{ab} g^{cd} (\tilde g^{ef}-g^{ef}) \tilde g_{gh,i}\tilde g_{jk,l} \\
+\, & g_{ab} g^{cd} g^{ef} (\tilde g_{gh,i} - g_{gh,i}) \tilde g_{jk,l} \\
+\, & g_{ab} g^{cd} g^{ef} g_{gh,i} (\tilde g_{jk,l} - g_{jk,l})
\end{align*}
The relevant estimates (where $q$ can be chosen at will for each term) are as follows:
\begin{enumerate}
 \item $\tilde g_{ab} - g_{ab} = O(\e^q \abso{x}^{-q})$ by \Autoref{lemma.6.5}.
 \item $\tilde g^{ab} - g^{ab} = O(\e^q \abso{x}^{-q})$ via the cofactor formula. In fact, $g^{ab} = (\det g)^{-1} C_{ba}$ where the $(b,a)$-cofactor $C_{ba}$ of $g$ is a polynomial in the components of $g$, and similarly $\tilde g^{ab} = (\det \tilde g)^{-1} \tilde C_{ab}$. this 
 \[ \tilde g^{ab} - g^{ab} = \left( \frac{1}{\det \tilde g} - \frac{1}{\det g} \right) C_{ba} + \frac{1}{\det g} \left( C_{ba} - \tilde C_{ba} \right). \]
 While $C_{ba}$ is bounded, $\det \tilde g - \det g = O(\e^q\abso{x}^{-q}$ by \Autoref{lemma.6.5} and hence
 \[ \frac{1}{\det \tilde g} - \frac{1}{\det g} = \frac{\det g - \det \tilde g}{\det \tilde g \cdot \det g} = O(\e^q \abso{x}^{-q}) \]
 by boundedness of $\det \tilde g$ and $\det g$. Similarly, one has
 \[ \frac{1}{\det g} ( C_{ba} - \tilde C_{ba} ) = O(\e^q \abso{x}^{-q}). \]
 \item $\tilde g_{ab,c} - g_{ab,c} = O(\e^q \abso{x}^{-(q+1)})$ by \Autoref{lemma.6.5}.
 \item $\tilde g^{cd} = O(1)$.
 \item $\tilde g_{gh,i} = O(\e^{-1})$.
 \item $g_{cd} = O(1)$.
 \item $\tilde g_{ab} = O(1)$.
 \item $g_{gh,i} = O(\abso{x}^{-1})$.
\end{enumerate}
Putting these estimates together (and choosing $q$ appropriately for the individual terms) we have established \eqref{equation.6.43}. Hence, the second integral of \eqref{equation.6.42} can be estimated in polar coordinates by
\[ \int_{2 \e C}^\lambda \e^q r^{-(q+2)} r^2 \,\ud r = \e^q \int_{2 \e C}^\lambda r^{-q} \,\ud r = \frac{\e^q}{q-1} \left( \frac{1}{(2\e C)^{q-1}} - \frac{1}{\lambda^{q-1}} \right) = O(\e). \]
This completes the proof of \Autoref{theorem.6.1}.

\subsection*{Acknowledgments}

This research was supported by grants P26859 and P30233 of the Austrian Science Fund (FWF).

\printbibliography

\end{document}